\documentclass[nonblindrev]{informs3_hide}

\OneAndAHalfSpacedXI 

\usepackage{mathabx, makecell}
\usepackage{natbib, multirow, multicol, xspace, enumitem, amsmath, pifont, algorithm, algpseudocode, subcaption}
\usepackage{accents}
\usepackage{algorithmicx}
\usepackage{graphicx}
\usepackage{multicol}
\usepackage{mdframed}
\usepackage{lmodern}
\usepackage{multirow}
\usepackage{makecell}
\bibpunct[, ]{(}{)}{,}{a}{}{,}%
\def\bibfont{\small}%
\usepackage{bm}
\usepackage[normalem]{ulem}
\usepackage[dvipsnames]{xcolor}
\usepackage[all]{xy}
\usepackage{hyperref}
\hypersetup{
    colorlinks=true,
    linkcolor=blue,
    filecolor=magenta,      
    urlcolor=cyan,
    citecolor=blue
    }
\usepackage{ulem}
\usepackage{framed}

\newcommand{\bp}{\bm{p}}
\newcommand{\bq}{\bm{q}}
\newcommand{\bu}{\bm{u}}
\newcommand{\Reg}{\mathrm{Regret}}

\newcommand{\defeq}{\overset{\text{\tiny def}}{=}}
\theoremstyle{TH}%

\TheoremsNumberedThrough     
\ECRepeatTheorems

\EquationsNumberedThrough    

\MANUSCRIPTNO{}

\begin{document}

\RUNAUTHOR{Chen et al.}

\RUNTITLE{Online Linear Programming with Replenishment}

\TITLE{Online Linear Programming with Replenishment}

\ARTICLEAUTHORS{
\AUTHOR{Yuze Chen\footnotemark[1]}
 \AFF{Qiuzhen College, Tsinghua University, Beijing 100084, China}
 \AUTHOR{Yuan Zhou\footnotemark[1]\footnotemark[2]}
 \AFF{Yau Mathematical Sciences Center \& Department of Mathematical Sciences, Tsinghua University, Beijing 100084, China, \\ Beijing Institute of Mathematical Sciences and Applications, Beijing 101408, China}
 \AUTHOR{Baichuan Mo}
 \AFF{Bytedance Inc., 1199 Coleman Ave, San Jose, CA 95110, USA}
 \AUTHOR{Jie Ying}
 \AFF{Qiuzhen College, Tsinghua University, Beijing 100084, China}
 \AUTHOR{Yufei Ruan}
 \AFF{Bytedance Inc., 1199 Coleman Ave, San Jose, CA 95110, USA}
 \AUTHOR{Zhou Ye}
 \AFF{Bytedance China, Aili Center, Cangqian Street, Yuhang District, Hangzhou, Zhejiang 310000, China}
}

\renewcommand{\thefootnote}{\fnsymbol{footnote}}
\footnotetext[1]{These authors contributed equally.}
\footnotetext[2]{Corresponding author. Email: \EMAIL{yuan-zhou@tsinghua.edu.cn}.}
\renewcommand{\thefootnote}{\arabic{footnote}}

\ABSTRACT{We study an online linear programming (OLP) model in which inventory is not provided upfront but instead arrives gradually through an exogenous stochastic replenishment process. This replenishment-based formulation captures operational settings, such as e-commerce fulfillment, perishable supply chains, and renewable-powered systems, where resources are accumulated gradually and initial inventories are small or zero. The introduction of dispersed, uncertain replenishment fundamentally alters the structure of classical OLPs, creating persistent stockout risk and eliminating advance knowledge of the total budget. 

We develop new algorithms and regret analyses for three major distributional regimes studied in the OLP literature: bounded distributions, finite-support distributions, and continuous-support distributions with a non-degeneracy condition. For bounded distributions, we design an algorithm that achieves $\widetilde{\mathcal{O}}(\sqrt{T})$ regret. For finite-support distributions with a non-degenerate induced LP, we obtain $\mathcal{O}(\log T)$ regret, and we establish an $\Omega(\sqrt{T})$ lower bound for degenerate instances, demonstrating a sharp separation from the classical setting where $\mathcal{O}(1)$ regret is achievable. For continuous-support, non-degenerate distributions, we develop a two-stage accumulate-then-convert algorithm that achieves $\mathcal{O}(\log^2 T)$ regret, comparable to the $\mathcal{O}(\log T)$ regret in classical OLPs. Together, these results provide a near-complete characterization of the optimal regret achievable in OLP with replenishment. Finally, we empirically evaluate our algorithms and demonstrate their advantages over natural adaptations of classical OLP methods in the replenishment setting.
}


\maketitle

\section{Introduction}
In recent years, online linear programming (OLP) has emerged as a central framework for modeling a wide range of online decision-making problems with resource constraints. An OLP describes a sequential process in which each arriving order consists of a vector of requested resources together with an associated reward. Upon observing an order, the decision maker must immediately and irrevocably determine whether to accept it, earning the reward provided sufficient inventory is available, or to reject it. Each order can be viewed as a column of the constraint matrix of an underlying linear program (LP), with the corresponding decision variable $x_t \in \{0, 1\}$ indicating acceptance or rejection. Since these columns arrive online and decisions must be made without knowledge of future arrivals, the problem is naturally referred to as an \emph{online linear program}. This influential abstraction captures numerous classical online optimization problems, including the online knapsack \citep{babaioff2007knapsack}, online advertising and ad allocation \citep{mehta2007adwords,devanur2009adwords,mahdian2011online}, online packing \citep{buchbinder2009online,cetin2025online}, network revenue management \citep{jasin2015performance}, the multi-secretary problem \citep{kleinberg2005multiple,arlotto2019uniformly}, online auctions \citep{buchbinder2007online}, and e-commerce fulfillment \citep{jasin2015lp}.

Starting from the seminal work of \citet{agrawal2014dynamic}, most prior work on OLP assumes a large initial resource endowment that scales with the time horizon $T$, and all online orders draw exclusively from this static initial inventory throughout the selling season.\footnote{Some works, such as \citet{li2022online}, also allow inventory to be replenished during the selling season through orders that require a negative amount of resources. However, their algorithm still fundamentally relies on a large initial inventory endowment.} In many operational environments, however, resources are generated and accumulated through uncertain, small, and dispersed replenishment events rather than deterministically concentrated in a single large initial stock. For example, e-commerce fulfillment centers face variable and unpredictable upstream supply; fresh-food and grocery items are harvested or delivered daily in uncertain quantities; and renewable-powered computing systems receive intermittently arriving energy. Because inventory capacity, such as warehouse storage limits or battery capacity, is typically constrained, these systems seldom begin a season with the full amount of resources needed over the entire horizon. In settings such as fresh-food and grocery distribution, resource perishability further renders large upfront stocking impractical, making gradual and uncertain replenishment a more accurate representation of operational reality.

In this paper, motivated by the above considerations, we investigate an OLP model in which inventory is not provided upfront but instead arrives gradually through an exogenous stochastic replenishment process. Starting from a zero initial inventory, the system receives a random replenishment at the beginning of each period; the decision maker then observes the arriving order and must immediately decide whether to accept it, earning the reward and consuming the required resources, or to reject it. The objective is to maximize the cumulative reward over the planning horizon.

The introduction of exogenous replenishment significantly alters the structure of the OLP. A central challenge is that resource scarcity may arise at any point in the process if too many orders are accepted before sufficient inventory has been replenished. This is in contrast to the classical setting, where shortages typically occur only near the end, as the large initial budget is exhausted. To address this challenge, our key algorithmic insight is to manage early-stage scarcity by strategically rejecting even profitable orders, thereby preserving an inventory buffer for potentially higher-value orders that may arrive later. This \emph{accumulate-then-utilize} principle is intuitive at a high level, but implementing it in a regret-optimal way requires new algorithmic ideas and analytical techniques that differ across the distributional settings we study. For bounded distributions, where the target is $\tilde{\mathcal{O}}(\sqrt{T})$ regret, we realize this principle by directly rejecting a fixed number of roughly $\tilde{\mathcal{O}}(\sqrt{T})$ early orders to build an inventory buffer large enough to accommodate the random fluctuations of the remaining $\Theta(T)$ stochastic orders while avoiding costly stockout events. In the other settings we study, where the regret targets are logarithmic, this straightforward approach is no longer viable, as rejecting that many early orders would itself incur more regret than the desired bound.

An additional difference between the classical OLP and our replenishment-based setting lies in the availability of information about the total budget. In the classical OLP, the total amount of inventory (i.e., the budget) is known upfront, enabling direct normalization and facilitating many existing algorithmic techniques. In our setting, however, the total resource budget is dispersed across multiple stochastic replenishment events, whose cumulative supply is unknown in advance, introducing an additional layer of uncertainty. To cope with this uncertainty, we develop new online learning mechanisms that estimate relevant quantities on the fly and support near-optimal decision making throughout the horizon.

In this work, we systematically characterize the optimal regret achievable in OLP with replenishment under several natural classes of order and replenishment distributions. These classes include general \emph{bounded distributions}, \emph{finite-support distributions}, and \emph{continuous-support and non-degenerate distributions}, each serving as a direct analogue of distributional families extensively studied in the classical OLP literature. For all these settings, we establish regret bounds that are optimal up to logarithmic factors. In many cases, despite the additional challenges created by gradual and uncertain replenishment, we develop new algorithms that (almost) match the optimal regret known for the classical OLP problem under the same distributional assumptions. In other cases, however, we prove regret lower bounds that are substantially larger than the upper bounds achieved in the classical setting, thereby demonstrating a clear separation between the two problems and highlighting the intrinsic difficulty introduced by stochastic and dispersed replenishment. Below, we report our results and technical contributions.

\subsection{Theoretical Results and Technical Contributions}

Let $\mathcal{P}$ denote the underlying order and replenishment distribution. We study the OLP with replenishment problem when $\mathcal{P}$ belongs to the following three distributional classes. A summary of our results, along with a comparison to the classical OLP setting, is provided in Table~\ref{tab:result-summary}.

\renewcommand{\arraystretch}{1.9}
\begin{table}[]
    \centering
    \resizebox{\columnwidth}{!}{%
    \begin{tabular}{c|c||c|c|c}
        \hline
         \multicolumn{2}{c||}{} & \multicolumn{2}{c|}{OLP with Replenishment} & Classical OLP \\ \hline\hline 
         \multirow{2}{*}[-4pt]{\makecell{(1) Bounded\\ Distributions}} & UB & \multicolumn{2}{c|}{$\widetilde{\mathcal{O}}(\sqrt{T})$ \footnotesize{(Algorithm~\ref{alg:bounded-distribution})}} & $\widetilde{\mathcal{O}}(\sqrt{T})$ \footnotesize{ \citep{agrawal2014dynamic,li2023simple}} \\ \cline{2-5} 
         & LB &  \multicolumn{2}{c|}{$\Omega(\sqrt{T})$ \footnotesize{(implied by Theorem~\ref{thm:finite-dis-degenerate-hard-instance})}} &  \makecell{$\Omega(\sqrt{T})$ \\ {\footnotesize \citep{agrawal2014dynamic,arlotto2019uniformly}}}  \\ \hline \hline
         \multirow{2}{*}[-6pt]{\makecell{(2) Finite-Support\\ Distributions}} & UB & \makecell{{\footnotesize Non-Deg.~Induced LP:} \\ $\mathcal{O}(\log T)$ \footnotesize{(Algorithm~\ref{alg:finite-support-distribution})}} & \makecell{{\footnotesize Degenerate Induced LP:} \\ $\widetilde{\mathcal{O}}(\sqrt{T})$ \footnotesize{(Algorithm~\ref{alg:bounded-distribution})}} & \makecell{$\mathcal{O}(1)$ \\ \footnotesize{\citep{xie2025benefits,li2024infrequent}}} \\ \cline{2-5}
         & LB & \makecell{{\footnotesize Non-Deg.~Induced LP:} \\ $\Omega(1)$ \footnotesize{(trivial)}} & \makecell{{\footnotesize Degenerate Induced LP:} \\ $\Omega(\sqrt{T})$ \footnotesize{(Theorem~\ref{thm:finite-dis-degenerate-hard-instance})}} & $\Omega(1)$ \footnotesize{(trivial)} \\ \hline \hline
         \multirow{2}{*}[-1pt]{\makecell{(3) Continuous-Support\\ \& Non-Degenerate\\ Distributions}} & UB & \multicolumn{2}{c|}{$\mathcal{O}(\log^2 T)$ \footnotesize{(Algorithm~\ref{alg:replenish-non-degenerate-dis-cts-combine})}} &  \makecell{$\mathcal{O}(\log T)$ \\\footnotesize{ \citep{li2022online,bray2025logarithmic,ma2025optimal}} }\\ \cline{2-5}
         & LB & \multicolumn{2}{c|}{$\Omega(1)$ \footnotesize{(trivial)} } & $\Omega(\log T)$ \footnotesize{\citep{bray2025logarithmic}} \\ \hline
    \end{tabular}
    }
    \caption{Summary of our theoretical results, and comparison to the classical OLP setting. (2) and (3) are special cases of (1). UB denotes ``upper bound'' and LB denotes ``lower bound''.}
    \label{tab:result-summary}
\end{table}
\renewcommand{\arraystretch}{1}

\medskip
\noindent\underline {(1) Bounded Distributions.} In Section~\ref{sec:bounded-distribution}, we assume that  $\mathcal{P}$ is a bounded distribution (Definition~\ref{def:bounded-dis}), which is the most general distributional class considered in the OLP literature. The best-known regret guarantee for classical OLPs under this class is $\widetilde{\mathcal{O}}(\sqrt{T})$ \citep{agrawal2014dynamic,arlotto2019uniformly}. In our replenishment-based setting, we propose Algorithm~\ref{alg:bounded-distribution}, which also achieves $\widetilde{\mathcal{O}}(\sqrt{T})$ regret. This bound is tight up to logarithmic factors via matching lower bounds presented later in this paper. As previously discussed, the algorithm follows an accumulate-then-utilize principle by introducing a \emph{warm-up phase} in which the first $\Theta(\sqrt{T \log T})$ orders are automatically rejected to build an initial inventory buffer. After this warm-up period, the algorithm transitions to the \emph{main decision phase}, which employs online gradient descent to estimate dual resource prices; an order is accepted only if its reward exceeds its estimated dual cost. 

As also noted earlier, another complication relative to classical OLPs is that the total resource budget in our replenishment-based setting is no longer fixed and known in advance, but instead depends on an unknown stochastic replenishment process. This makes direct gradient-based dual-price estimation (as in, e.g., \citet{li2023simple}) infeasible, since the gradient depends on the total budget. To overcome this issue, we replace the unknown budget with realized replenishment observations, which introduces additional noise into the gradient estimator. Nevertheless, our analysis shows that this noise can be controlled, allowing us to establish the desired $\widetilde{\mathcal{O}}(\sqrt{T})$ regret guarantee.

\medskip
\noindent\underline{(2) Finite-Support Distributions.} In Section~\ref{sec:finite-support-distribution}, we introduce Algorithm~\ref{alg:finite-support-distribution}. This algorithm operates under the assumption that $\mathcal{P}$ is bounded and also has finite support (Definition~\ref{def:finite-support-distribution}, i.e., there are only finitely many order types). In additions, the algorithm requires the induced linear program (LP) to be non-degenerate (Definitions~\ref{def:induced-lp} and \ref{def:finite-support-non-degenerate}). Under these conditions, we prove that the algorithm achieves $\mathcal{O}(\log T)$ regret. Similar to the bounded-distribution case, the algorithm follows an accumulate-then-utilize principle, but now rejects only the first $\Theta(\log T)$ orders, building a much smaller inventory buffer, thanks to the additional finite-support and non-degenerate induced-LP conditions. The remaining time periods are partitioned into short batches of length $\Theta(\log T)$. At the beginning of each batch, the algorithm estimates an induced LP using both the current inventory level and the samples collected in the preceding batch, and then constructs the decision policy for the current batch from the solution of this estimated LP.

In the classical OLP setting, \citet{xie2025benefits,li2024infrequent} design algorithms that achieve $\mathcal{O}(1)$ regret for every finite-support distribution. Our result for the replenishment-based setting is comparable up to logarithmic factors, but requires the additional non-degeneracy assumption on the induced LP. To complement our upper bound, we construct a hard instance with finite support but a degenerate induced LP, and show that any algorithm must incur $\Omega(\sqrt{T})$ regret on this instance. This lower bound underscores the critical role of induced-LP non-degeneracy in determining the optimal regret for OLP with replenishment under finite-support distributions, and further demonstrates a clear separation from the classical OLP setting, highlighting the fundamental difficulty introduced by gradual, stochastic replenishment. Moreover, since finite-support distributions form a special case of general bounded distributions, our lower bound implies that the $\widetilde{\mathcal{O}}(\sqrt{T})$ regret achieved by Algorithm~\ref{alg:bounded-distribution} is minimax-optimal up to logarithmic factors.

\medskip
\noindent\underline{(3) Continuous-Support and Non-Degenerate Distributions.} In Section~\ref{sec:non-degenerate-dis}, we assume that the unknown distribution $\mathcal{P}$ is bounded, has continuous support, and additionally satisfies the non-degeneracy condition (Definition~\ref{def:cts-non-degeneracy}). Since the non-degeneracy condition automatically implies that $\mathcal{P}$ is bounded and has continuous support, we also refer to this case as the non-degenerate distribution for simplicity. Under this assumption, we design Algorithm~\ref{alg:replenish-non-degenerate-dis-cts-combine} and show that it achieves $\mathcal{O}(\log^2 T)$ regret. For comparison, \citet{li2022online} develop an algorithm that attains $\mathcal{O}(\log T)$ regret for classical OLPs under continuous-support and non-degenerate distributions. Our regret bound therefore matches the classical benchmark up to logarithmic factors.

Our algorithmic design, however, differs significantly from that of \citet{li2022online} because we do not have access to a large initial inventory. It is also the most technically involved algorithm in this paper. To achieve the desired logarithmic regret, our algorithm needs to accumulate a reserve of $\Omega(\sqrt{T})$ inventory. However, unlike the previous algorithms, it cannot simply reject the first $\Omega(\sqrt{T})$ orders to build this reserve, as doing so would itself incur much more regret than our target bound. Instead, the algorithm follows the accumulate-then-utilize principle in a more nuanced manner:  it employs a time-varying cost-effectiveness threshold that evolves over the horizon, rejecting orders whose cost-efficiency falls below this dynamic threshold. More specifically, the algorithm splits the horizon into two stages of roughly equal length. During the first stage, it runs an \emph{inventory accumulation procedure} that mildly suppresses resource consumption by inflating the estimated dual prices, thereby allowing the system to gradually build a reserve of $\Omega(\sqrt{T})$ inventory. Both the suppression of the resource budget and the inflation of the dual-price estimates are dynamically adjusted with the time index $t$, each decaying at roughly the rate $\mathcal{O}(\sqrt{(\log T)/ t})$. In the second stage, the algorithm transitions to an \emph{inventory conversion procedure} that aims to convert the accumulated inventory, together with ongoing replenishment, into realized rewards. This is accomplished by slightly enlarging the resource-consumption budgets. The amount of enlargement is carefully tuned according to the remaining horizon, scaling at roughly $\mathcal{O}(\sqrt{(\log T)/ (T-t)})$. This dynamic adjustment allows the algorithm to effectively leverage the inventory accumulated in the first stage, avoid costly stockout events, and fully utilize the available inventory to realize rewards.

To analyze the algorithm and establish its $\mathcal{O}(\log^2 T)$ regret, we consider a different performance metric, the \emph{dual-valued gain}, which accounts for both the immediate reward earned and the value of the remaining inventory evaluated at the optimal dual prices. We show that the regret with respect to the dual-valued gain is $\mathcal{O}(\log^2 T)$ in both stages of the algorithm. This, in turn, implies an $\mathcal{O}(\log^2 T)$ regret bound for the cumulative reward, since most of the inventory is consumed by the end of the horizon.

\medskip
\noindent{\bf Additional Considerations: Initial Inventory, Capacity Constraints, and Holding Time.} All of our algorithmic results extend naturally to settings with both non-zero initial inventory and per-period stochastic replenishment, without affecting the regret guarantees. Moreover, our algorithms for the general bounded-distribution case and for non-degenerate distributions can be adapted, by discarding excess inventory, to operate under an inventory capacity of only $\widetilde{\mathcal{O}}(\sqrt{T})$, while maintaining their existing regret bounds. This requirement is much smaller than the $\Theta(T)$ capacity needed in classical OLPs. For finite-support distributions, the required capacity can be reduced even further, to $\mathcal{O}(\log T)$. Similarly, the maximum holding time can be bounded by $\widetilde{\mathcal{O}}(\sqrt{T})$ for the general bounded and non-degenerate distribution classes, and by $\mathcal{O}(\log T)$ for finite-support distributions. The details and analysis of these generalizations are conceptually straightforward and are omitted due to space constraints. Nevertheless, these considerations highlight the flexibility of our algorithmic frameworks and their ability to accommodate a wide range of operational constraints.

\section{Literature Review}

\underline{OLP with Unknown Distributions.} The online linear programming (OLP) problem with an unknown underlying distribution was first introduced by \citet{agrawal2014dynamic}, whose results apply to both the \emph{i.i.d.} arrival model and the random permutation model. In what follows, we concentrate on prior work under the \emph{i.i.d.} arrival model, which is also the setting of our paper. For completeness, we note that there is also a line of research on the random permutation model (e.g., \cite{gupta2014experts,molinaro2014geometry}).

For \emph{i.i.d.} orders, the original work \citet{agrawal2014dynamic} propose an algorithm to achieve $\widetilde{\mathcal{O}}(\sqrt{T})$ regret for the general bounded distributions. Much of the subsequent literature distinguishes between two classes of underlying distributions: those with finite support and those with continuous support. In the finite-support setting, the order distribution is supported on finitely many bounded types, and many results rely on assuming that the induced linear program is non-degenerate. For example, \citet{jasin2015performance} establishes an $\mathcal{O}(\log^{2} T)$ regret bound, which is later improved to $\mathcal{O}(1)$ by \citet{wei2023constant} and \citet{chen2024improved}. More recently, \citet{xie2025benefits} and \citet{li2024infrequent} remove the non-degeneracy assumption entirely and propose algorithms that achieve $\mathcal{O}(1)$ regret for general finite-support distributions.
In the continuous-support setting, the order distribution is assumed to lie on a bounded continuous domain with a positive and bounded density. Under the additional assumption of non-degeneracy, \citet{li2022online} propose a dual-price-based LP algorithm that achieves regret $\mathcal{O}(\log T \log\log T)$. Subsequent refinements by \citet{bray2025logarithmic} and \citet{ma2025optimal} show that this algorithm in fact achieves $\mathcal{O}(\log T)$ regret, and \citet{bray2025logarithmic} further establishes a matching lower bound of $\Omega(\log T)$ for this setting. 

All the aforementioned algorithmic results rely on repeatedly solving linear programs, which can be computationally expensive. In the finite-support setting, \citet{li2024infrequent} design LP-based algorithms that attain $\mathcal{O}(1)$ regret while requiring only $\mathcal{O}(\log\log T)$ LP solves, thereby reducing computational cost. A parallel line of work also investigates LP-free approaches based on subgradient descent. \citet{balseiro2020dual,balseiro2022analysis,li2023simple,gao2023solving} propose algorithms that achieve $\mathcal{O}(\sqrt{T})$ regret for general bounded distributions. More recently, sharper regret guarantees have been obtained for continuous-support distributions under additional structural assumptions---for example, $\mathcal{O}(T^{1/3})$ regret under a quadratic-growth condition \citep{gao2024decoupling}, and $\mathcal{O}(\log^{2} T)$ regret under the non-degeneracy assumption \citep{ma2025optimal}.

In this work, we study OLP with replenishment across the major distributional classes examined in the OLP literature: bounded distributions, finite-support distributions, and continuous-support distributions satisfying the non-degeneracy condition. Without knowing the distribution, our regret bounds match the best-known guarantees for classical OLP in many of these settings; in others, we establish lower bounds showing that the replenishment-based formulation is inherently more challenging than the classical model. In addition, our algorithm for bounded distributions is LP-free, while our algorithm for continuous-support, non-degenerate distributions requires a modest number of LP solves (on the order of $\mathcal{O}(\log T)$).

\medskip
\noindent\underline{OLP with Known Distributions.} Another important stream of literature studies the setting where the underlying order distribution is known to the decision-maker. This setting is closely connected to the classical Network Revenue Management (NRM) problem. Research on NRM has advanced substantially over the past three decades, originating from the airline industry \citep{talluri1998analysis} and expanding to a wide variety of capacity‐constrained service and platform applications. A large body of work focuses on static policies, which compute a fluid LP solution at the outset and implement a time‐stationary control rule, where prominent examples include overbooking control \citep{beckmann1958decision,rothstein1971airline}, bid-price control \citep{talluri1998analysis}, booking-limit control \citep{cooper2002asymptotic}, probabilistic allocation control \citep{reiman2008asymptotically}, and price-discrimination control \citep{littlewood1972forecasting}. Motivated by the potential benefits of dynamic re-optimization, subsequent work examines LP-based resolving policies, which periodically (or infrequently) re-solve an LP based on the known arrival distribution and the remaining resources. Under the non-degeneracy assumption for the induced LP in the finite-support case, \citet{jasin2012re,jasin2013analysis} establish $\mathcal{O}(1)$ regret guarantees using such resolving methods. Later, \citet{bumpensanti2020re} obtain $\mathcal{O}(1)$ regret without requiring the non-degeneracy condition. For arbitrary bounded distributions (including both finite- and continuous-support cases), \citet{chen2025beyond} show that the certainty equivalent heuristic achieves $\mathcal{O}(\log^2 T)$ regret under the (reverse) $H\Ddot{o}lder$ condition.

Although all of our algorithms for OLP with replenishment operate without access to the underlying distribution, our $\Omega(\sqrt{T})$ lower bound for the finite-support case with a non-degenerate induced LP continues to hold even when the algorithm fully knows the distribution. This stronger result highlights the substantial additional challenges introduced by gradual and stochastic replenishment, even under complete distributional information.

\medskip
\noindent\underline{OLP with replenishment.} Most existing works on OLP assume a large initial inventory (typically scaling linearly with the horizon~$T$), and do not incorporate replenishment during the process. In some works~\citep{li2022online,li2023simple,bumpensanti2020re,gao2024decoupling}, replenishment is introduced implicitly by allowing orders with negative resource consumption, so that accepting such requests increases the available inventory. This mechanism results in endogenous replenishment driven by the arrival stream. Nevertheless, these works still require a large initial inventory.

\citet{vera2024dynamic} consider a more general framework that incorporates exogenous replenishment and does not require any initial inventory. However, their approach assumes full knowledge of the underlying distribution and focuses only on finite-support distributions. In addition, both resource consumption and replenishment amounts are restricted to take values in $\{0,1\}$.

In contrast, we study OLP with exogenous replenishment in a setting where the system may start with zero initial inventory. Our algorithms also accommodate endogenous replenishment through orders with negative resource consumption, while operating without access to the underlying distribution.

\section{Problem Formulation and Preliminaries}
\label{sec:problem-formulation}
We study an online resource allocation problem in which a firm sequentially processes orders involving $m$ resources indexed by $\{1,2,\dots,m\}$. The initial inventory of each resource is zero. At each discrete period $t$, each resource $j$ is first replenished by an exogenous amount $b_{jt} \in \mathbb{R}_{\geq 0}$. Then an order arrives with reward $r_t \in \mathbb{R}$ and requirement vector $\bm{a}_t=(a_{1t},\ldots,a_{mt}) \in \mathbb{R}^m$, where $a_{jt}\geq 0$ denotes a demand for $a_{jt}$ units of resource $j$ and $a_{jt}<0$ denotes a restock of $-a_{jt}$ units of resource $j$. The vector $\bm{a}_t$ can encode a pure purchase ($a_{jt}\geq 0$ for all $j$), a pure sale ($a_{jt}\leq 0$ for all $j$), or a mixed transaction. Upon observing the order, the decision-maker must immediately and irrevocably accept ($x_t=1$) or reject ($x_t=0$). If accepted, the firm earns $r_t$ and inventory updates by subtraction of the net requirements; otherwise, neither reward nor inventory changes. In other words, the inventory evolves as
\[
\ell_{j,t} \leftarrow \ell_{j,t-1} + b_{jt} -  a_{jt} x_t, \qquad \qquad \forall j \in [m],
\]
subject to the no-short-selling constraint $\ell_{j,t}\geq 0$, i.e., the inventory level of  resource $j$ at time $t$ should be nonnegative for all $j,t$. The objective is to maximize cumulative reward $\sum_{t=1}^{T} r_t x_t$ over the horizon $T$. 

We formulate the problem as the following integer linear program, denoted OLP-R (online linear programming with replenishment):
\begin{align}
    &\text{maximize} & &\sum_{t=1}^T r_t x_t \label{eq:basic-offline-LP} \\
    &\text{subject to} & &\sum_{k=1}^t a_{jk} x_k\leq \sum_{k=1}^t b_{jk}, \quad &\forall j\in[m],~t\in[T],\label{eq:basic-offline-LP-nonnegative-inventory} \\
    && &  x_t \in \{0, 1\},\quad &\forall t\in[T],  \notag
\end{align}
where $r_t\in\mathbb{R}$, $\bm{a}_t=(a_{1t},a_{2t},\dots,a_{mt})\in\mathbb{R}^m$, and $\bm{b}_t=(b_{1t},b_{2t},\dots,b_{mt})\in\mathbb{R}_{\geq 0}^m$, and the prefix constraints \eqref{eq:basic-offline-LP-nonnegative-inventory} encode nonnegativity of inventories for each resource $j$ and time $t$.

Suppose the entire sequence $(r_t, \bm{a}_t, \bm{b}_t)_{t=1}^{T}$ is known beforehand. Solving the integer linear program \eqref{eq:basic-offline-LP} yields the \emph{hindsight optimal solution} $\bm{x}^* = (x_1^*, x_2^*, \dots, x_T^*)$, and the \emph{hindsight optimal reward} is $R^* = \sum_{t=1}^T r_t x_t^*$. 

In the online setting, we assume that the period-$t$ tuple $(r_t, \bm{a}_t, \bm{b}_t)$ is drawn \emph{i.i.d.}~from an unknown distribution $\mathcal{P}$. The decision-maker must decide whether to accept the order at time $t$ based only on the \emph{history} $\mathcal{H}_t = ((r_k, \bm{a}_k, \bm{b}_k, x_k)_{k=1}^{t-1}, (r_t, \bm{a}_t, \bm{b}_{t}))$. A (possibly randomized) online policy is a vector $\bm{\pi}=(\pi_1,\dots,\pi_T)$ with 
\begin{align}
    \pi_{t}:\left[\prod_{k=1}^{t-1} \left(\mathbb{R}\times\mathbb{R}^m\times\mathbb{R}^m_{\geq 0}\times\mathbb{R}\right)\right]\times \left(\mathbb{R}\times\mathbb{R}^m\times\mathbb{R}^m_{\geq 0}\right)\to [0, 1] ,
\end{align}
so that, given $\mathcal{H}_t$, the acceptance decision is drawn as $x_t \sim \mathrm{Bernoulli}(\pi_t(\mathcal{H}_t))$.  The online policy $\bm{\pi}$ must also be \emph{feasible} throughout this paper: for any sequence $(r_t, \bm{a}_t, \bm{b}_t)_{t=1}^{T}$, the induced decisions $\bm{x} = (x_1, \dots, x_T)$ satisfy the nonnegativity of inventory constraints \eqref{eq:basic-offline-LP-nonnegative-inventory}. The reward achieved by the policy is $R^{\bm{\pi}} = \sum_{t=1}^T r_t x_t$. For any distribution $\mathcal{P}$, the \emph{regret} of an online policy $\bm{\pi}$ is defined as 
\begin{align}
    \Reg_{\mathcal{P}}(\bm{\pi})\defeq \mathbb{E}_{\mathcal{P}}\left[R^*-R^{\bm{\pi}}\right],
\end{align}

\subsection{Linear Programming Relaxation and its Dual}
The following linear program relaxes the integer linear program~\eqref{eq:basic-offline-LP} by dropping the integrality of decision variables and enforcing only the terminal (cumulative) inventory constraints, thereby yielding an upper bound on the hindsight optimum $R^*$:
\begin{align}
    &\text{maximize} &&\sum_{t=1}^T r_t x_t \label{eq:basic-offline-relaxation-LP} \\
    &\text{subject to} &&\sum_{t=1}^T a_{jt} x_t\leq \sum_{t=1}^T b_{jt}, \quad &\forall j\in[m], \label{eq:basic-offline-relaxation-LP-nonnegativity-inventory} \\
    &&& 0\leq x_t\leq 1,\quad &\forall t\in[T].  \notag
\end{align}
This LP relaxation is in a similar form to the standard relaxation for the classical OLP in \citet{agrawal2014dynamic}: there, capacities are given by initial inventories (typically scaling with $T$) and no exogenous replenishments, so the right-hand side of the constraint~\eqref{eq:basic-offline-relaxation-LP-nonnegativity-inventory} is replaced by the initial inventory level.

We will also need to study the dual of the LP relaxation~\eqref{eq:basic-offline-relaxation-LP}, which is expressed below.
\begin{align}
    &\text{minimize} & & \sum_{j=1}^m p_j\cdot\left(\sum_{t=1}^T b_{jt}\right)+ \sum_{t=1}^T s_t \label{eq:basic-offline-relax-dual} \\
    &\text{subject to} & &s_t +\sum_{j=1}^m a_{jt}p_j\geq r_t,\qquad &\forall t\in[T], \notag\\
    &&&p_{j},s_t\geq 0,\qquad&\forall j\in[m],~t\in[T].\notag
\end{align}
If we denote the optimal solution to Eq.~\eqref{eq:basic-offline-relaxation-LP} by $\bm{x}^*$ and the optimal dual solution to Eq.~\eqref{eq:basic-offline-relax-dual} by $(\bp^*, \bm{s}^*)$, then by complementary slackness, we have
\begin{equation}
    x_t^*=
    \begin{cases}
        1\qquad \text{if }r_t > \langle \bm{a}_t, \bp^*\rangle \\
        0\qquad \text{if }r_t <  \langle\bm{a}_t, \bp^*\rangle
    \end{cases},
    \label{eq:complemantary-condition}
\end{equation}
for $t\in[T]$. When $r_t=\langle \bm{a}_t, \bp^*\rangle$, the optimal solution $x_j^*$ may be a non-integer value. The implication of Eq.~\eqref{eq:complemantary-condition} is that the decision variable $x_t$ can be largely determined by the dual solution $\bp^*_T$, which is also referred to as the \emph{dual price} in the literature. 

\subsection{Bounded Distributions}

In this paper, we assume $\mathcal{P}$ lies in a broad class of \emph{bounded distributions} defined below.
\begin{definition}[Bounded Distribution]
\label{def:bounded-dis}
Let $\mathcal{P}$ be a distribution on $\mathbb{R}\times\mathbb{R}^m\times \mathbb{R}^m_{\geq 0}$ and $(r,\bm{a},\bm{b})$ be a random sample from $\mathcal{P}$ such that $r\in\mathbb{R}$, $\bm{a}\in\mathbb{R}^m$, and $\bm{b}\in\mathbb{R}^m_{\geq 0}$. We say $\mathcal{P}$ is \emph{bounded} if:
\begin{enumerate}[label=(\alph*)]
    \item There exist constants $\bar{r}, \bar{a}, \bar{b} > 0$ such that $|r| < \bar{r}$, $\|\bm{a}\|_{2} < \bar{a}$, and $\|\bm{b}\|_{\infty} < \bar{b}$ hold almost surely. 
    \item There exists $\underline{b} > 0$ such that $\mathbb{E}\left[b_{j}\right] > \underline{b}$ for all $j\in[m]$.
\end{enumerate}
\end{definition}
The boundedness assumption on $\mathcal{P}$ is standard in many influential works on the classical OLP problem (see, e.g., \citet{li2022online, li2023simple, gao2023solving, gao2024decoupling, ma2025optimal}).

\subsection{Finite-Support Distributions and the Induced Linear Program}

\emph{Finite-support distributions} form an important subclass of bounded distributions widely studied in the classical OLP literature \citep{bumpensanti2020re,chen2024improved, li2024infrequent, xie2025benefits,banerjee2025good}. We extend this notion naturally to the setting of OLP with replenishment as follows.
\begin{definition}[Finite-Support Distribution]
\label{def:finite-support-distribution}
Let $\mathcal{P}$ be a distribution on $\mathbb{R}\times\mathbb{R}^m\times \mathbb{R}^m_{\geq 0}$ and $(r,\bm{a},\bm{b})$ be a random sample from $\mathcal{P}$. We say  $\mathcal{P}$ has \emph{finite support} if:
\begin{enumerate}[label=(\alph*)]
    \item There exist $n$ support pairs $(R_i,\bm{A}_i)$ ($i\in[n]$) with $R_i\in\mathbb{R}$ and $\bm{A}_i=(A_{1i},A_{2i},\dots,A_{mi})\in\mathbb{R}^m$, such that $\mu_i\defeq \Pr\left((r,\bm{a})=(R_i,\bm{A}_i)\right)\geq\underline{\mu}>0$ and $\sum_{i=1}^n \mu_i=1$.
    \item There exist constants $\underline{b},\bar{b} > 0$ such that $B_j \defeq \mathbb{E}\left[b_{j}\right]>\underline{b}$ and $\Pr(b_{j}<\bar{b})=1$ for all $j\in[m]$.
\end{enumerate}
\end{definition}
If $\mathcal{P}$ is finite-support, it can be verified that $\mathcal{P}$ is also bounded in the sense of Definition~\ref{def:bounded-dis} with parameters $\bar{r} > \max_{i \in [n]} |R_i|$, $\bar{a} > \max_{i \in [n]} \|\bm{A}_i\|_2$,  and the same $\bar{b}$, $\underline{b}$ from above. 

Now, consider a realization of the sequence $\{(r_t, \bm{a}_t, \bm{b}_t)\}_{t=1}^T$ drawn \emph{i.i.d.}~from $\mathcal{P}$. Let $\theta_t \in [n]$ denote the type index of the order at time $t$, i.e., $(r_t, \bm{a}_t) = (R_{\theta_t}, A_{\theta_t})$. We can then reformulate the LP relaxation~\eqref{eq:basic-offline-relaxation-LP} into the following equivalent form, where $X_i$ represents the normalized total amount of accepted type-$i$ orders, i.e., $X_i = T^{-1} \cdot \sum_{t=1}^T \mathbb{I}(\theta_t = i) \cdot x_t$:
\begin{align}
    &\text{maximize} &&T\cdot\sum_{i=1}^n R_i X_i \label{eq:finite-support-rewritten-lp} \\
    &\text{subject to} &&T\cdot\sum_{i=1}^n A_{ji}X_i \leq \sum_{t=1}^T b_{jt},\qquad&\forall j\in[m], \notag \\
    &&&0\leq X_i\leq \frac{1}{T}\cdot\sum_{t=1}^T \mathbb{I}\left( \theta_t=i \right),&\forall i\in[n].\notag
\end{align}
\begin{definition}[Induced Linear Program]
\label{def:induced-lp}
Denote $\bm{B} = (B_j)_{j \in [m]}$, $\bm{\mu} = (\mu_i)_{i \in [n]}$, where $B_j$ and $\mu_i$ are defined in Definition~\ref{def:finite-support-distribution}. By replacing the stochastic constraints in \eqref{eq:finite-support-rewritten-lp} with their expectations, we obtain the following \emph{induced linear program}, denoted $\mathrm{LP}^{\mathrm{induced}}(\bm{B}, \bm{\mu})$, which serves as a fluid approximation of the LP relaxation under a finite-support distribution:
\begin{align}
    &\text{maximize} &&\sum_{i=1}^n R_i X_i \label{eq:finite-support-fluid-lp} \\
    &\text{subject to} &&\sum_{i=1}^n A_{ji}X_i +S_j = B_j,\qquad&\forall j\in[m], \notag \\
    &&&X_i+V_i= \mu_i,&\forall i\in[n],\notag \\
    &&&X_i,V_i,S_j\geq 0,\qquad&\forall i\in[n],j\in[m]. \notag
\end{align}
\end{definition}

The following proposition about the induced LP can be easily verified.
\begin{proposition}\label{fact:induced-lp-scale}
Let $\{X^*_i, V^*_i\}_{i\in[n]}, \{S_j^*\}_{j\in[m]}$ be an optimal solution to $\mathrm{LP}^{\mathrm{induced}}(\bm{B}, \bm{\mu})$. Then, for any $\alpha > 0$, $\{\alpha X^*_i, \alpha V^*_i\}_{i\in[n]}, \{\alpha S_j^*\}_{j\in[m]}$ is an optimal solution to $\mathrm{LP}^{\mathrm{induced}}(\alpha \bm{B}, \alpha \bm{\mu})$, and $\alpha \cdot \mathrm{LP}^{\mathrm{induced}}(\bm{B}, \bm{\mu}) = \mathrm{LP}^{\mathrm{induced}}(\alpha\bm{B}, \alpha \bm{\mu})$.
\end{proposition}

This following notion of non-degeneracy is standard in the linear programming literature. In particular, any degenerate LP can be perturbed by an arbitrarily small amount to become non-degenerate; see~\citet{megiddo1989varepsilon}.
\begin{definition}[non-degeneracy of induced linear program and optimal basis]
\label{def:finite-support-non-degenerate}
We say that LP~\eqref{eq:finite-support-fluid-lp} is \emph{non-degenerate} if:
\begin{enumerate}[label=(\alph*)]
\item The LP has a unique optimal solution, denoted by $\{X_i^*\}_{i=1}^n,\{V_i^*\}_{i = 1}^n,\{S_j^*\}_{j=1}^m$.
\item The number of non-zero variables in the optimal solution is exactly $m+n$, i.e., 
\[
\left| \{i\in[n]:X^*_i>0\}  \right| + \left| \{i\in[n]: V_i^*>0\}  \right| + \left| \{j\in[m]: S_j^*>0\}  \right|= m + n.
\]
\end{enumerate}
When LP~\eqref{eq:finite-support-fluid-lp} is non-degenerate, we define its \emph{optimal basis} as the set of variable indices corresponding to strictly positive entries in the solution:
\[
\left\{(X,i):i\in[n],X_i^*>0\right\} \cup  
    \left\{(V,i):i\in[n],V_i^*>0\right\} \cup  
    \left\{(S,j):j\in[m],S_j^*>0\right\}  .
\]    
\end{definition}

\subsection{Non-Degenerate Distributions} \label{prelim:non-degenerate-distribution}
\emph{Non-degenerate distributions} in the OLP context were first proposed and studied by \citet{li2022online}. It can be viewed as a generalization of the non-degeneracy condition introduced in \citet{jasin2012re, jasin2015performance}, as well as a stochastic analogue of the general-position condition used in \citet{devanur2009adwords, agrawal2014dynamic}. This formulation has since been widely adopted in subsequent work \citep{gao2024decoupling, bray2025logarithmic, ma2025optimal}. In this subsection, we extend this notion in a natural way to the replenishment-based setting.

Let $\mathcal{S} = \{(r_t,\bm{a}_t, \bm{b}_t)\}$ be a sequence of order and replenishment tuples, and let $\bm{B}' = ({B}'_1, \dots B'_m)\in \mathbb{R}^m_{\geq 0}$. For any $\bp = (p_1, \dots, p_m) \in \mathbb{R}^m_{\geq 0}$, define
\begin{align} \label{eq:def-f}
    f(\bp; \bm{B}',\mathcal{S}) \defeq \left\langle \bp, \bm{B}'\right\rangle + \frac{1}{\left|\mathcal{S}\right|}\sum_{(r,\bm{a},\bm{b}) \in \mathcal{S}}   \left[r - \langle\bm{a}, \bp\rangle\right]^+  .
\end{align}
Using this formulation, the dual of the relaxed offline LP~\eqref{eq:basic-offline-relax-dual} can be equivalently rewritten as:
\begin{align}
\min_{\bp \geq \bm{0}} \left\{ f\left(\bp; \frac{1}{T}\sum_{t=1}^T \bm{b}_t,\{(r_t,\bm{a}_t, \bm{b}_t)\}_{t=1}^T \right)\right\} .\label{eq:basic-offline-relax-dual-simple} 
\end{align}
We also need to extend the definition in Eq.~\eqref{eq:def-f} from empirical samples to distributions. Specifically, for any $\bm{B}' \in \mathbb{R}^m_{\geq 0}$ and any distribution $\mathcal{P}$ over $(r,\bm{a},\bm{b})$, define
\begin{align}
    f(\bp; \bm{B}', \mathcal{P})\defeq\langle \bp, \bm{B}'\rangle +\mathbb{E}_{(r,\bm{a},\bm{b})\sim\mathcal{P}} \left[ r- \langle \bm{a}, \bp \rangle  \right]^+.\label{eq:def-f-ext-distribution}
\end{align}
Let $\Omega_b \defeq \times_{i=1}^m (\underline{b},\bar{b}) \subseteq \mathbb{R}_{\geq 0}^m$ be the set of admissible replenishment levels, and define the feasible set for dual variables as $\Omega_p\defeq \{ \bp\in\mathbb{R}_{\geq 0}^m:\sum_{j=1}^m p_j\leq {\bar{r}}/{\underline{b}} \} \subseteq \mathbb{R}_{\geq 0}^m$. We now formalize the definition of non-degeneracy for a distribution $\mathcal{P}$.
\begin{definition}[Non-Degenerate Distribution]
\label{def:cts-non-degeneracy}
    The distribution $\mathcal{P}$ over $(r,\bm{a},\bm{b})$ is said to be \emph{non-degenerate} if it is bounded (in the sense of Definition~\ref{def:bounded-dis}) and satisfies the following conditions:
    \begin{enumerate}
        \item The second-order moment matrix $\mathbb{E}_{(r,\bm{a},\bm{b})}\left[\bm{a}\bm{a}^{\top}\right]$ is positive-definite. Denote its minimum eigenvalue by $\lambda_{\min} > 0$.
        \item There exist constants $\lambda,\mu>0$ such that for any $\bp'\in\Omega_p$, $\bm{B}'\in\Omega_b$, and any $\bp^* \in\arg\min_{\bp\geq\bm{0}} f(\bp; \bm{B}', \mathcal{P})$, let $(\widetilde{r},\widetilde{\bm{a}},\widetilde{\bm{b}})\sim\mathcal{P}$ and it holds that
        \begin{align*}
            \lambda\left|\langle \widetilde{\bm{a}},\bp'\rangle -\langle \widetilde{\bm{a}}, \bp^*\rangle \right|  \leq \left| \Pr\left(r>\langle \widetilde{\bm{a}}, \bp'\rangle \middle|\widetilde{\bm{a}}  \right) - \Pr\left(r>\langle \widetilde{\bm{a}}, \bp^*\rangle \middle|\widetilde{\bm{a}}  \right)  \right|\leq  \mu\left|\langle \widetilde{\bm{a}},\bp'\rangle -\langle \widetilde{\bm{a}}, \bp^*\rangle \right|.
        \end{align*}
        \item For any $\bm{B}'\in\Omega_b$, and any $\bp^*\in\arg\min_{\bp\geq\bm{0}} f(\bp; \bm{B}', \mathcal{P})$, it satisfies that for each $j\in[m]$
        \begin{align}
            p_j^*=0\iff B'_j>\mathbb{E}_{(r,\bm{a},\bm{b})\sim\mathcal{P}}\left[ a_j\cdot\mathbb{I}\left( r>\langle \bm{a}, \bp^*\rangle \right) \right], \label{eq:cts-no-de-no-bingding-def} \\
            p_j^*>0\iff B'_j=\mathbb{E}_{(r,\bm{a},\bm{b})\sim\mathcal{P}}\left[ a_j\cdot\mathbb{I}\left( r>\langle \bm{a}, \bp^*\rangle  \right) \right]. \label{eq:cts-no-de-bingding-def} 
        \end{align}
    \end{enumerate}
\end{definition}
\begin{remark}
\label{rmk:cts-no-de-auto-hold-equal}
Equation~\eqref{eq:cts-no-de-bingding-def} is logically implied by Part 1, Part 2, and Eq.~\eqref{eq:cts-no-de-no-bingding-def} of Definition~\ref{def:cts-non-degeneracy} (please refer to Section~\ref{sec:proof-rmk-cts-no-de-auto-hold-equal} for the proof). However, we retain it explicitly in the definition for completeness and to highlight the intuitive connection to complementary slackness conditions in linear programming.
\end{remark}

\section{$\widetilde{\mathcal{O}}(\sqrt{T})$ Regret for Bounded Distributions}
\label{sec:bounded-distribution}

In this section, we assume only that the distribution $\mathcal{P}$ is bounded and impose no other conditions. We design an algorithm to achieve $\widetilde{\mathcal{O}}(\sqrt{T})$ regret.

\subsection{The Algorithm}
\label{sec:main-alg-bounded-distribution}

Our algorithm is described in Algorithm~\ref{alg:bounded-distribution}. At a high level, the algorithm starts with a conservative \emph{warm-up phase}, during which all orders are rejected for the first $\kappa = \Theta(\sqrt{T \ln T})$ periods to accumulate resources. The exact quantity of $\kappa$, along with the appropriate scaling factors $W$ and $C$, is computed in Line~\ref{line:bounded-distribution-initial-parameters}.
After this initial phase, the algorithm enters the \emph{main decision phase}, employing an adaptive acceptance rule based on estimated dual prices $\bp^t = (p^t_1, \dots, p^t_m)$. Specifically, at each period $t > \kappa$, an order is accepted if its reward $r_t$ exceeds the estimated dual cost $\langle \bp^t, \bm{a}_t\rangle$, provided that sufficient resources are available. The dual prices $\bp^t$ are updated iteratively via projected gradient descent (Line~\ref{line:bounded-update-dual}), ensuring non-negativity and gradually refining the estimation of the optimal dual prices.

In the classical online linear programming (OLP) setting where the inventory is fully stocked at the beginning, algorithms that estimate the dual prices via gradient descent have been investigated (see, e.g., \citet{li2023simple}). The difference between their setting and ours lies in the availability of inventory information: while their algorithms assume the total inventory (or budget) is known in advance, our setting features \emph{unknown and stochastic} replenishment over time. To address this, we replace the unknown expected replenishment with observed realizations in the online gradient descent update of the dual prices (Line~\ref{line:bounded-update-dual}). Despite this uncertainty, we are still able to establish a regret bound of $\tilde{\mathcal{O}}(\sqrt{T})$ for our algorithm.

\begin{remark}\label{remark:alg-require-problem-parameters}
The algorithm takes several problem parameters, such as $\bar{r}$, $\bar{a}$, $\bar{b}$, and $\underline{b}$, as inputs. These parameters may be estimated on the fly while still preserving the target regret guarantees. However, for simplicity and due to space limitations, we assume access to these parameters and focus on presenting the core ideas of the algorithm. 
\end{remark}

\algtext*{EndWhile}
\algtext*{EndIf}
\algtext*{EndFor}

\begin{algorithm}[t]
\caption{Algorithm for bounded distributions}
\label{alg:bounded-distribution}
\begin{algorithmic}[1]
\Require Boundedness parameters $\bar{r}$, $\bar{a}$, $\bar{b}$, $\underline{b}$, time horizon $T$.
\State Set the parameter $W\leftarrow 2+\left\lceil\max\left(\frac{8\sqrt{m}\bar{r}}{\underline{b}},\frac{24\sqrt{m}\left(\bar{b}+\bar{a}\right)\bar{b}^2}{\underline{b}^2},\frac{\bar{r}+2m\left(\bar{b}^2+\bar{a}^2 \right)}{\sqrt{m}\left(\bar{b}+\bar{a}\right)}\right)\right\rceil$, $C\leftarrow 9$, $\kappa\leftarrow \lceil \frac{4W\sqrt{CT\ln{T}}}{\underline{b}}\rceil$. \label{line:bounded-distribution-initial-parameters}
\State Initialize $p^1_j\leftarrow 0,\ell_{j,1}\leftarrow 0$ for all $j\in[m]$.
\For{$t = 1$ to $T$} 
    \State Observe $(r_t,\bm{a}_t,\bm{b}_t)$.
    \State \textbf{if} {$t\leq \kappa$} \textbf{then} Set $x_t\leftarrow 0$ and  $p_j^{t+1}\gets 0$ for all $j \in [m]$. \Comment{warm-up phase}
    \State ~~~~~~~~~~~~\textbf{else}  Set 
            $x_t \leftarrow
                    \mathbb{I}\left(\forall j\in[m],~\ell_{j,t} \geq a_{jt}\right) \cdot \mathbb{I}(r_t>\langle \bp^t, \bm{a}_t\rangle)$, \Comment{main decision phase}
        \State\label{line:bounded-update-dual} ~~~~~~~~~~~~~~~~~~~and $p_{j}^{t+1}\gets \left[p_{j}^t-\frac{1}{\sqrt{CT\ln{T}}}\cdot\left[b_{jt}-a_{jt}\cdot\mathbb{I}\left(r_t>\langle \bp^t, \bm{a}_t\rangle \right)\right]\right]^+$  for all $j \in [m]$.
    \State\label{line:bounded-update-left-inventory} Set $\ell_{j,t+1}\gets \ell_{j,t} + b_{jt}-a_{jt}x_t$ for all $j \in [m]$.
\EndFor
\end{algorithmic}
\end{algorithm}

\subsection{Regret Analysis}

We begin by analyzing the dual price sequence $\{\bp^t\}$. The following lemma provides a high-probability upper bound on the norm of the estimated dual prices throughout the time horizon. The proof is deferred to Section~\ref{sec:pf-bounded-distribution-bounded-dual-price}.
\begin{lemma}
\label{lem:bounded-distribution-bounded-dual-price}
If $C\geq 9$ and $T\geq m\cdot\left(\bar{b}+\bar{a}\right)^2$, then 
$
\Pr\left(\forall t\in[T+1],~\left\|\bp^t\right\|_2\leq 2W\right)\geq 1-\frac{m}{T}$.
\end{lemma}

The next lemma shows that, with high probability, making decisions solely based on the estimated dual prices avoids any out-of-stock violations after the initial warm-up phase. The proof is deferred to Section~\ref{sec:pf-of-bounded-distribution-no-out-of-stock}.
\begin{lemma}
\label{lem:bounded-distribution-no-out-of-stock}
For $T\geq \bar{a}^2+3$, we have
$        \Pr\left(\forall t>\kappa, r_t>\langle \bp^t, \bm{a}_t\rangle \Rightarrow \forall j \in [m],  \ell_{j,t} \geq a_{jt}  \right)\geq 1-\frac{2m}{T}$.
\end{lemma}

Combining Lemmas~\ref{lem:bounded-distribution-bounded-dual-price} and~\ref{lem:bounded-distribution-no-out-of-stock}, we obtain the following regret bound for Algorithm~\ref{alg:bounded-distribution}. 
\begin{theorem}
\label{thm:bounded-dis-upper}
Let $\bm{\pi}_{\mathrm{bounded}}$ denote the policy defined by Algorithm~\ref{alg:bounded-distribution}. If the underlying distribution $\mathcal{P}$ is bounded, then the regret satisfies $\Reg_{\mathcal{P}}(\bm{\pi}_{\mathrm{bounded}})=\widetilde{\mathcal{O}}(\sqrt{T})$.
\end{theorem}
\proof{Proof.} Consider $f(\bp) = f(\bp; \bm{B}, \mathcal{P})$ as defined in Eq.~\eqref{eq:def-f-ext-distribution}
and any $\bp^* \in \arg\min_{\bp\geq\bm{0}} f(\bp; \bm{B}, \mathcal{P})$. We first prove that $Tf(\bm{p^{*}})$ upper bounds $\mathbb{E}\left[R^{*}\right]$. This is because
\begin{align}
    \mathbb{E}\left[R^* \right]\leq \mathbb{E}\left[ f(\bp^*;\{(r_t, \bm{a}_t, \bm{b}_t\}_{t=1}^T\right]
    = \sum_{t=1}^T \mathbb{E}_{(r_t,\bm{a}_t, \bm{b}_t) \sim \mathcal{P}}\left[\langle \bp^*, \bm{b}_{t}\rangle+\left[r_t-\langle \bm{a}_{t}, \bp^*\rangle \right]^+\right]=Tf(\bp^*),
\end{align}
where the inequality is because of strong duality and that the dual LP~\eqref{eq:basic-offline-relax-dual-simple} is a relaxation of LP~\eqref{eq:basic-offline-LP}. 

Next, we have that
\begin{align}
    \left\|\bp^{t+1}\right\|_2^2 &\leq \left\|\bp^{t} - \frac{1}{\sqrt{CT\ln{T}}}\left[ \bm{b}_{t}-\bm{a}_{t}\cdot\mathbb{I}\left(r_t>\langle \bp^t, \bm{a}_t \rangle\right) \right]\right\|_2^2 \notag \\
    & \leq \left\|\bp^t\right\|_2^2 + \frac{\left( \bar{b}+\bar{a} \right)^2}{CT\ln{T}}-\frac{2}{\sqrt{CT\ln{T}}}\left[\langle \bm{b}_t,\bp^t \rangle-\langle \bp^t, \bm{a}_t \rangle\cdot\mathbb{I}\left( r_t>\langle \bp^t, \bm{a}_t \rangle \right) \right],
\end{align}
which indicates that
\begin{align}
    \sum_{t=\kappa+1}^T \left[\langle \bm{b}_t,\bp^t \rangle-\langle \bp^t, \bm{a}_t \rangle\cdot\mathbb{I}\left( r_t>\langle \bp^t, \bm{a}_t \rangle \right) \right]\leq \frac{\sqrt{CT\ln{T}}}{2}\cdot T\cdot  \frac{\left( \bar{b}+\bar{a} \right)^2}{CT\ln{T}} \leq \left(\bar{b}+\bar{a}\right)^2 \sqrt{T}. \label{eq:bounded-dis-control-left-resource}
\end{align}

Finally, we estimate the expectation of the reward obtained by Algorithm~\ref{alg:bounded-distribution}:
\begin{align}
   &\quad\mathbb{E}\left[R^{\bm{\pi}_{\mathrm{bounded}}}\right] \geq \mathbb{E}\sum_{t=\kappa+1}^T r_t\cdot\mathbb{I}\left(r_t> \langle \bp^t, \bm{a}_t\rangle\right)- \mathbb{E}\sum_{t=\kappa+1}^T r_t\cdot\mathbb{I}\left(x_t\neq\mathbb{I}\left(r_t> \langle \bp^t, \bm{a}_t\rangle\right)\right)  \notag \\
   &\geq \mathbb{E}\sum_{t=\kappa+1}^T r_t\cdot\mathbb{I}\left(r_t> \langle \bp^t, \bm{a}_t\rangle\right)-2m\bar{r}\notag\\
   &= \mathbb{E}\sum_{t=\kappa+1}^T \left[\langle \bm{b}_t,\bp^t \rangle + \left(r_t-\langle \bp^t, \bm{a}_t\rangle \right) \cdot \mathbb{I}\left(r_t> \langle \bp^t, \bm{a}_t\rangle \right) \right]  -  \mathbb{E}\sum_{t=\kappa+1}^T \left[\langle \bm{b}_t,\bp^t \rangle-\langle \bp^t, \bm{a}_t\rangle \cdot \mathbb{I}\left( r_t>\langle \bp^t, \bm{a}_t\rangle \right)\right]  - 2m\bar{r} \notag\\
   &\geq \mathbb{E}\sum_{t=\kappa+1}^T f\left(\bp^t\right)  - \left(\bar{b}+\bar{a}\right)^2 \sqrt{T}-2m\bar{r}, \label{eq:bounded-dis-last-eq}
\end{align}
where the second inequality is due to Lemma~\ref{lem:bounded-distribution-no-out-of-stock}, and the last inequality is by the definition of $f(\cdot)$ and Eq.~\eqref{eq:bounded-dis-control-left-resource}. Notice that $f(\bm{0})\leq \bar{r}$ and $\min(f(\bm{0}),f(\bp^t))\geq f(\bp^*)$. Continuing with Eq.~\eqref{eq:bounded-dis-last-eq}, we have
\begin{align*}
     \mathbb{E}\left[R^{\bm{\pi}_{\mathrm{bounded}}}\right]&\geq \mathbb{E}\sum_{t=1}^T f\left(\bp^*\right)  - \kappa\bar{r}- \left(\bar{b}+\bar{a}\right)^2 \sqrt{T}-2m\bar{r} \notag \\
     &\geq  \mathbb{E}\left[R^* \right]- 2m\bar{r}-\frac{4W\sqrt{CT\ln{T}}}{\underline{b}}\bar{r}-\bar{r}-\left(\bar{b}+\bar{a}\right)^2 \sqrt{T}= \mathbb{E}\left[R^*\right]-\widetilde{\mathcal{O}}(\sqrt{T}). \Halmos
\end{align*}

\section{Regret Bounds for Finite-Support Distributions}
\label{sec:finite-support-distribution}
In this section, we further assume the underlying distribution $\mathcal{P}$ has finite support. When the induced linear program (Definition~\ref{def:induced-lp}) is non-degenerate, we design an algorithm that achieves $\mathcal{O}(\log{T})$ regret. In contrast, we also construct a hard instance in which the induced LP is degenerate, and show that \emph{any} algorithm must incur at least $\Omega(\sqrt{T})$ regret in this case. The comparison between our upper and lower bounds highlights the critical role of the non-degeneracy property in determining the optimal regret for online linear programming with replenishment under finite-support distributions.

Moreover, this lower bound result stands in sharp contrast to the classical online linear programming (OLP) problem \emph{without} replenishment, where \citet{xie2025benefits,li2024infrequent} design algorithms that achieve $\mathcal{O}(1)$ regret for all finite-support distributions. The comparison highlights a key insight: replenishing inventory gradually over time, as is common in practice (e.g., daily restocking), rather than stocking it fully at the outset, as assumed in classical OLP, may introduce additional challenges for the decision-maker in the presence of demand uncertainty.

\subsection{Logarithmic Regret Upper Bound for the Non-Degenerate Case}
\label{sec:finite-dis-log-no-de}

We begin by establishing a stability result for the induced LP under the non-degeneracy assumption.

\begin{lemma}
\label{lem:finite-support-non-degenerate-unique}
Let $\mathcal{P}$ be a finite-support distribution, and suppose that its induced linear program~\eqref{eq:finite-support-fluid-lp} is non-degenerate. Then there exists a constant $L > 0$, depending only on $\{(R_i,\bm{A}_i,\mu_i)\}_{i=1}^n$ and $\{B_j\}_{j=1}^m$, such that for any perturbed vectors $\widehat{\bm{B}}=(\widehat{B}_j)_{j\in[m]},\widehat{\bm{\mu}}=(\widehat{\mu}_i)_{i\in[n]}$ satisfying
$\max_{j\in[m]}\left|\widehat{B}_j - B_j\right|\leq L$ and $\max_{i\in[n]}\left| \widehat{\mu}_i-\mu_i \right|\leq L$,
we have that the perturbed $\mathrm{LP}^{\mathrm{induced}}(\widehat{\bm{B}}, \widehat{\bm{\mu}})$ admits a unique optimal solution and preserves the same optimal basis as the original $\mathrm{LP}^{\mathrm{induced}}({\bm{B}}, {\bm{\mu}})$.
\end{lemma}
Let $\left(X_i^*(\widehat{\bm{B}},\widehat{\bm{\mu}}), V_i^*(\widehat{\bm{B}},\widehat{\bm{\mu}}),S_i^*(\widehat{\bm{B}},\widehat{\bm{\mu}}) \right)$ denote the optimal solution to the perturbed $\mathrm{LP}^{\mathrm{induced}}(\widehat{\bm{B}},\widehat{\bm{\mu}})$. The following lemma shows that the solution varies linearly with respect to the perturbations.
\begin{lemma}
\label{lem:finite-support-non-degeneracy-linear-solution}
Under the same assumptions as Lemma~\ref{lem:finite-support-non-degenerate-unique}, let  $\{\widehat{\bm{B}}^k=(\widehat{B}_j^k)_{j\in[m]},\widehat{\bm{\mu}}^k=(\widehat{\mu}_i^k)_{i\in[n]}\}_{k =1}^K$ ($k\in[K]$) be a collection of perturbed parameters such that
$\max_{j\in[m]}\left|\widehat{B}_j^k - B_j\right| \leq L$ and $\max_{i\in[n]}\left| \widehat{\mu}_i^k-\mu_i \right|\leq L$ hold for all $k \in [K]$, and let $\{w_k\}_{k=1}^K$ be any set of nonnegative weights satisfying $\sum_{k=1}^K w_k=1$. Then the optimal solution to the LP with averaged parameters satisfies:
\begin{align}
    X_i^* \left(\sum_{k\in[K]}w_k\widehat{\bm{B}}^k,\sum_{k\in[K]}w_k\widehat{\bm{\mu}}^k\right) &=\sum_{k\in[K]}w_k X_i^*\left(\widehat{\bm{B}}^k,\widehat{\bm{\mu}}^k\right), \qquad\forall i\in[n],\notag\\
    V_i^* \left(\sum_{k\in[K]}w_k\widehat{\bm{B}}^k,\sum_{k\in[K]}w_k\widehat{\bm{\mu}}^k\right) &=\sum_{k\in[K]}w_k V_i^*\left(\widehat{\bm{B}}^k,\widehat{\bm{\mu}}^k\right), \qquad\forall i\in[n],\notag\\
    S_j^* \left(\sum_{k\in[K]}w_k\widehat{\bm{B}}^k,\sum_{k\in[K]}w_k\widehat{\bm{\mu}}^k\right) &=\sum_{k\in[K]}w_k S_j^*\left(\widehat{\bm{B}}^k,\widehat{\bm{\mu}}^k\right),\qquad\forall j\in[m].\notag
\end{align}
\end{lemma} 
The proofs of Lemmas~\ref{lem:finite-support-non-degenerate-unique} and~\ref{lem:finite-support-non-degeneracy-linear-solution} are deferred to Section~\ref{sec:pf-finite-support-non-degenerate-unique}.

We now present Algorithm~\ref{alg:finite-support-distribution}, which is designed for finite-support distributions whose induced linear program is non-degenerate. The algorithm takes a few problem parameters as input, and the discussion in Remark~\ref{remark:alg-require-problem-parameters} applies here as well. The time horizon is partitioned into batches of $\kappa = \Theta(\log T)$ consecutive periods, such that each batch $z \in \{1, 2, 3, \dots \}$ consists of time periods $t \in ((z-1) \kappa, z\kappa)]$.

Similar to Algorithm~\ref{alg:bounded-distribution}, Algorithm~\ref{alg:finite-support-distribution} begins with a conservative \emph{warm-up phase} during which all orders are rejected. However, unlike the longer warm-up phase in Algorithm~\ref{alg:bounded-distribution}, the warm-up phase here is significantly shorter, lasting only $W$ batches, i.e.,  $W\kappa =\Theta(\log T)$, where $W$ is a constant defined by the algorithm at the beginning.

After the warm-up, the algorithm enters the \emph{main decision phase}. During this phase, the algorithm maintains: (i) an \emph{always-accept set} of order types, denoted $\widehat{\mathcal{I}}_{z}^{*}$ for each batch $z$; and (ii) an \emph{acceptance budget}, denoted $\Phi_{i,t}$, for each order type $i$ at each time period $t$. At time $t$, the algorithm accepts an order if either: (i) its type belongs to the always-accept set $\widehat{\mathcal{I}}_{z}^{*}$ (Line~\ref{line:finite-start-decide}), or (ii) it has remaining acceptance budget $\Phi_{\theta_t, t}$ (Line~\ref{line:finite-dis-no-degenerate-i-not-in-satisfy}). The always-accept sets $\widehat{\mathcal{I}}_{z}^{*}$ and the corresponding acceptance budget increments are computed using the solution to $\mathrm{LP}^{\mathrm{induced}}(\widehat{\bm{B}}_{z}, \widehat{\bm{\mu}}_{z})$, where $\widehat{\bm{B}}_{z}$ and $\widehat{\bm{\mu}}_{z}$ are empirical estimates computed from the data collected in the preceding batch ending at time $t' = (z-1)\kappa$ (Line~\ref{line:finite-dis-no-degenerate-compute-optimal-batch}). Based on this LP solution, denoted $\{X_{i, z}^{*}, V_{i,z}^{*}\}_{i \in [n]}$, $\{S_{j,z}^{*}\}_{j \in [m]}$, the always-accept set $\widehat{\mathcal{I}}_{z}^{*}$ includes all types $i$ such that the corresponding slack variable  $V_{i,z}^*$ is strictly positive, i.e., those that are not part of the optimal basis. The variable $X_{i,z}^*$ indicates the expected number of accepted type-$i$ orders in batch $z$. Accordingly, the acceptance budget $\Phi_{i,t}$ is incremented by $X_{i,z}^*$ at the beginning of the batch (Line~\ref{line:finite-dis-no-degenerate-batch-update-i-not-in}).

\algtext*{EndWhile}
\algtext*{EndIf}
\algtext*{EndFor}

\begin{algorithm}[h]
\caption{Algorithm for finite-support distributions with non-degenerate induced LPs}
\label{alg:finite-support-distribution}
\begin{algorithmic}[1] 
\Require Distribution support $\{R_i, \bm{A}_i\}_{i\in[n]}$, parameters $\underline{\mu}$, $\underline{b}$, $\bar{b}$ in Def.~\ref{def:finite-support-distribution}, parameter $L$ in Lemma~\ref{lem:finite-support-non-degenerate-unique}.
\State Set the parameter $C\leftarrow 2\left(1+\frac{8+8\bar{b}}{\min(L,\underline{b},\underline{\mu})}\right)^2$, $\kappa \leftarrow \lceil C\ln{T} \rceil$, $W\leftarrow \left\lceil \frac{4n\max_{i\in[n],j\in[m]}|A_{ji}|}{\underline{b}} \right\rceil$.
\State Initialize: $\Phi_{i, 1}\leftarrow 0$ for each $i\in[n]$, $\ell_{j,1} \gets 0$ for each $j\in[m]$.
\For{$t = 1$ to $T$}
    \State Observe $(r_t,\bm{a}_t,\bm{b}_t)$ and $\theta_t$, carry over acceptance budget to next period: $\Phi_{i,t+1}\leftarrow \Phi_{i,t}$, $\forall i\in [n]$.
    \State \textbf{if} {$t\leq W \kappa$} \textbf{then}  Set $x_{t}\gets 0$. \Comment{warm-up phase} 
    \State \textbf{else} \Comment{main decision phase}
    \State ~~~~\textbf{if} {$\theta_t\in\mathcal{\widehat{I}}^*_{\lceil t/\kappa\rceil}$} \label{line:finite-start-decide} \textbf{then} Set $x_t\gets \mathbb{I}\left(\forall j\in[m],~ \ell_{j,t} \geq a_{jt}\right)$. \Comment{current batch index is $\lceil t/\kappa\rceil$}
    \State \textbf{~~~~~~~~~~~~~~~~~~else} \label{line:finite-dis-no-degenerate-i-not-in-satisfy} ~Set $x_t\gets\mathbb{I}(\Phi_{\theta_t, t}\geq 1)\cdot\mathbb{I}\left(\forall j\in[m],~ \ell_{j,t} \geq a_{jt}\right)$, update $\Phi_{\theta_t,t + 1}\gets \Phi_{\theta_t, t+1}-x_t$.
    \If{$t \geq W\kappa$ \textbf{and} $z = t/\kappa+1 \in \mathbb{Z}$} \Comment{prepare for new batch $z$ in decision phase}
        \State Compute $\widehat{\bm{B}}_{z} \leftarrow (\sum_{\tau=t-\kappa+1}^{t}b_{j\tau})_{j \in [m]}$ and $\widehat{\bm{\mu}}_{z} \leftarrow  (\sum_{\tau=t-\kappa+1}^{t}\mathbb{I}(\theta_{\tau}=i))_{i \in [n]}$.
        \State\label{line:finite-dis-no-degenerate-compute-optimal-batch} Find the optimal solution, $\{X_{i, z}^{*}, V_{i,z}^{*}\}_{i \in [n]}$, $\{S_{j,z}^{*}\}_{j \in [m]}$, to $\mathrm{LP}^{\mathrm{induced}}(\widehat{\bm{B}}_{z}, \widehat{\bm{\mu}}_{z})$.
        \State\label{line:finite-dis-no-degenerate-batch-update-i-not-in} Update $\Phi_{i,t+1} \gets \Phi_{i,t+1}+X^{*}_{i,z}$ for each $i\in[n]$, and $\mathcal{\widehat{I}}^*_{z}\gets\{i\in[n]:V_{i,z}^{*}=0\}$.
    \EndIf
    \State Set $\ell_{j,t+1}\gets \ell_{j,t} +b_{jt}-a_{jt}x_t$ for each $j\in[m]$.
\EndFor
\end{algorithmic}
\end{algorithm}

The main theorem of this subsection is as follows.
\begin{theorem}
\label{thm:finite-dis-upp}
    If the unknown distribution $\mathcal{P}$ is \emph{finite-support} with non-degenerate \emph{induced} linear program~\eqref{eq:finite-support-fluid-lp}, then the regret of Algorithm~\ref{alg:finite-support-distribution} (denoted by $\bm{\pi}_{\mathrm{finite}}$) satisfies
    $   \Reg_{\mathcal{P}}(\bm{\pi}_{\mathrm{finite}})=\mathcal{O}\left(\log{T}\right)$.
\end{theorem}
\medskip
\proof{Proof.} Without loss of generality, we assume that $\kappa$ divides $T$, i.e., $\kappa \mid T$. The analysis readily extends to arbitrary $T$, with at most an additional $\mathcal{O}(\log T)$ increase in the regret bound. Let $K\defeq T/\kappa$, and define $t_z \defeq z \cdot \kappa$ for each $0\leq z\leq K$. We also define $\widehat{L}\defeq\min\left(L,\underline{b}/4,\underline{\mu}/4\right)$,   
where $L$, $\underline{\mu}$ and $\underline{b}$ are defined in the algorithm description. We also naturally extend the definitions of the batch-wise LP solutions $\{X_{i, z}^{*}, V_{i,z}^{*}\}_{i \in [n]}$, $\{S_{j,z}^{*}\}_{j \in [m]}$ to the batches $z$ in the warm-up phase, i.e., for all $2\leq z\leq W$. 

We consider the desired event $\mathcal{E}\defeq\mathcal{E}_{\text{source}}\cap\mathcal{E}_{\text{request}}$, where 
\begin{align}
    &\mathcal{E}_{\text{source}}\defeq\left\{\forall z\in[K]\text{ and }j\in[m],~\frac{1}{\kappa}{(\widehat{\bm{B}}_{z+1})_j}\in[B_j-\widehat{L},B_j+\widehat{L}] \right\}, \\
    &\mathcal{E}_{\text{request}}\defeq\left\{  \forall z\in[K]\text{ and }i\in[n],~\frac{1}{\kappa}{(\widehat{\bm{\mu}}_{z+1})_i} \in[\mu_i-\widehat{L},\mu_i+\widehat{L}] \right\}, 
\end{align}
where $B_j=\mathbb{E}[b_{jt}]$ is the expected per-period replenishment amount for resource $j$ and $\mu_i = \Pr((r_t,\bm{a}_t)=(R_i, A_i))$ is the probability of observing a type-$i$ order, as defined in Definition~\ref{def:finite-support-distribution}. The event $\mathcal{E}_{\mathrm{resource}}$ captures the case where the average resource replenishment in each batch interval $(t_z, t_{z+1}]$ is close to its expectation; and $\mathcal{E}_{\mathrm{request}}$ similarly ensures that the empirical frequency of each order type is concentrated around its mean. The following lemma, whose proof is deferred to Section~\ref{sec:pf-finite-support-rare-event}, shows that $\mathcal{E}$ happens with high probability.

\begin{lemma}
\label{lem:finite-support-rare-event} It holds that
   $\Pr\left(\mathcal{E}\right)\geq 1-\frac{4m+4n}{T^3}$.
\end{lemma}

We next define the set of \emph{binding} order types $\mathcal{I}^*\defeq\left\{i\in[n]:V_i^*\left(\bm{B},\bm{\mu}\right)=0\right\}$, where $\bm{B}=(B_j)_{j\in[m]}$, $\bm{\mu}=\left(\mu_i\right)_{i\in[n]}$, and $V_i^*(\cdot,\cdot)$ denotes the unique optimal solution to $\mathrm{LP}^{\mathrm{induced}}(\widehat{\bm{B}}, \widehat{\bm{\mu}})$, guaranteed under the non-degeneracy condition by Lemma~\ref{lem:finite-support-non-degenerate-unique}. At a high level, our goal is to accept as many orders of types in $\mathcal{I}^*$  as possible. When the desired event $\mathcal{E}$ holds, Lemma~\ref{lem:finite-support-non-degenerate-unique} directly implies the following lemma, which shows that the algorithm can correctly recover the binding type set $\mathcal{I}^*$ using the empirical estimates computed from any batch, with high probability.
\begin{lemma}
\label{lem:finite-dis-est-right}
Under the event $\mathcal{E}$, we have that for all $z\in[K]$, it holds that $\widehat{\mathcal{I}}_{z+1}^*=\mathcal{I}^*$.
\end{lemma}

The following lemma, whose proof is deferred in Section~\ref{sec:pf-finite-dis-estimation-right-no-out-of-stock}, shows that under the desired event, all binding orders are accepted, while the acceptance budgets for non-binding types ($\Phi_{i, t}, i \notin \mathcal{I}^*$) remain consistently small throughout the horizon.
\begin{lemma}
\label{lem:finite-dis-estimation-right-no-out-of-stock-1}
    Under the event $\mathcal{E}$, the following hold: (1) for each $t>W\kappa$, the decision satisfies
        $x_t =\mathbb{I}\left(\theta_t\in\mathcal{I}^* ~\mathrm{or}~\Phi_{\theta_t,t}\geq 1\right)$; and (2) for each $t\in[T+1]$, and each $i\notin\mathcal{I}^*$, the acceptance budget $\Phi_{i,t}\leq \kappa$.
\end{lemma}

Finally, we estimate the total reward. Under the event $\mathcal{E}$, invoking Lemma~\ref{lem:finite-dis-est-right}, Lemma~\ref{lem:finite-dis-estimation-right-no-out-of-stock-1} and Line~\ref{line:finite-dis-no-degenerate-batch-update-i-not-in}, we have that for each $i\notin\mathcal{I}^*$,
\begin{align}
     \sum_{t=1}^T x_t\cdot\mathbb{I}\left( \theta_t=i \right) =\sum_{v=W}^K X_{i,v+1}^*-\Phi_{i,T+1},\label{eq:finite-dis-no-degenerate-i-not-in-satisfy-number-1}
\end{align}
and for each $i\in\mathcal{I}^*$,
\begin{align}
        \sum_{t=1}^T x_t\cdot\mathbb{I}\left( \theta_t=i \right)= \sum_{t'=W\kappa+1}^T\mathbb{I}(\theta_{t'}=i)=\sum_{z=W+1}^K (\widehat{\bm{\mu}}_{z+1})_i=\sum_{v=W+1}^K X_{i,v+1}^* ,\label{eq:finite-dis-no-degenerate-i-in-satisfy-number-1}
\end{align}
while the last equality is due to the LP constraints and $i \in \mathcal{I}^*$. Combining Eqs.~(\ref{eq:finite-dis-no-degenerate-i-not-in-satisfy-number-1},\ref{eq:finite-dis-no-degenerate-i-in-satisfy-number-1}), we have
\begin{align}
    &R^{\bm{\pi}_{\mathrm{finite}}}=\sum_{t=1}^T r_tx_t \geq \sum_{z=W}^K\sum_{i=1}^n R_i X_{i,z+1}^*- \sum_{i\notin\mathcal{I}^*}|R_i|\Phi_{i,T+1}-\sum_{i\in\mathcal{I}^*}|R_i| X^*_{i,W+1} \notag \\
    &\qquad \geq \sum_{z=1}^{K}\sum_{i=1}^n R_i X_{i,z+1}^*-\sum_{i=1}^n |R_i| \left(\kappa+\sum_{z=1}^{W-1} X_{i,z+1}^* \right)\geq   \sum_{z=1}^{K}\sum_{i=1}^n R_i X_{i,z+1}^* - n\max_{i\in[n]}|R_i|\cdot(W+3)\kappa,\label{eq:finite-dis-no-degenerate-online-reward-lower}
\end{align}
where the second inequality is due to Lemma~\ref{lem:finite-dis-estimation-right-no-out-of-stock-1}. On the other hand, under the event $\mathcal{E}$, we can estimate the offline optimal reward $R^*$ by
\begin{align}
    R^* &= T\cdot \sum_{i=1}^n R_i X_i^*\left( \sum_{z=1}^{K}\frac{\widehat{\bm{B}}_{z+1}}{T} , ~\sum_{z=1}^{K}\frac{\widehat{\bm{\mu}}_{z+1}}{T} \right)=\sum_{i=1}^n R_i \sum_{z=1}^{K} \kappa \cdot X_i^*\left( \frac{\widehat{\bm{B}}_{z+1}}{\kappa},\frac{\widehat{\bm{\mu}}_{z+1}}{\kappa} \right) =\sum_{i=1}^n R_i\sum_{z=1}^{K} X_{i,z+1}^*, \label{eq:finite-dis-no-degenerate-optimal-estimation}
\end{align}
where the first and last equality is due to Fact~\ref{fact:induced-lp-scale}, and the second one is due to the event $\mathcal{E}$ and Lemma~\ref{lem:finite-support-non-degeneracy-linear-solution}. Combining Eqs.~(\ref{eq:finite-dis-no-degenerate-online-reward-lower},\ref{eq:finite-dis-no-degenerate-optimal-estimation}), and Lemma~\ref{lem:finite-support-rare-event}, we have
\begin{align*}
    \Reg_{\mathcal{P}_{\mathrm{finite}}}(\bm{\pi})\leq 2T\max_{i\in[n]}|R_i|\cdot\left(1-\Pr\left( \mathcal{E} \right)\right) +n\max_{i\in[n]}|R_i|\cdot(W+3)\kappa= \mathcal{O}(\log T). \Halmos
\end{align*}

\subsection{$\Omega(\sqrt{T})$ Regret Lower Bound for the Degenerate Case}
\label{sec:finite-hard-instance-de}
We consider the case where the induced linear program associated with a finite-support distribution $\mathcal{P}$ is \emph{degenerate}. To illustrate the inherent difficulty in this setting, we construct a hard instance, denoted $\mathcal{P}_{\mathrm{finite}\text{-}\mathrm{hard}}$, such that $\mathcal{P}_{\mathrm{finite}\text{-}\mathrm{hard}}$ is finite-support, and even when $\mathcal{P}_{\mathrm{finite}\text{-}\mathrm{hard}}$ is known to the decision-maker, \emph{any} algorithm must incur $\Omega(\sqrt{T})$ regret. In this construction, we set $m=2$ (resources), $n=6$ (order types). The six order types are defined as follows.
\begin{align}
    \bm{A}_1=(2,2),~R_1=5,~\mu_1=\frac{1}{4};\qquad \bm{A}_2=(2,0),~R_2=3,~\mu_2=\frac{1}{8}; \qquad \bm{A}_3=(0,2),~R_3=3,~\mu_3=\frac{1}{8};\notag \\
    \bm{A}_4=(2,0),~R_4=4,~\mu_4=\frac{1}{8};\qquad \bm{A}_5=(0,2),~R_5=4,~\mu_5=\frac{1}{8};\qquad \bm{A}_6=(0,0),~R_6=0,~\mu_6=\frac{1}{4}.\notag
\end{align}
We further assume that the replenishment vector is deterministic: $\Pr(\bm{b}=(1,1))=1$. That is, in each time period, one unit of each resource is replenished. We present the main lower bound theorem of this subsection. The proof is deferred to Section~\ref{sec:pf-of-finite-dis-degenerate-hard-instance}.
\begin{theorem}
\label{thm:finite-dis-degenerate-hard-instance}
    There exist universal constants  $T_{\mathrm{finite}\text{-}\mathrm{hard}} > 0$ and $c>0$ such that when $T > T_{\mathrm{finite}\text{-}\mathrm{hard}}$, for any online algorithm $\bm{\pi}$, it holds that $\Reg_{\mathcal{P}_{\mathrm{finite}\text{-}\mathrm{hard}}}(\bm{\pi}) \geq c \sqrt{T}$. 
\end{theorem}

\section{Regret Bounds for Continuous-Support and Non-Degenerate Distributions}
\label{sec:non-degenerate-dis}
In this section, we further assume that the unknown distribution $\mathcal{P}$ is \emph{non-degenerate} (as defined in Definition~\ref{def:cts-non-degeneracy}, which implies that $\mathcal{P}$ is bounded and has continuous support). Under this setting, we design an algorithm that achieves a regret bound of $\mathcal{O}(\log^2 T)$. The key challenge in OLP with replenishment is the elevated risk of stockout due to the absence of initial inventory, which can substantially increase regret. To address this, our algorithm partitions the time horizon $T$ into two stages of approximately equal length and adopts an \emph{accumulate-then-convert} strategy, which is an instance of the accumulate-then-utilize principle, to secure sufficient inventory for most periods. 

In the first stage, the algorithm executes an \emph{inventory accumulation procedure} that collects a sufficient amount of inventory (on the order of $\Omega(\sqrt{T})$) while incurring at most $\mathcal{O}(\log^2 T)$ regret with respect to a specially defined reward metric that incorporates both immediate rewards and the value of remaining inventory (referred to as the \emph{dual-valued gain}). This is accomplished by slightly reducing the effective per-period resource budgets, thereby lowering acceptance probabilities, gradually building up inventory, and avoiding stockouts despite starting with zero inventory. In the second stage, the algorithm invokes an \emph{inventory conversion procedure} that aims to convert as much of the accumulated inventory and ongoing replenishment as possible into realized rewards. This is done by slightly increasing the effective per-period resource budgets. A carefully designed multi-batch budget adjustment scheme, combined with the $\Omega(\sqrt{T})$ inventory accumulated during the first stage, ensures that the risk of stockout remains minimal. Consequently, the regret with respect to the dual-valued gain in this stage is also bounded by $\mathcal{O}(\log^2 T)$.

The remainder of this section is organized as follows. Section~\ref{sec:non-degenerate-dis-preliminary} introduces several definitions and properties of the optimal dual prices that will be instrumental to our analysis. Section~\ref{sec:non-degenerate-dis-dual-price-framework} presents the dual-price-based decision and analysis framework that underlies both stages of our algorithm. Sections~\ref{sec:non-degenerate-dis-inventory-accumulation}–\ref{sec:non-degenerate-dis-inventory-conversion} describe the three components of the algorithm: the first and third correspond to the two stages outlined above, while the second computes empirical estimates that support the inventory conversion procedure. Finally, Section~\ref{sec:non-degenerate-dis-main-alg} integrates all components to form the complete algorithm and establishes its regret guarantee. 

Throughout the algorithm design and analysis in this section, we use constants $C_0, C_1, \dots, C_{11}$, defined in Section~\ref{sec:const-def}. These constants depend polynomially on the problem parameters and do not depend on the time horizon $T$. The algorithms in this section take a few problem parameters as input. Again, the discussion in Remark~\ref{remark:alg-require-problem-parameters} applies here as well.

\subsection{Properties about Optimal Dual Prices and Definition of Binding Resources}
\label{sec:non-degenerate-dis-preliminary}
Consider an underlying non-degenerate distribution $\mathcal{P}$. For simplicity, we write $f(\cdot;\cdot) = f(\cdot; \cdot,\mathcal{P})$, and  define the gradient of $f$ with respect to $\bp$ as $\nabla f(\bp;\bm{B}')\defeq \mathbb{E}_{(r,\bm{a},\bm{b})\sim\mathcal{P}}\left[ \bm{B}' -\bm{a} \cdot\mathbb{I}\left( r> \left\langle\bm{a} , \bp \right\rangle \right)  \right]$.

We summarize several useful lemmas from \cite{li2022online} below.

\begin{lemma}[{Proposition~1 in \cite{li2022online}}]
\label{lem:cts-non-de-bounded-opt}
    For any $\bm{B}'\in\Omega_b$, we have $\left\| \bp^*(\bm{B}') \right\|_1< \frac{\bar{r}}{\underline{b}}$.
\end{lemma}

\begin{lemma}[{Proposition~2 in \cite{li2022online}}]
\label{lem:cts-no-de-unique}
    For any $\bm{B}' \in \Omega_b$, let $\bp' \in \arg\min_{\bp \geq \bm{0}} f(\bp; \bm{B}')$, and let $\bp \in \Omega_p$. Then,
    $
        \frac{\lambda\lambda_{{\min}}}{2}\left\|\bp-\bp'\right\|_2^2 \leq f(\bp;\bm{B}')-f(\bp';\bm{B}')-\left\langle \nabla f(\bp';\bm{B}') , \bp-\bp'\right\rangle\leq \frac{\mu\bar{a}^2}{2} \left\|\bp-\bp'\right\|_2^2$. 
    Moreover, the minimizer $\arg\min_{\bp\geq\bm{0}}f(\bp;\bm{B}')$ is unique; we denote it by $\bp^*(\bm{B}')$. For brevity, we also write $\bp^* \defeq \bp^*(\bm{B})$.
\end{lemma}

\begin{definition}[Binding Resources] \label{def:binding-resource-set}
Given $\bm{B}' \in \Omega_b$, we define the binding set of resources as $I_B(\bm{B}') \defeq \{j\in[m]: p_j^*(\bm{B}')>0\}$. For simplicity, we denote $I_B \defeq I_B(\bm{B})$.
\end{definition}
Intuitively, and consistent with the third part of Definition~\ref{def:cts-non-degeneracy}, the binding set $I_B$ captures the bottleneck resources that constrain the optimal solution. 

The following lemma establishes the stability of the optimal dual price vector $\bp^*(\cdot)$ with respect to perturbations in the per-period replenishment vector.
\begin{lemma}[{Lemma 12 and Lemma 13 in \cite{li2022online}}]
\label{lem:cts-no-de-d-robust}
    Let $\bm{B}',\bm{B}''\in\Omega_b$. Then, 
    $    \left\| \bp^*(\bm{B}')-\bp^*(\bm{B}'') \right\|_2^2 \leq \frac{1}{\lambda^2\lambda_{{\min}}^2}\left\|\bm{B}'-\bm{B}''\right\|_2^2$.
    Moreover, if $I_B(\bm{B}')=I_B(\bm{B}'')$, then,  
    $   \left\| \bp^*(\bm{B}')-\bp^*(\bm{B}'') \right\|_2^2 \leq \frac{1}{\lambda^2\lambda_{\text{min}}^2}\sum_{j\in I_B(\bm{B}')}\left(B'_j-B''_j \right)^2$.
    Furthermore, there exists a constant $\delta_B>0$ such that for all $\bm{B}'\in\times_{j=1}^m \left(B_j-\delta_B, B_j+\delta_B\right)\subset\Omega_b$, the binding set remains unchanged: $I_B(\bm{B}')=I_B$.
\end{lemma}

\subsection{The Dual-Price-Based Algorithmic Framework}
\label{sec:non-degenerate-dis-dual-price-framework}

We formalize an algorithmic framework in which decisions are guided by dual prices. As shown in Algorithm~\ref{alg:dual-price-based-alg}, the algorithm operates over a time interval $[\tau_1, \tau_2]$. At each time step $t \in [\tau_1, \tau_2]$, the decision-maker selects a per-unit price vector $\bp^t$ for the $m$ resources (Line~\ref{line:dual-price-based-alg-decide-dual-price}). This step can be instantiated using various strategies, depending on the specific algorithm. The vector $\bp^t$ is commonly referred to as the \emph{dual price}, as it is typically derived from the dual LP (e.g., LP~\eqref{eq:basic-offline-relax-dual}). Once the dual price $\bp^t$ is determined, the algorithm observes the current order and makes the binary decision (accept or reject) based on a simple comparison involving the dual price (Line~\ref{line:dual-price-based-alg-observe-and-decide}). 

Dual-price-based decision rules have been widely used in the classical OLP setting without replenishment (e.g., \cite{agrawal2014dynamic,li2022online}). The main decision phase of our Algorithm~\ref{alg:bounded-distribution} for general bounded distributions also fits within the dual-price-based Algorithmic Framework~\ref{alg:dual-price-based-alg}. 

In this subsection, we aim to establish a more refined performance guarantee for Algorithm~\ref{alg:dual-price-based-alg} in the more structured setting of non-degenerate distributions. To state this guarantee, recall from Lemma~\ref{lem:cts-no-de-unique} that  $\bp^*$ is the unique minimizer in $\arg\min_{\bp\geq \bm{0}} f(\bm{p}; \bm{B})$. Then, the quantity $\sum_{t=\tau_1}^{\tau_2} r_tx_t + \left\langle \bm{\ell}_{\tau_2+1}-\bm{\ell}_{\tau_1},\bp^* \right\rangle$ represents the total gain (consisting of both the direct rewards from accepted orders and the value of inventory carried over) achieved by the algorithm when pricing resources at $\bp^*$. We refer to this expression, which was also examined in \cite{li2022online}, as the \emph{dual-valued gain}. The following lemma, which will play a central role in our analysis of non-degenerate distributions, compares this dual-valued gain to the benchmark $(\tau_2 - \tau_1 + 1) \cdot f(\bp^*; \bm{B})$, which serves as an upper bound on the expected reward of the hindsight optimal solution over $(\tau_2 - \tau_1 + 1)$ time steps.

\begin{algorithm}[t]
\caption{Dual-price-based algorithmic framework}
\label{alg:dual-price-based-alg}
\begin{algorithmic}[1] 
\Require Time interval $[\tau_1, \tau_2]$, initial inventory $\{\ell_{j,\tau_1} \geq 0\}_{j\in[m]}$ at time $\tau_1$.
\For{$t=\tau_1$ to $\tau_2$}
\State Use all information collected up to time $(t-1)$ to decide the dual price vector $\bp^t$. \label{line:dual-price-based-alg-decide-dual-price}
\State Observe $(r_t,\bm{a}_t,\bm{b}_t)$, and decide $x_t\gets \mathbb{I}\left(r_t>\left\langle \bm{a}_t,\bp^t \right\rangle\right)\cdot\mathbb{I}\left( \ell_{j,t}+b_{jt}\geq a_{jt} \right)$. \label{line:dual-price-based-alg-observe-and-decide}
\State Set $\ell_{j,t+1}\gets \ell_{j,t}+b_{jt}-a_{jt}\cdot x_t$ for all $j\in[m]$. \label{line:dual-price-based-alg-update-inventory}
\EndFor
\end{algorithmic}
\end{algorithm}

\begin{lemma}
\label{lem:dual-price-based-est}
For any non-degenerate distribution $\mathcal{P}$, any time interval $[\tau_1, \tau_2]$ and any initial inventory $\{\ell_{j,\tau_1}\geq 0\}_{j\in[m]}$, define the event  $\mathcal{E}\defeq \left\{ \forall \tau_1\leq t\leq\tau_2,  \bp^t\in\Omega_p  \wedge  x_t= \mathbb{I}\left(r_t> \langle \bm{a}_t, \bp^t \rangle \right) \right\}$ which represents the desirable scenario where all dual prices lie within the feasible region $\Omega_p$ (as defined in Section~\ref{prelim:non-degenerate-distribution}), and the decision to accept an order, solely based on the dual price, is always feasible with respect to the available inventory. We then have the following lower bound on the expected dual-valued gain:
\begin{align}
\mathbb{E} \left[\sum_{t=\tau_1}^{\tau_2} r_tx_t + \left\langle \bm{\ell}_{\tau_2+1}-\bm{\ell}_{\tau_1},\bp^* \right\rangle\right]&\geq (\tau_2-\tau_1+1)\cdot f\left(\bp^*;\bm{B}\right) - \mu\bar{a}^2\cdot\mathbb{E}  \sum_{t=\tau_1}^{\tau_2}\left\|\bp^t-\bp^*\right\|_2^2\cdot\mathbb{I}\left( \bp^t\in\Omega_p \right)\notag \\
&\qquad\qquad\qquad\qquad\qquad  - 2(\tau_2-\tau_1+1)\cdot\left(\bar{r}+\left(\bar{b}+\bar{a}\right)\cdot\frac{\bar{r}}{\underline{b}}\right)\cdot\Pr\left( \bar{\mathcal{E}} \right). \notag 
\end{align}
\end{lemma}

\proof{Proof.} We define a function $g(\bp'):\mathbb{R}_{\geq 0}^m \to \mathbb{R}$ to represent the expected dual-price gain obtained in a single time step when decisions are made using a dual price vector $\bp'$:
\begin{align}
    g(\bp')&\defeq \mathbb{E}_{(r,\bm{a},\bm{b})\sim\mathcal{P}}\left[ r\cdot\mathbb{I}\left(  r> \left\langle \bm{a},\bp'  \right\rangle\right) + \left\langle \bm{b}-\bm{a}\cdot\mathbb{I}\left(r>\left\langle\bm{a},\bp'\right\rangle\right),\bp^* \right\rangle \right]. \notag
\end{align}
In Section~\ref{sec:pf-cts-non-de-decomposition-basic-item}, we prove the following properties of $g$ under the assumption that the underlying distribution $\mathcal{P}$ is non-degenerate:
\begin{lemma}
\label{lem:cts-non-de-decomposition-basic-item}
For any $\bp'\in\mathbb{R}_{\geq 0}^m$, we have $\left|g(\bp')\right| \leq\bar{r}+\bar{r}\cdot\frac{\bar{b}+\bar{a}}{\underline{b}}$ and 
    $
        0 \leq g\left( \bp^* \right) - g\left( \bp' \right) \leq \mu \bar{a}^2 \cdot \left\| \bp'-\bp^* \right\|_2^2$.
\end{lemma}
Using the definition of $\mathcal{E}$, We lower bound the expected dual-valued gain as follows:
\begin{align}
    &\mathbb{E} \left[\sum_{t=\tau_1}^{\tau_2}  r_t x_t + \left\langle\bm{\ell}_{\tau_2+1}-\bm{\ell}_{\tau_1},\bp^*\right\rangle\right]=\mathbb{E} \sum_{t=\tau_1}^{\tau_2} \left(r_t x_t + \left\langle \bm{b}_t-\bm{a}_t x_t, \bp^*\right\rangle\right),\notag \\
    &\qquad\qquad\geq \mathbb{E}\sum_{t=\tau_1}^{\tau_2} \left[(r_t\cdot\mathbb{I}\left( r_t > \left\langle \bm{a}_t, \bp^t\right\rangle \right)+\left\langle\bm{b}_t-\bm{a}_t\cdot\mathbb{I}\left( r_t>\left\langle \bm{a}_t, \bp^t \right\rangle \right) , \bp^*\right\rangle\right] \cdot\mathbb{I}(\bp^t\in\Omega_p) \notag \\
    &\qquad\qquad\qquad\qquad\qquad\qquad\qquad\qquad\qquad\qquad - (\tau_2-\tau_1+1)\cdot(\bar{r}+(\bar{b}+\bar{a})\cdot\left\|\bp^*\right\|_1)\cdot\Pr\left( \bar{\mathcal{E}} \right),\notag \\ 
    &\qquad\qquad \geq \mathbb{E}\sum_{t=\tau_1}^{\tau_2} g(\bp^t)\cdot\mathbb{I}(\bp^t\in\Omega_p)- (\tau_2-\tau_1+1)\cdot\left(\bar{r}+\left(\bar{b}+\bar{a}\right)\cdot\frac{\bar{r}}{\underline{b}}\right)\cdot\Pr\left(\bar{\mathcal{E}}\right),\label{eq:dual-price-base-eq-1}
\end{align}
where the last inequality is due to that $\bp^t$ is independent of $(r_t,\bm{a}_t,\bm{b}_t)$ and Lemma~\ref{lem:cts-non-de-bounded-opt}. Invoking Lemma~\ref{lem:cts-non-de-decomposition-basic-item}, we have
\begin{align}
  &\mathbb{E}\sum_{t=\tau_1}^{\tau_2} g(\bp^t)\cdot\mathbb{I}(\bp^t\in\Omega_p)
     \geq \mathbb{E}\sum_{t=\tau_1}^{\tau_2} \left[g(\bp^*)\cdot\mathbb{I}(\bp^t\in\Omega_p)-\mu\bar{a}^2\left\|\bp^t-\bp^*\right\|_2^2\cdot\mathbb{I}\left( \bp^t\in\Omega_p\right)\right] \notag  \\ 
     &\geq (\tau_2-\tau_1+1)\cdot f(\bp^*;\bm{B})-\mu\bar{a}^2\cdot\mathbb{E}\sum_{t=1}^T \left\|\bp^t-\bp^*\right\|_2^2\cdot\mathbb{I}(\bp^t\in\Omega_p) -(\tau_2-\tau_1+1)\cdot \left(\bar{r}+\left(\bar{b}+\bar{a}\right)\cdot\frac{\bar{r}}{\underline{b}}\right)\cdot\Pr(\bar{\mathcal{E}}). \label{eq:dual-price-base-eq-2}
\end{align}
Combining Eq.~\eqref{eq:dual-price-base-eq-1} and Eq.~\eqref{eq:dual-price-base-eq-2}, we prove the lemma.
\hfill\Halmos

\subsection{Component I: The Inventory Accumulation Procedure}
\label{sec:non-degenerate-dis-inventory-accumulation}

\begin{algorithm}[t]
\caption{Inventory accumulation for non-degenerate distributions}
\label{alg:replenish-non-degenerate-dis-cts-saving}
\begin{algorithmic}[1] 
\Require Parameters $\bar{r},\bar{a},\bar{b},\underline{b},\lambda_{\text{min}},\lambda,\mu,\delta_B$ in Definitions~\ref{def:bounded-dis}, \ref{def:cts-non-degeneracy}, and Lemma~\ref{lem:cts-no-de-d-robust}, time horizon $H$.
\State Initialize: $p^1_{j}\gets 0$, $\ell_{j,1}\gets 0$ for all $j\in[m]$, and $\kappa\gets C_0\ln^2 H$. \Comment{warm-up phase}
\State \textbf{for} $t\leftarrow 1$ \textbf{to} $\kappa$ \textbf{do} Observe $(r_t,\bm{a}_t,\bm{b}_t)$, set $x_t \leftarrow 0$, and $\ell_{j, t+1} \leftarrow \ell_{j, t} + b_{jt}$ for all $j \in [m]$ .
\State Set $V_1\gets \kappa+1$, $V_w\gets\kappa+ C_1 \ln H\cdot2^{w-1}$ for $2\leq w\leq N_V\defeq \lfloor\log_2\left( \frac{H-\kappa}{C_1\ln H} \right)\rfloor$, and $V_{N_V+1}\gets H+1$. \label{line:alg-cts-no-de-first-phase-divide-horizon-1}
\For{$w=1$ to $N_V$} \Comment{main decision phase divided into $N_V$ batches}
    \State Set $\bp^{t} \leftarrow \bp$ for all $t\in [V_w, V_{w+1})$, where we choose an arbitrary minimizer \label{line:inventory-accumulation-decide-dual-price}
        \begin{align}
            \bp \in \arg\min_{\bp \geq \bm{0}} f\left( \bp ;  \frac{\sum_{k=1}^{V_w-1} \bm{b}_{k}}{V_{w}-1} - C_2\sqrt{\frac{\ln H}{V_{w+1}-V_w}}\cdot\bm{1} , \left\{  (r_t,\bm{a}_t,\bm{b}_t)\right\}_{t=1}^{V_w-1} \right). \label{eq:cts-no-de-choose-pt-saving-period-1}
        \end{align} 
    \For{$t=V_w$ to $V_{w+1}-1$}
        Execute Lines~\ref{line:dual-price-based-alg-observe-and-decide}-\ref{line:dual-price-based-alg-update-inventory} in the framework Algorithm~\ref{alg:dual-price-based-alg}.
    \EndFor
\EndFor
\end{algorithmic}
\end{algorithm}
The first component of our algorithm for non-degenerate distributions is the \emph{inventory accumulation procedure}, which governs decision-making during the first $H$ time steps (where $H$ will be set to $T/2$ in the final algorithm). The objective of this phase is twofold: (i) to accumulate a sufficient amount of inventory (on the order of $\Omega(\sqrt{H})$) for each resource with high probability, and (ii) to ensure that the expected dual-valued gain $\mathbb{E} [\sum_{t=1}^H x_t r_t + \langle \bm{\ell}_{H+1}, \bp^*\rangle ]$ closely tracks the benchmark $H \cdot f(\bp^*; \bm{B},)$, incurring at most $O(\log^2 H)$ regret.

The full description of the inventory accumulation procedure is provided in Algorithm~\ref{alg:replenish-non-degenerate-dis-cts-saving}, which is divided into a brief warm-up phase followed by the main decision phase. The warm-up phase spans $\kappa = \Theta(\log^2 H)$ time steps, during which all incoming orders are rejected to allow for initial inventory buildup. The main decision phase conforms to the dual-price-based framework (Algorithm~\ref{alg:dual-price-based-alg}), with Line~\ref{line:inventory-accumulation-decide-dual-price} specifying how the dual price is computed in this context. The main phase is further partitioned into $N_V = \Theta(\log H)$ batches: the first batch lasts $\Theta(\log^2 H)$ steps, and each subsequent batch has twice the length of its predecessor. Within each batch, a fixed dual price vector is computed as the minimizer of the dual objective function $f$, using empirical estimates of the replenishment and order distributions based on all observations up to the current batch (see Eq.~\eqref{eq:cts-no-de-choose-pt-saving-period-1}). A key innovation here is the introduction of the bias term, $-C_2 \sqrt{\frac{\ln H}{V_{w+1}-V_{w}}} \cdot \bm{1}$, added to the unbiased replenishment estimate $\frac{1}{V_w - 1}\sum_{k=1}^{V_w - 1} \bm{b}_k$. This carefully designed adjustment slightly increases the dual price $\bm{p}^t$, thereby reducing the acceptance probability and helping to meet the inventory accumulation target of $\Omega(\sqrt{H})$. Crucially, our analysis shows that this conservative bias only incurs an additional regret of at most $O(\log^2 H)$ in the total dual-valued gain.

To analyze the performance of Algorithm~\ref{alg:replenish-non-degenerate-dis-cts-saving}, we first present the following lemma, whose proof is deferred to Section~\ref{sec:pf-cts-no-de-first-t-2-performance}.

\begin{lemma}
\label{lem:cts-no-de-inventory-accumulation}
There exist constants $C_6, C_7, C_8, C_{11}$, defined in Section~\ref{sec:const-def}, such that when $H\geq C_{11}$, the following events hold simultaneously with probability at least $1-\frac{27m^2+22m}{H^4}$:
\begin{enumerate}
    \item\label{itm:inventory-accumulation-1} For all $t \in (\kappa, H]$, we have $\bp^t\in\Omega_p$ and $
            \left\|\bp^t-\bp^*\right\|_2 \leq C_6\sqrt{\frac{\ln H}{t-\kappa}}$, i.e., the computed dual price $\bp^t$ tracks the optimal dual price $\bp^*$, with an approximation error decaying at a rate of $\sqrt{(\ln H)/t}$.
        \item\label{itm:inventory-accumulation-2} For all $t \in (\kappa, H]$ such that $r_t>\langle \bm{a}_t, \bp^t\rangle$, the order is accepted, i.e., $x_t=1$, meaning the inventory is sufficient: $\bm{\ell}_t+\bm{b}_t\geq \bm{a}_t$.
        \item\label{itm:inventory-accumulation-3} For all $j\in[m]$, the inventory level at time $H+1$ satisfies $\ell_{j,H+1}\geq C_7\sqrt{H\ln{H}}$; moreover, for all $j\in I_B$  (the set of binding resources), we also have $\ell_{j,H+1}\leq C_8\sqrt{H\ln{H}}$.
    \end{enumerate}
\end{lemma}

The third part of Lemma~\ref{lem:cts-no-de-inventory-accumulation} establishes the accumulation of resources. Next, by combining the first two parts of Lemma~\ref{lem:cts-no-de-inventory-accumulation} with Lemma~\ref{lem:dual-price-based-est}, we derive the following lemma, which provides an upper bound on the regret in terms of the total dual-valued gain.
\begin{lemma}
\label{lem:reward-saving}
    For any non-degenerate distribution $\mathcal{P}$,  Algorithm~\ref{alg:replenish-non-degenerate-dis-cts-saving} has the following guarantee:
    \begin{align*}
        \mathbb{E}\left[ \sum_{t=1}^H r_t x_t+\langle \bm{\ell}_{H+1}-\bm{\ell}_1 , \bp^* \rangle \right]\geq H \cdot f(\bp^*;\bm{B})- \mathcal{O}(\log^2 H).
    \end{align*}
\end{lemma}
\proof{Proof.} Note that the main decision phase of Algorithm~\ref{alg:replenish-non-degenerate-dis-cts-saving} fits into the framework of Algorithm~\ref{alg:dual-price-based-alg}. We invoke Lemma~\ref{lem:cts-non-de-decomposition-basic-item} by setting $\tau_1=\kappa+1$ and $\tau_2=H$, and have
\begin{align}
    &\mathbb{E}\left[ \sum_{t=\kappa+1}^H r_t x_t+\left\langle \bm{\ell}_{H+1}-\bm{\ell}_{\kappa+1} , \bp^* \right\rangle \right] \notag\\
    &\qquad \geq (H-\kappa)\cdot f(\bp^*;\bm{B})-\mu\bar{a}^2\cdot\mathbb{E}\left[ \sum_{t=\kappa+1}^H \left\| \bp^t - \bp^* \right\|_2^2\cdot\mathbb{I}(\bp^t\in\Omega_p) \right] \notag \\
    &\qquad\qquad\qquad\qquad  - \mathcal{O}\left((H-\kappa)\cdot\left[1 - \Pr\left(\forall t \in (\kappa, H], \bp^t\in\Omega_p \wedge x_t= \mathbb{I}(r_t>\langle \bm{a}_t,\bp^t \rangle)\right)\right]\right) \notag  \\
    &\qquad \geq (H-\kappa)\cdot f(\bp^*;\bm{B})-\mathcal{O}\left(\sum_{t=1}^H\frac{\ln H}{t}\right)-\mathcal{O}(1)\geq H\cdot f(\bp^*;\bm{B})-\mathcal{O}(\log^2 H),  \notag
\end{align}
where the second inequality is due to the first two parts of Lemma~\ref{lem:cts-no-de-inventory-accumulation}. Finally, we account for the regret incurred during the first $\kappa = \mathcal{O}(\log^2 H)$ time steps, and derive the following bound:
\begin{align*}
    \mathbb{E}\left[ \sum_{t=1}^H r_t x_t+\left\langle \bm{\ell}_{H+1}-\bm{\ell}_{1} , \bp^* \right\rangle \right] &\geq \mathbb{E}\left[ \sum_{t=\kappa+1}^H r_t x_t+\left\langle \bm{\ell}_{H+1}-\bm{\ell}_{\kappa+1} , \bp^* \right\rangle \right] \geq  H\cdot f(\bp^*;\bm{B})-\mathcal{O}(\log^2 H).\Halmos
\end{align*}

\subsection{Component II: The Detection and Estimation Procedure}
\label{sec:non-degenerate-dis-detection-estimation}

The second component of our algorithm for non-degenerate distributions is the \emph{detection and estimation procedure}, which leverages $H$ independent observations $\{r_t, \bm{a}_t, \bm{b}_t\}_{t \in [H]}$ to (i) identify the set of binding resources ($\widehat{I}_B$), (ii) estimate the mean per-period replenishment $\widehat{B}_j$ for each resource $j \in [m]$, and (iii) construct an empirical estimate $\widehat{f}(\cdot;~\cdot)$ of the dual objective function $f(\cdot;\cdot,)$. 

The full procedure is detailed in Algorithm~\ref{alg:detect-binding-est-resource}. While $\widehat{\bm{B}}$ and $\widehat{f}(\cdot;~\cdot)$ are obtained via direct sample averaging, identifying the binding resource set requires a more careful design. To this end, the $H$ observations are split into two halves of approximately equal size. The first half is used to estimate the optimal dual price vector $\bp^*$ (Line~\ref{line:detect-binding-est-resource-hat-p}). Then, using $\widehat{\bp}^*$ and the second half of the data, we compute $\{\widehat{\ell}_j\}_{j \in [m]}$ to estimate the remaining inventory under a dual-price-based decision rule with the optimal dual price vector (Line~\ref{line:detect-binding-est-resource-hat-ell}). This data splitting ensures that the estimation of $\widehat{\bp}^*$ and the evaluation of its impact on inventory are based on independent samples, thereby avoiding bias due to data reuse. Finally, resources with low estimated remaining inventory are identified as the estimated binding set $\widehat{I}_B$ (Line~\ref{line:cts-no-de-define-IBIN-1}). The following lemma provides theoretical guarantees for Algorithm~\ref{alg:detect-binding-est-resource}; its proof is deferred to Section~\ref{sec:pf-cts-no-de-second-t-2-performance-pre}.

\begin{algorithm}[t]
\caption{Detection and estimation}
\label{alg:detect-binding-est-resource}
\begin{algorithmic}[1] 
\Require Same set of parameters as required in Algorithm~\ref{alg:replenish-non-degenerate-dis-cts-saving}, $H$ observations $\{r_t,\bm{a}_t,\bm{b}_t\}_{t\in[H]}$.
\State Choose any $\widehat{\bp}^* \in
   \arg\min_{\bp\geq\bm{0}} f\left( \bp; \frac{1}{\lfloor 0.5 H \rfloor}{\sum_{k=1}^{\lfloor 0.5 H\rfloor} \bm{b}_k }, \{(r_t,\bm{a}_t,\bm{b}_t)\}_{t\in[\lfloor 0.5 H\rfloor]} \right)$. \label{line:detect-binding-est-resource-hat-p}
\State For each $j\in[m]$, calculate $
    \widehat{\ell}_j\gets \sum_{k=\lfloor 0.5 H \rfloor +1}^{ H } \left[ b_{jk} - a_{jk}\cdot\mathbb{I}\left( r_k>\left\langle\bm{a}_k, \widehat{\bp}^*\right\rangle\right)\right]$. \label{line:detect-binding-est-resource-hat-ell}
\State Detect  the binding set: $\widehat{I}_B\gets\left\{j\in[m]: \widehat{\ell}_j\leq C_3\sqrt{H\ln H}\right\}$. \label{line:cts-no-de-define-IBIN-1} 
\State Estimate the mean per-period replenishment: 
$\widehat{B}_j\gets \frac{1}{H}\sum_{k=1}^H b_{jk}$ for each $j \in [m]$.
\State Define the following function $\widehat{f}(\cdot;~\cdot):\mathbb{R}_{\geq 0}^m\times \mathbb{R}_{\geq 0}^m \to \mathbb{R}$ which serves as an estimate of $f(\cdot;~\cdot,)$:
\begin{align*}
    \widehat{f}(\bp';\bm{B}')\defeq f\left(\bp';\bm{B}',\{r_t,\bm{a}_t,\bm{b}_t\}_{t\in[H]}\right).
\end{align*}
\State Output: $\widehat{I}_B$, $\widehat{\bm{B}} = (\widehat{B}_j)_{j\in[m]}$, and $\widehat{f}(\cdot;~\cdot)$.
\end{algorithmic}
\end{algorithm}

\begin{lemma}
\label{lem:cts-no-de-detection-estimation}
There exist constants $C_5$ and $C_{11}$, defined in Section~\ref{sec:const-def}, such that when $H\geq C_{11}$, the following events hold simultaneously with probability at least $1-\frac{9m+16}{H^5}$:
\begin{enumerate}
\item\label{itm:cts-no-de-detect-est-1} The binding resource set is identified correctly:  $\widehat{I}_B=I_B$.
\item\label{itm:cts-no-de-detect-est-2} The estimation error of per-period replenishment is bounded by:  $\left\|\widehat{\bm{B}}_j-\bm{B}\right\|_\infty\leq \sqrt{\frac{5\bar{b}^2\ln H}{H}}$.
\item\label{itm:cts-no-de-detect-est-3} For any $\bm{B}'\in\Omega_b$ and any $\bp'\in\arg\min_{\bp\geq\bm{0}}\widehat{f}(\bp;\bm{B}')$, it holds that $
            \left\| \bp'-\bp^*(\bm{B}') \right\|_2\leq C_5\sqrt{(\ln H)/{H}}$.
\end{enumerate}
\end{lemma}
While the first two parts of Lemma~\ref{lem:cts-no-de-detection-estimation} establish the accuracy of detecting the binding set, as well as estimating $\bm{B}$, the third part characterizes the accuracy of the estimated function $\widehat{f}(\cdot;~\cdot)$.
Specifically, it shows that for any per-period mean replenishment vector $\bm{B}' \in \Omega_b$, the optimal dual price vector of the estimated problem $\widehat{f}(\cdot;\bm{B}')$ is close to that of the true problem $f(\cdot;\bm{B}',)$.

\subsection{Component III: The Inventory Conversion Procedure}
\label{sec:non-degenerate-dis-inventory-conversion}

The third component of our algorithm for non-degenerate distributions is the \emph{inventory conversion procedure}, which governs the decision-making process during the final $H$ time steps (where, as before, $H$ is set to $T/2$ in the final algorithm). The procedure assumes a non-zero initial inventory $\bm{\ell}_1 = (\ell_{j,1})_{j \in [m]}$, which, in the final algorithm, is provided by the preceding inventory accumulation procedure. (For clarity, the time steps in this phase are indexed relatively: from $1$ to $H$ based on their offset from the start of the procedure rather than the final algorithm.) Its goal is to convert as much of the accumulated inventory, along with any replenishment received during the process, into rewards, thereby maximizing the total realized revenue.

The full details of this procedure are presented in Algorithm~\ref{alg:replenish-non-degenerate-dis-cts-consumption}. Similar to the inventory accumulation phase, the inventory conversion procedure also follows the dual-price-based framework (Algorithm~\ref{alg:dual-price-based-alg}), with Line~\ref{line:inventory-conversion-decide-dual-price} specifying the method for computing the dual price. The total duration of $H$ time steps is partitioned into $N_U = \Theta(\log H)$ batches, where the final batch contains only one time step, and each preceding batch is three times as long as its successor. In each batch $s$, a fixed dual price vector is determined as the minimizer of the estimated dual objective function $\widehat{f}(\cdot;\widetilde{\bm{B}}_s)$, where $\widetilde{\bm{B}}_s$ is computed at Line~\ref{line:inventory-conversion-decide-usable-inventory-level}. This quantity represents the effective per-period resource budget available for reward conversion and comprises two components: the estimated per-period replenishment $\widehat{\bm{B}}$, and a surplus term for each resource $j \in [m]$, defined as $\frac{\mathbb{I}(j \in\widehat{I}_B)}{U_{s+1}-U_s}\cdot (\ell_{j,U_s}-C_4\sqrt{3^{N_U+1-s}\cdot\ln{H}})$.
This surplus term designates the amount of accumulated inventory that should be consumed during batch $s$, helping to systematically convert inventory into rewards. For each binding resource, the target surplus consumption per batch is roughly $C_4 \sqrt{3^{N_U+1-s}\cdot \ln H}$, which scales as $\sqrt{H\cdot 3^{-s}}$ (up to logarithmic factors, and see the third part of Lemma~\ref{lem:cts-no-de-inventory-conversion}) and gradually decreases across successive batches. For non-binding resources, no surplus term is included, as the per-period replenishment is already sufficient to sustain the desired level of reward realization.

\begin{algorithm}[t]
\caption{Inventory conversion for non-degenerate distributions}
\label{alg:replenish-non-degenerate-dis-cts-consumption}
\begin{algorithmic}[1] 
\Require Same set of parameters as required in Algorithm~\ref{alg:replenish-non-degenerate-dis-cts-saving}, time horizon $H$, initial inventory $(\ell_{j,1})_{j \in [m]}$, and the estimates  $\widehat{I}_B$,  $(\widehat{B}_j)_{j\in[m]}$, $\widehat{f}(\cdot;~\cdot)$ computed by Algorithm~\ref{alg:detect-binding-est-resource}.
\State Initialize: $N_U\gets \lfloor\log_{3} H \rfloor$,  $U_1\gets 1$,  and $U_s\gets H+2-3^{N_U+1-s}$ for each $s \in \{2, 3, \dots, N_U\}$. \label{line:alg-cts-no-de-second-phase-divide-horizon-1}
\For{$s = 1$ to $N_U$}
    \State Set $\widetilde{B}_{j,s}\gets \widehat{B}_j+
    \frac{\mathbb{I}(j \in\widehat{I}_B)}{U_{s+1}-U_s}\cdot (\ell_{j,U_s}-C_4\sqrt{3^{N_U+1-s}\cdot\ln{H}}), \forall j \in [m]$; denote  $\widetilde{\bm{B}}_s = (\widetilde{B}_{j,s})_{j\in[m]}$.\label{line:inventory-conversion-decide-usable-inventory-level}
    \State Set $\bp^t \gets \bp$ for all $t\in[U_s,U_{s+1})$, where we choose any $\bp \in \arg \min_{\bp\geq\bm{0}}\widehat{f}\left(\bp;\widetilde{\bm{B}}_s\right)$. \label{line:inventory-conversion-decide-dual-price}
    \For{$t=U_s$ to $U_{s+1}-1$} Execute Lines~\ref{line:dual-price-based-alg-observe-and-decide}-\ref{line:dual-price-based-alg-update-inventory} in the framework Algorithm~\ref{alg:dual-price-based-alg}.
    \EndFor
\EndFor
\end{algorithmic}
\end{algorithm}

Let ${\color{black} N_U^\flat = N_U-\left\lceil \log_3\left( \max\left\{\frac{4(C_{10})^2\bar{a}^4\mu^2}{3\delta_B^2} , \frac{5 (\bar{b}+2\bar{a})^2}{\delta_B^2},  \frac{12(C_8+C_9)^2}{\delta_B^2} \right\} \cdot\ln H\right)\right\rceil}$. We now state a theoretical guarantee for Algorithm~\ref{alg:replenish-non-degenerate-dis-cts-consumption}; the proof is deferred to Section~\ref{sec:pf-cts-no-de-second-t-2-performance}.
\begin{lemma}
\label{lem:cts-no-de-inventory-conversion}
Let the constants $C_4, C_7, C_9,C_{10},C_{11}$ be defined as in Section~\ref{sec:const-def}. Suppose $H \geq C_{11}$, and the conditions in Lemma~\ref{lem:cts-no-de-detection-estimation} hold for $\widehat{I}_B$, $\widehat{\bm{B}}$, and $\widehat{f}(\cdot;~\cdot)$. If the initial inventory satisfies $\ell_{j,1}\geq C_7\sqrt{H\ln H}$ for all $j\in [m]$ and $\ell_{j,1}\leq C_8\sqrt{H\ln H}$ for $j\in I_B$, then with probability at least $1-\frac{10m}{H^3}$, the following statements hold simultaneously for every batch $s\in{1,2,\dots, N_U^\flat}$ and every time step $t\in[U_s, U_{s+1})$ within that batch:
\begin{enumerate}
\item\label{itm:cts-no-de-inventory-conversion-1} We have $\widetilde{\bm{B}}_s\in\Omega_b$,  $\bp^{t}\in\Omega_p$, and $\left\|\bp^{t}-\bp^*\right\|_2\leq C_{10}\cdot\sqrt{\frac{\ln H}{3^{N_U+1-s}}}$. That is, the constructed $\widetilde{\bm{B}}_s$ and computed dual price $\bp^t$ lie within their respective admissible and feasible sets, and the deviation of $\bp^t$ from the optimal dual price $\bp^*$ is bounded by an error rate on the order of $\sqrt{(\ln H)/(H - t + 1)}$.
\item\label{itm:cts-no-de-inventory-conversion-2} Whenever $r_t>\langle \bm{a}_t, \bp^t\rangle$, the available inventory is sufficient to accept the order, i.e., $\bm{\ell}_t+\bm{b}_t\geq \bm{a}_t$.
\item\label{itm:cts-no-de-inventory-conversion-3} We have {\color{black} $C_4 \geq 3 \cdot C_9$}, and for every binding resource $j\in I_B$, the inventory at the end of batch $s$ satisfies $\left|\ell_{j,U_{s+1}}-C_4\sqrt{3^{N_U+1-s}\cdot\ln H}\right| \leq C_9\sqrt{3^{N_U-s}\cdot\ln H}$. For every non-binding resource $j \in [m] \setminus I_B$, the inventory satisfies $\ell_{j,U_{s+1}}\geq C_7\sqrt{H\ln H}$.
\end{enumerate}
\end{lemma}

Combining Lemma~\ref{lem:dual-price-based-est} and Lemma~\ref{lem:cts-no-de-inventory-conversion}, we can derive the following lemma regarding the total dual-valued gain of Algorithm~\ref{alg:replenish-non-degenerate-dis-cts-consumption}.
\begin{lemma}
\label{lemma:reward-consumption}
Under the same conditions in Lemma~\ref{lem:cts-no-de-inventory-conversion}, we have 
\begin{align}
    \mathbb{E}\left[ \sum_{t=1}^{H} r_t x_t - \left\langle \bm{\ell}_1, \bp^* \right\rangle  \right]  \geq H \cdot f(\bp^*;\bm{B}) - \mathcal{O}(\log^2 H). \notag
\end{align}
\end{lemma}
\proof{Proof.} Let $H^\flat = U_{N_{U}^\flat+1}-1$ (the last time step of batch $N_U^\flat$. We first establish that
\begin{align}\label{eq:lem-reward-consumption-1}
\mathbb{E}\left[\sum_{t=1}^{H^\flat} r_t x_t + \left\langle \bm{\ell}_{H^\flat}-\bm{\ell}_1, \bp^* \right\rangle\right] \geq H \cdot f(\bp^*;\bm{B}) - \mathcal{O}(\log^2 H).
\end{align}
This is because, by invoking Lemma~\ref{lem:dual-price-based-est} with $\tau_1=1$ and $\tau_2=H^\flat$, we have
\begin{align*}
&\mathbb{E}\left[\sum_{t=1}^{H^\flat} r_t x_t + \left\langle \bm{\ell}_{H^\flat+1}-\bm{\ell}_1, \bp^* \right\rangle\right] \geq H^\flat \cdot f(\bp^*;\bm{B})-\mu\bar{a}^2\cdot\mathbb{E}\left[ \sum_{t=1}^{H^\flat} \left\| \bp^t - \bp^* \right\|_2^2\cdot\mathbb{I}(\bp^t\in\Omega_p) \right] \\
&\qquad\qquad\qquad\qquad\qquad\qquad\qquad\quad\quad - \mathcal{O}\left( H^\flat \cdot \left[1-\Pr\left(\forall t \in [1,H^\flat], \bp^t\in\Omega_p \wedge x_t= \mathbb{I}\left(r_t>\left\langle \bm{a}_t,\bp^t \right\rangle\right)\right)\right]\right) \\
&\qquad\qquad\qquad\qquad\qquad\geq H^\flat \cdot f(\bp^*;\bm{B})-\mathcal{O}\left(\sum_{t=1}^{H}\frac{\ln H}{t}\right)-\mathcal{O}(1)\geq H\cdot f(\bp^*;\bm{B})-\mathcal{O}(\log^2 H),
\end{align*}
where the second inequality is due to the first two parts of Lemma~\ref{lem:cts-no-de-inventory-conversion}, and the third inequality is due to $H-H^{\flat}=\mathcal{O}(\log H)$.

Next, by the definition of the binding set (Definition~\ref{def:binding-resource-set}), we have that
\begin{align}\label{eq:lem-reward-consumption-2}
\langle \bm{\ell}_{H^\flat+1}, \bp^* \rangle = \sum_{j\in I_B} \ell_{j, H^\flat+1} \cdot p^*_j \leq \mathcal{O}(\sqrt{\log H \cdot \log H}) \cdot p_j^* \leq \mathcal{O}(\log H),
\end{align}
where the first inequality is due to the third part of Lemma~\ref{lem:cts-no-de-inventory-conversion} and the second inequality is due to Lemma~\ref{lem:cts-non-de-bounded-opt}. Combining Eq.~\eqref{eq:lem-reward-consumption-1} and Eq.~\eqref{eq:lem-reward-consumption-2}, we prove the lemma.
\hfill\Halmos

\subsection{The Main Algorithm}
\label{sec:non-degenerate-dis-main-alg}

We now present our main Algorithm~\ref{alg:replenish-non-degenerate-dis-cts-combine} for  non-degenerate distributions by putting together the three components introduced above. We combine Lemma~\ref{lem:reward-saving}, Lemma~\ref{lem:cts-no-de-detection-estimation}, and Lemma~\ref{lemma:reward-consumption} to derive the following theorem on the regret of our main algorithm.

\begin{algorithm}[t]
\caption{Main algorithm for non-degenerate distributions}
\label{alg:replenish-non-degenerate-dis-cts-combine}
\begin{algorithmic}[1] 
\State For the first $\lceil \frac{T}{2} \rceil$ time periods, invoke Algorithm~\ref{alg:replenish-non-degenerate-dis-cts-saving} with the horizon parameter $H=\lceil \frac{T}{2} \rceil$.
\State Using the  $\lceil \frac{T}{2} \rceil$ observations $\left\{(r_t,\bm{a}_t,\bm{b}_t)\right\}_{t\in\lceil {T}/{2} \rceil}$ obtained in the previous step as the input, invoke Algorithm~\ref{alg:detect-binding-est-resource} to obtain the estimates  $\widehat{I}_B$, $\widehat{\bm{B}}$, and $\widehat{f}(\cdot;~\cdot)$.
\State For the last $\lfloor \frac{T}{2} \rfloor$ time periods, invoke Algorithm~\ref{alg:replenish-non-degenerate-dis-cts-consumption} with the initial inventory $\bm{\ell}_{\lceil T/2\rceil+1}$, the horizon parameter $H'=\lceil \frac{T}{2} \rceil$ and the estimates $\widehat{I}_B$, $\widehat{\bm{B}}$, and $\widehat{f}(\cdot;~\cdot)$.
\end{algorithmic}
\end{algorithm}

\begin{theorem}
\label{thm:regret-non-degenerate-cts-dis-1}
    Let $\mathcal{P}$ be a non-degenerate distribution. The regret of Algorithm~\ref{alg:replenish-non-degenerate-dis-cts-combine} (denoted by $\bm{\pi}_{\mathrm{non}\text{-}\mathrm{deg}}$) satisfies $\Reg_{\mathcal{P}}(\bm{\pi}_{\mathrm{non}\text{-}\mathrm{deg}})=\mathcal{O}(\log^2 T)$.
\end{theorem}
\proof{Proof.} In this proof, we assume w.l.o.g.~that $T\geq 2 C_{11}$ (with the justification provided in Section~\ref{sec:const-def}), which guarantees that all technical lemmas established earlier can be applied.  We use $\bm{\ell}_t$ to denote the inventory level at the beginning of period $t$ and $x_t$ as the decision made at time $t$ in Algorithm~\ref{alg:replenish-non-degenerate-dis-cts-combine}. Recall that we set $H=\frac{T}{2}$. We define the event $\mathcal{E}$
\begin{align*}
&\mathcal{E}\defeq \left\{\widehat{I}_B, \widehat{\bm{B}}, \widehat{f}(\cdot;~\cdot) \text{~satisfies the desired conditions in Lemma~\ref{lem:cts-no-de-detection-estimation}}\right\} \\
&\qquad \qquad\qquad \bigcap
    \left\{\forall j\in[m],~\ell_{j,H+1}\geq C_7\sqrt{H\ln H}\right\}\bigcap \left\{\forall j\in I_B,~\ell_{j,H+1}\leq C_8\sqrt{H\ln H}\right\}.
\end{align*}
By Lemma~\ref{lem:cts-no-de-detection-estimation} and the third part of Lemma~\ref{lem:cts-no-de-inventory-accumulation}, we have $\Pr(\mathcal{E})=1-\mathcal{O}(\frac{1}{H^2})$. Invoking Lemma~\ref{lem:reward-saving} and Lemma~\ref{lemma:reward-consumption}, we have 
\begin{align}
    &\mathbb{E}\left[\sum_{t=1}^H  \left(r_tx_t+\left\langle \bm{\ell}_{H+1}-\bm{\ell}_1 , \bp^* \right\rangle \right) + \sum_{t=H+1}^{T} \left( r_t x_t - \left\langle \bm{\ell}_{H+1} , \bp^* \right\rangle\right)\cdot\mathbb{I}(\mathcal{E})\right] \geq T \cdot  f(\bp^*;\bm{B})  - \mathcal{O}(\log^2 H). \label{eq:thm-regret-non-degenerate-1}
\end{align}
On the other hand, since $\left|\sum_{t=H+1}^{T} r_t x_t - \left\langle \bm{\ell}_{H+1} , \bp^* \right\rangle\right|\leq T\bar{r}\cdot\frac{1+\bar{b}+\bar{a}}{\underline{b}}=\mathcal{O}(T)$, we have 
\begin{align}
    \text{LHS of Eq.~\eqref{eq:thm-regret-non-degenerate-1}}
    &\leq \mathbb{E}\left[\sum_{t=1}^H  \left(r_tx_t+\left\langle \bm{\ell}_{H+1}-\bm{\ell}_1 , \bp^* \right\rangle \right) + \sum_{t=H+1}^{T} \left( r_t x_t - \left\langle \bm{\ell}_{H+1} , \bp^* \right\rangle\right)\right]+ \mathcal{O}\left(T \cdot \left(1-\Pr(\mathcal{E})\right)\right) \notag \\
    &= \mathbb{E}\left[\sum_{t=1}^T r_tx_t\right] +\mathcal{O}(1). \label{eq:thm-regret-non-degenerate-2}
\end{align}
Combining Eq.~\eqref{eq:thm-regret-non-degenerate-1} and Eq.~\eqref{eq:thm-regret-non-degenerate-2}, we conclude that 
$\mathbb{E}\left[\sum_{t=1}^T r_tx_t\right] \geq T \cdot f(\bp^*;\bm{B})-\mathcal{O}(\log^2 T)\geq \mathbb{E}\left[ R^* \right]-\mathcal{O}(\log^2 T)$. \hfill \Halmos

\section{Numerical Experiments}
\label{sec:experiment}
In this section, we implement a streamlined variant of our Algorithm~\ref{alg:replenish-non-degenerate-dis-cts-combine} for non-degenerate distributions and evaluate its empirical performance. The implementation preserves the core mechanism (inventory accumulation and inventory conversion), while adopting simplified parameter choices for practical efficiency. For comparison, we also implement the algorithm adapted from \citet{li2022online} as the baseline. The work of \citet{li2022online} studies ordinary OLP with a large initial inventory. In their algorithm, the decision $x_t$ is determined based on the current inventory level $\bm{\ell}_t^{\mathrm{OLP}}$ through a decision rule of the form $x_t = x_t^{\mathrm{OLP}}(\bm{\ell}_t^{\mathrm{OLP}}, \bm{a}_t, r_t)$. To accommodate the replenishment setting, we endow this baseline with additional knowledge of $\mathbb{E}[\bm{b}_1]$ and modify the first argument of $x_t^{\mathrm{OLP}}(\cdot, \cdot, \cdot)$ by replacing the current inventory level with the sum of the current inventory $\bm{\ell}_t^{\mathrm{OLP\text{-}R}}$ and the expected future replenishment $(T-t)\cdot \mathbb{E}[\bm{b}_1]$. Specifically, at time $t$, we accept the request if
$
x_t^{\mathrm{OLP}}\left(\bm{\ell}_t^{\mathrm{OLP\text{-}R}} + (T-t)\cdot
\mathbb{E}[\bm{b}_1], \bm{a}_t, r_t\right) = 1
$
and the available resources are sufficient.

We implement the algorithms under two non-degenerate input models adapted from \citet{li2022online}, with model details summarized in Table~\ref{tab:input}. In the first model (Random Input~I), the requirement vector $\bm{a}_t = (a_{1t}, \dots, a_{mt})$ and the reward $r_t$ are drawn \emph{i.i.d.}~from uniform distributions. This model shares the same structure as Random Input~I in \citet{li2022online}, differing only in the distribution parameters. In the second model (Random Input~II), the requirement vector $\bm{a}_t = (a_{1t}, \dots, a_{mt})$ is generated i.i.d. from a normal distribution, and the reward is given by $r_t = \epsilon_t + \sum_{j=1}^m a_{jt}$, where $\epsilon_t$ is \emph{i.i.d.}\ noise drawn from a normal distribution. This model is almost identical to Random Input~II in \citet{li2022online}, except for the additional noise term $\epsilon_t$ in the reward. In both models, the replenishment vector $\bm{b}_t = (b_{1t}, \dots, b_{mt})$ is generated \emph{i.i.d.}~from a uniform distribution; in contrast, \citet{li2022online} consider a setting with an initial inventory of $\frac{T}{4}$ for all resources, which matches the expected total replenishment in our setting.

\begin{table}[htbp]
\centering
\renewcommand{\arraystretch}{1.7}
\scalebox{1.0}{
\begin{tabular}{c|c|c|c}
\hline
\makecell[c]{Model}& \makecell[c]{$\bm{a}_t$} & \makecell[c]{$r_t$} & \makecell[c]{$\bm{b}_t$} \\ \hline
\makecell[c]{Random Input I} & \makecell[c]{$a_{jt}\overset{\textit{i.i.d.}}{\sim} \mathcal{U}[0,1]$ $^*$}  & \makecell[c]{$r_{t}\overset{\textit{i.i.d.}}{\sim} \mathcal{U}[0,10]$}  &  \makecell[c]{$b_{jt}\overset{\textit{i.i.d.}}{\sim} \mathcal{U}[0,0.5]$}  \\ \hline
\makecell[c]{Random Input II} & \makecell[c]{$a_{jt}\overset{\textit{i.i.d.}}{\sim} \mathcal{N}(0.5,1)$ $^{**}$}  & \makecell[c]{$r_{t}=\epsilon_t+\sum_{j=1}^m a_{jt},\quad \epsilon_t\overset{\textit{i.i.d.}}{\sim} \mathcal{N}(0,5)$}  &  \makecell[c]{$b_{jt}\overset{\textit{i.i.d.}}{\sim} \mathcal{U}[0,0.5]$}  \\ \hline
\multicolumn{4}{l}{\makecell[l]{\small $^{*~}$ Here $\mathcal{U}[a,b]$ denotes the uniform distribution on the interval $[a,b]$. \\
\small $^{**}$ Here $\mathcal{N}(\mu,\sigma^2)$ denotes the normal distribution with mean $\mu$ and variance $\sigma^2$.}}
\end{tabular}
}
\caption{Input Models}
\label{tab:input}
\end{table}

\begin{figure}[htbp]
    \centering
        \begin{subfigure}{\textwidth}
        \centering
        \includegraphics[width=0.32\textwidth]{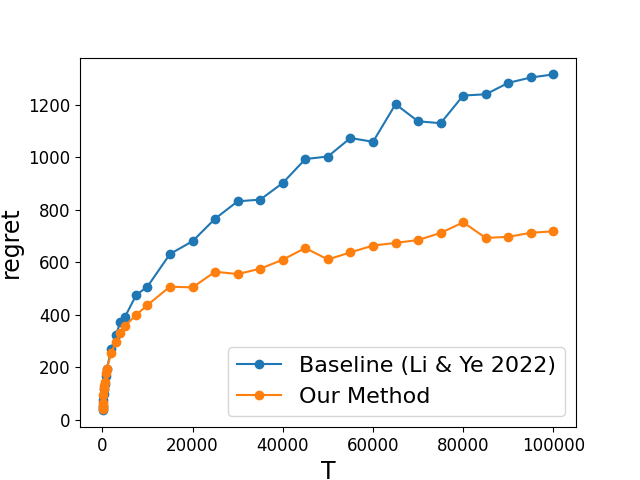}
        \includegraphics[width=0.32\textwidth]{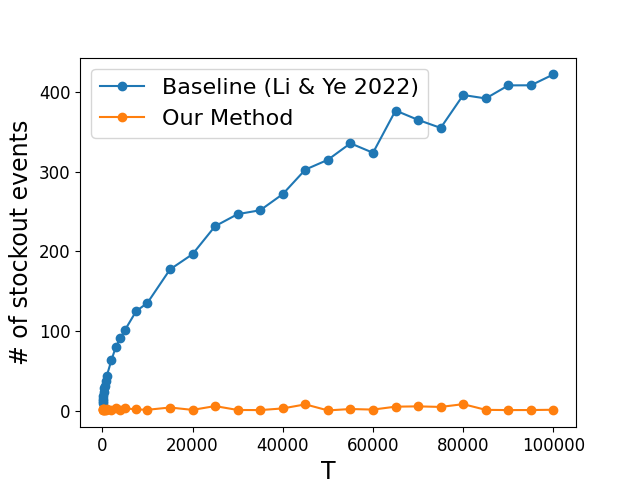}
        \includegraphics[width=0.32\textwidth]{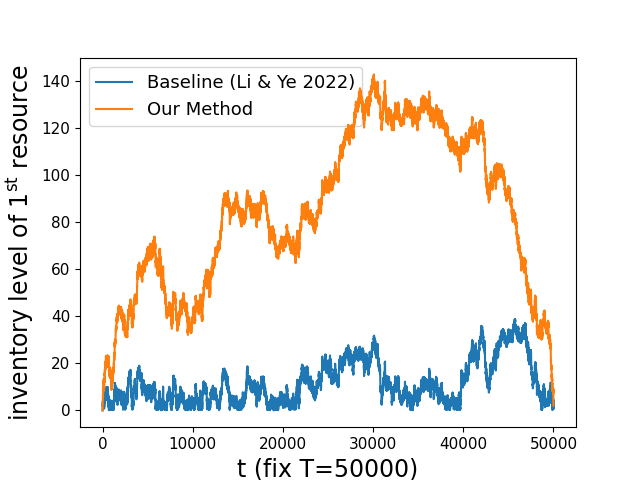}
        \caption{Random Input I with $m=5$}
    \end{subfigure}

    \vspace{0.3em}
    
    \begin{subfigure}{\textwidth}
        \centering
        \includegraphics[width=0.32\textwidth]{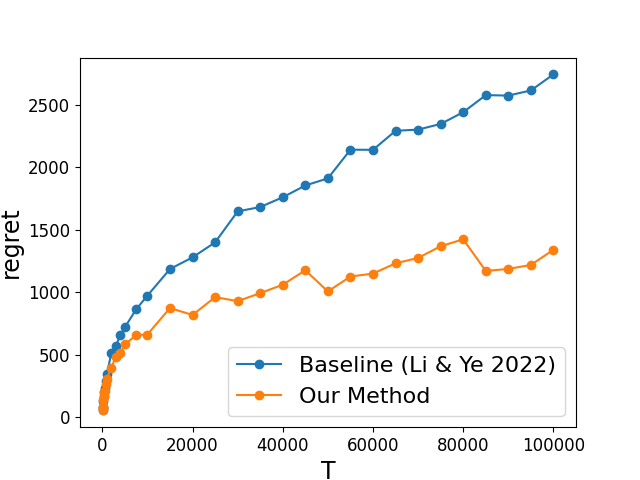}
        \includegraphics[width=0.32\textwidth]{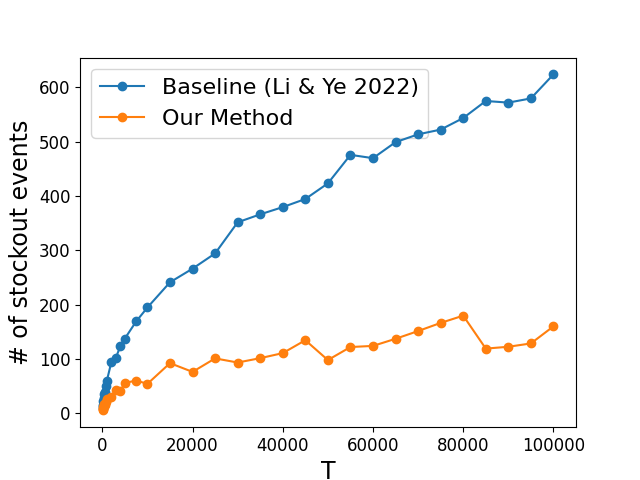}
        \includegraphics[width=0.32\textwidth]{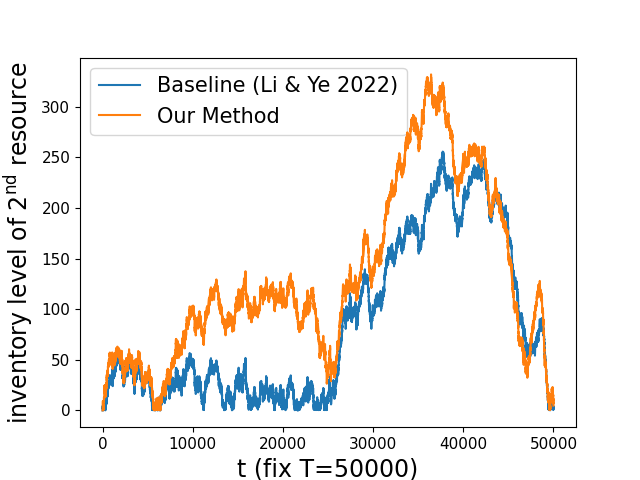}
        \caption{Random Input II with $m=5$}
    \end{subfigure}
    \caption{Empirical comparison between our algorithm and \cite{li2022online}}
    \label{fig:random-all}
\end{figure}

We present the results in Figure~\ref{fig:random-all}, where the blue curves correspond to the algorithm adapted from \citet{li2022online} and the orange curves correspond to our algorithm. In this figure, panel~(a) corresponds to Random Input~I with $m=5$, while panel~(b) corresponds to Random Input~II with $m=5$. Each panel consists of three subpanels. We additionally consider the cases of Random Input~I and Random Input~II with $m=10$; the corresponding results are reported in Section~\ref{sec:add-experiment}.

In each panel, the left figure shows the empirical regret as a function of the total time horizon $T$, with results averaged over $500$ simulation trials. In each trial, a problem instance (i.e., the sequences $\{a_{jt}\}$, $\{r_t\}$, and $\{b_{jt}\}$) is generated according to the corresponding input model. Overall, we observe that our algorithm consistently achieves lower empirical regret than the algorithm adapted from \citet{li2022online}. Although the performance improvement is not quite significant for small values of $T$, this is expected, as the theoretical gap between the two regret bounds ($\sqrt{T}$ for \citet{li2022online} and $\ln^2 T$ for our algorithm) is relatively modest in this regime. As $T$ increases, however, the empirical regret of our algorithm grows much slower, in line with the theoretical predictions. We also remark that the non-smooth appearance of the curves may be attributed the batched nature of both algorithms.

The middle figure in each panel reports the total number of stockout events across different time horizons, with results averaged over 500 simulation trials. A stockout event is defined as a situation in which a request is accepted by the dual-price-based decision rule but must ultimately be rejected due to insufficient available resources. From a theoretical perspective, one can show that the algorithm adapted from \citet{li2022online} incurs $\Omega(\sqrt{T})$ stockout events, whereas our algorithm incurs only a constant number of stockouts. The empirical results are consistent with this theoretical prediction: our algorithm experiences substantially fewer stockout events than the baseline across all considered time horizons.

Finally, the right figure in each panel depicts the inventory evolution of a representative resource over a fixed time horizon ($T=50000$) in a single trial. Under our algorithm, the inventory level initially increases to $\Theta(\sqrt{T})$ and then gradually decreases to zero. In contrast, under the algorithm adapted from \citet{li2022online}, the inventory level spends a substantial portion of time fluctuating near zero, which results in a higher likelihood of stockout events. We note that the plot shows the inventory trajectory of only one representative resource; therefore, the absence of stockouts for this particular resource at a given time does not preclude stockouts occurring in other resources. Overall, the observed inventory dynamics are consistent with our theoretical analysis and align with the expected behavior of the two algorithms.

\section{Conclusion and Future Directions}
In this work, we study a variant of the classical OLP problem in which there is no initial inventory and resources are replenished stochastically over the time horizon. We design and analyze regret-optimal algorithms (up to logarithmic factors) under several distributional settings commonly considered in the OLP literature. For future directions, a natural goal is to close the remaining logarithmic gaps between our upper and lower bounds. Another promising direction is to extend the analysis to more general input models beyond the \emph{i.i.d.} setting, such as random-permutation arrivals or non-stationary input processes.

\section*{Acknowledgements}

The authors would like to express their sincere gratitude to Prof.~Yinyu Ye for insightful and valuable discussions.

\bibliographystyle{informs2014} 
\bibliography{main} 

\ECSwitch

\ECHead{Electronic Companion}

\section{Omitted Proofs in Section~\ref{sec:bounded-distribution}}

\subsection{Proof of Lemma~\ref{lem:bounded-distribution-bounded-dual-price}}
\label{sec:pf-bounded-distribution-bounded-dual-price}
We first establish that several desired events happen with high probability. For all $t,t+s\in[T]$ with $s\geq\frac{W\sqrt{CT\ln{T}}}{12\sqrt{m}\left(\bar{b}+\bar{a}\right)}$, we first define the event 
\begin{align}
    \mathcal{E}_{j,t,t+s}\defeq\left\{ \sum_{k=t}^{t+s}b_{js}\geq\frac{4(s+1)\sqrt{m}\bar{r}}{W}\right\}. \notag
\end{align}
For all $t,t+s\in[T]$ with $s\geq\frac{W\sqrt{CT\ln{T}}}{12\sqrt{m}\left(\bar{b}+\bar{a}\right)}$, using Hoeffding's inequality, we have 
\begin{align}
    \Pr\left( \mathcal{E}_{j,t,t+s}\right)=& 1-\Pr\left(\frac{1}{s+1}\sum_{k=t}^{t+s}\left(b_{jt}-\underline{b}\right)\leq- \left(\underline{b}-\frac{4\sqrt{m}\bar{r}}{W}\right)\right) \notag \\
    \geq &1-\exp\left(-\frac{2(s+1)^2\left(\underline{b}-\frac{4\sqrt{m}\bar{r}}{W}\right)^2}{(s+1)\bar{b}^2}\right)\notag \\
    \geq &1-\exp\left(-\frac{W\sqrt{CT\ln{T}}}{6\sqrt{m}\left(\bar{b}+\bar{a}\right)}\cdot\frac{\underline{b}^2}{4\bar{b}^2}\right)
    \geq  1-\exp\left(-3\sqrt{T\ln{T}}\right)\geq1-\exp(-3\ln{T})\geq 1-\frac{1}{T^3},\label{eq:bounded-dis-bjt-sum-event-single}
\end{align}
where the second inequality is due to the condition that {\color{black} $W\geq\frac{8\sqrt{m}\bar{r}}{\underline{b}}$}, and the third inequality is due to the condition that {\color{black} $W\geq\frac{24\sqrt{m}\left(\bar{b}+\bar{a}\right)\bar{b}^2}{\underline{b}^2}$ and $C\geq 9$}. We further define 
\begin{align}
\mathcal{E}\defeq\bigcap_{t,t+s\in[T],s\geq\frac{W\sqrt{CT\ln{T}}}{12\sqrt{m}\left(\bar{b}+\bar{a}\right)}}\bigcap_{j\in[m]}\mathcal{E}_{j,t,t+s}. \notag
\end{align}
By Eq.~\eqref{eq:bounded-dis-bjt-sum-event-single}, we have that $\Pr(\mathcal{E})\geq 1-\frac{m}{T}$. 

The following lemma establishes the boundness of $\bp^t$ under the event $\mathcal{E}$; its proof is deferred to Section~\ref{sec:pf-of-bounded-dis-bounded-dual-under-event-E}.
\begin{lemma}
\label{lem:bounded-dis-bounded-dual-under-event-E}
Under the event $\mathcal{E}$, there do not exist time steps 
$K_1, K_2 \in [\kappa, T+1]$ such that  
\[
K_1 < K_2, \quad
\|\bp^{K_1}\|_2 \le W, \quad 
\|\bp^{K_2}\|_2 \ge 2W, \quad \text{and} \quad
\|\bp^{t}\|_2 \in [W,\,2W] \text{ for all } t \in (K_1, K_2).
\]
\end{lemma}

Finally, we show that under the event $\mathcal{E}$, the sequence satisfies $\|\bp^t\|_2 \le 2W$ for all $t \in [T+1]$. 
Consider any realization under the event $\mathcal{E}$, and for the purpose of contradiction, let $K_2$ denote the first index such that $\|\bp^{K_2}\|_2 \ge 2W$. 
If such an index $K_1$ exists, note that $\|\bp^t\|_2 = 0$ for all $t \le \kappa$. 
Hence, we can identify an index
\[
K_1 = \max\{ k < K_2 : \|\bp^k\|_2 \le W \} \ge \kappa.
\]
The indices $K_1$ and $K_2$ then satisfy
\[
K_2 > K_1 \ge \kappa, 
\qquad 
\|\bp^{K_1}\|_2 \le W, 
\qquad 
\|\bp^{K_2}\|_2 \ge 2W, 
\qquad 
\text{and} \quad 
\|\bp^t\|_2 \in [W, 2W] \ \text{for all } K_1 < t < K_2,
\]
which contradicts Lemma~\ref{lem:bounded-dis-bounded-dual-under-event-E}. 
Therefore, under the event $\mathcal{E}$, such $K_1$ does not exist. Thus,
\[
\Pr\!\left( \forall t \in [T+1],~ \|\bp^t\|_2 \le 2W \right)
\ge \Pr(\mathcal{E})
\ge 1 - \frac{m}{T}. \hfill \Halmos
\]

\subsection{Proof of Lemma~\ref{lem:bounded-distribution-no-out-of-stock}}
\label{sec:pf-of-bounded-distribution-no-out-of-stock}
By Line~\ref{line:bounded-update-dual} of the algorithm, we have
\begin{align}
    p_j^{s+1}
    \ge 
    p_j^{s}
    - \frac{1}{\sqrt{C T \ln T}}
    \left[
        b_{js} - a_{js} \cdot 
        \mathbb{I}\!\left(
            r_s >\left\langle \bp^s,\bm{a}_s \right\rangle 
        \right)
    \right], \notag
\end{align}
which implies that, for each $t > \kappa$,
\begin{align}
    \sum_{s = \kappa + 1}^{t}
    \left[
        b_{js} - a_{js} \cdot 
        \mathbb{I}\!\left(
            r_s >\left\langle \bp^s,\bm{a}_s \right\rangle
        \right)
    \right]
    \ge 
    -\sqrt{C T \ln T} \, p_j^{t+1}.
    \label{eq:bounded-dis-total-consumption}
\end{align}
We further define the event $\mathcal{E}_{\text{out}}$ by
\begin{align}
    \mathcal{E}_{\text{out}}
    \defeq
    \left\{
        \forall j \in [m],~
        \ell_{j, \kappa + 1} \ge 3W \sqrt{C T \ln T}
    \right\}
    \bigcap
    \left\{
        \forall t \in [T+1],~
        \|\bp^{t}\|_2 \leq 2W
    \right\}. \notag
\end{align}

We next show that the event $\mathcal{E}_{\text{out}}$ occurs with high probability; in particular,
\[
\Pr\!\left(\mathcal{E}_{\text{out}}\right) \ge 1 - \frac{2m}{T}.
\]
By Lemma~\ref{lem:bounded-distribution-bounded-dual-price}, it suffices to establish that
\begin{align}
    \Pr\!\left(\ell_{j,\kappa+1} \ge 3W\sqrt{C T \ln T}\right)
    \ge 1 - \frac{1}{T},
    \qquad \forall j \in [m]. \notag
\end{align}
Since no resources are consumed at and before time $\kappa$, we have
\[
\ell_{j,\kappa+1} = \sum_{t=1}^{\kappa} b_{jt}.
\]
Applying Hoeffding’s inequality yields
\begin{align}
    \Pr\!\left(\ell_{j,\kappa+1} < 3W\sqrt{C T \ln T}\right)
    &= \Pr\!\left(\sum_{t=1}^{\kappa} b_{jt} < 3W\sqrt{C T \ln T}\right) \notag\\
    &\le 
    \Pr\!\left(
        \frac{1}{\kappa}\sum_{t=1}^{\kappa}\!\left(b_{jt} - \underline{b}\right)
        < -\frac{\underline{b}}{4}
    \right) \notag\\
    &\le 
    \exp\!\left(-\frac{(\kappa \underline{b})^2}{8\kappa \bar{b}^2}\right)
    \le 
    \exp\!\left(-\frac{\underline{b}^2}{8\bar{b}^2}
        \cdot \frac{4W\sqrt{C T \ln T}}{\underline{b}}\right)
    \le \frac{1}{T}, \notag
\end{align}
where the first and third inequalities follow from the definition 
  $\kappa = \left\lceil \frac{4W\sqrt{C T \ln T}}{\underline{b}} \right\rceil$; the second inequality follows from Hoeffding’s inequality and the condition that $\mathbb{E}[b_{jt}] \ge \underline{b}$; the last inequality holds because $\frac{W \underline{b} \sqrt{C}}{2 \bar{b}^2} \ge 1$.

Finally, we show that under the event $\mathcal{E}_{\text{out}}$, it holds that 
\[
x_t = \mathbb{I}\!\left(r_t > \left\langle \bp^t,\bm{a}_t \right\rangle\right), 
\qquad \forall~ t \ge \kappa + 1.
\]
We establish this property by induction on $t$.

\noindent \underline{Induction basis.}  
    When $t = \kappa + 1$, because $T \ge \bar{a}^2 + 3$, we have
    $
    \ell_{j,t} \ge 3W\sqrt{C T \ln T} \ge \bar{a} \ge a_{jt}$,
    which implies that 
    $
    x_t = \mathbb{I}\!\left(r_t >\left\langle \bp^t,\bm{a}_t \right\rangle\right)$.

\noindent \underline{Induction step.}  
    Suppose that 
    $
    x_t = \mathbb{I}\!\left(r_t > \left\langle \bp^t,\bm{a}_t \right\rangle\right)$ holds for all 
    $ t = \kappa + 1, \ldots, \widehat{T} - 1$.
    Then, for $t = \widehat{T}$, we have
    \begin{align}
        \ell_{j,\widehat{T}}
        &= \ell_{j,\kappa+1}
        + \sum_{t=\kappa+1}^{\widehat{T}-1}
        \left[
            b_{jt}
            - a_{jt} \cdot 
            \mathbb{I}\!\left(
                r_t >  \left\langle \bp^t,\bm{a}_t \right\rangle
            \right)
        \right] \notag \ge 3W\sqrt{C T \ln T}
        - \sqrt{C T \ln T} \, p_j^{\widehat{T}}
        \ge W\sqrt{C T \ln T}
        \ge \bar{a}, \notag
    \end{align}
    where the first inequality follows from the definition of $\mathcal{E}_{\text{out}}$ and Eq.~\eqref{eq:bounded-dis-total-consumption};  
    the second inequality follows again from the definition of $\mathcal{E}_{\text{out}}$;  
    and the last inequality holds because $T \ge \bar{a}^2 + 3$.

Combining the above arguments completes the proof of 
Lemma~\ref{lem:bounded-distribution-no-out-of-stock}.
\hfill \Halmos

\subsection{Proof of Lemma~\ref{lem:bounded-dis-bounded-dual-under-event-E}}
\label{sec:pf-of-bounded-dis-bounded-dual-under-event-E}
Suppose for contradiction that there exist time steps $K_1, K_2 \in [\kappa, T+1]$ such that 
\[
K_1 < K_2, \quad
\|\bp^{K_1}\|_2 \le W, \quad 
\|\bp^{K_2}\|_2 \ge 2W, \quad \text{and} \quad
\|\bp^{t}\|_2 \in [W,\,2W] \text{ for all } t \in (K_1, K_2).
\]
By Line~\ref{line:bounded-update-dual} of the algorithm, for each time step $s \geq \kappa$, we have
\begin{align}
    \|\bp^{s+1}\|_2^2
    &\le 
    \left\|
        \bp^{s}
        - \frac{1}{\sqrt{C T \ln T}}
        \bigl[\bm{b}_s - \bm{a}_s x_s\bigr]
    \right\|_2^2 \notag\\
    &= 
    \|\bp^{s}\|_2^2
    + \frac{\|\bm{b}_s - \bm{a}_s x_s\|_2^2}{C T \ln T}
    - \frac{2}{\sqrt{C T \ln T}}
    \left[
        \langle \bp^s, \bm{b}_s \rangle
        - 
        \langle \bp^s, \bm{a}_s \rangle 
        \cdot 
        \mathbb{I}\!\left(
            r_s > \langle \bp^s, \bm{a}_s \rangle
        \right)
    \right] \notag\\
    &\le 
    \|\bp^{s}\|_2^2
    + \frac{2m(\bar{b}^2 + \bar{a}^2)}{C T \ln T}
    - \frac{2}{\sqrt{C T \ln T}}
    \bigl(
        \langle \bp^s, \bm{b}_s \rangle - \bar{r}
    \bigr),
    \label{eq:bounded-dis-near-dual-price-upper}
\end{align}
where the last inequality uses the boundedness assumptions 
$\|\bm{a}_s\|_2 \le \bar{a}$, 
$\|\bm{b}_s\|_\infty \le \bar{b}$, 
and $|r_s| \le \bar{r}$. Meanwhile, by the same update rule, we also have
\begin{align}
    p_j^{s+1}
    \ge 
    p_j^{s}
    - \frac{1}{\sqrt{C T \ln T}}(b_{js} - a_{js} x_s)
    \ge 
    p_j^{s}
    - \frac{\bar{b} + \bar{a}}{\sqrt{C T \ln T}},
    \label{eq:bounded-dis-near-dual-price-lower}
\end{align}
for all $j \in [m]$. Applying Eq.~\eqref{eq:bounded-dis-near-dual-price-upper} recursively over $s = K_1, \dots, K_2 - 1$, we obtain
\begin{align}
    4W^2 
    \le 
    \|\bp^{K_2}\|_2^2
    &\le 
    \|\bp^{K_1}\|_2^2
    + 
    \frac{2m (K_2 - K_1)(\bar{b}^2 + \bar{a}^2)}{C T \ln T}
    - 
    \frac{2}{\sqrt{C T \ln T}}
    \left(
        \sum_{s=K_1 + 1}^{K_2-1} 
        \langle \bp^s, \bm{b}_s \rangle 
        - (K_2 - K_1) \bar{r}
    \right), \notag
\end{align}
which further implies that
\begin{align}
    \sum_{s=K_1+1}^{K_2-1} 
    \langle \bp^s, \bm{b}_s \rangle
    \le 
    (K_2 - K_1) \bar{r}
    + 
    \frac{m(K_2 - K_1)(\bar{b}^2 + \bar{a}^2)}{\sqrt{C T \ln T}}
    - 
    \frac{3W^2}{2}\sqrt{C T \ln T}.
    \label{eq:bounded-dis-bounded-dual-price-sum-leq}
\end{align}

Next, we partition the time periods $[K_1+1, K_2]$ into consecutive batches. 
The $i$-th batch is defined as $[t_i, t_{i+1}-1]$, where
\begin{equation}
    t_i\defeq
    \begin{cases}
        K_1+1, & i=1, \\[12pt]
        t_{i-1}+\left\lfloor
            \dfrac{W\sqrt{CT\ln{T}}}{6\sqrt{m}\left(\bar{b}+\bar{a}\right)}
        \right\rfloor, 
        & i>2,~ t_{i-1}+3\cdot\left\lfloor
            \dfrac{W\sqrt{CT\ln{T}}}{6\sqrt{m}\left(\bar{b}+\bar{a}\right)}
        \right\rfloor\leq K_2, \\[12pt]
        K_2, & \text{otherwise}.
    \end{cases} \notag
\end{equation}
Let $N$ denote the total number of non-empty batches in $[K_1,K_2]$. 
It is straightforward to verify that if 
$
K\geq \left\lfloor
    \frac{W\sqrt{CT\ln{T}}}{6\sqrt{m}\left(\bar{b}+\bar{a}\right)}
\right\rfloor$,
then for each non-empty batch $i$, 
\begin{align}
(\text{length of batch $i$}) \in \left[\left\lfloor
    \frac{W\sqrt{CT\ln{T}}}{6\sqrt{m}\left(\bar{b}+\bar{a}\right)}
\right\rfloor,
3\left\lfloor
    \frac{W\sqrt{CT\ln{T}}}{6\sqrt{m}\left(\bar{b}+\bar{a}\right)}
\right\rfloor\right].
\label{eq:bounded-dis-bounded-dual-batch-size-estimation}  
\end{align}
For each non-empty batch $i$, we consider a resource index $v_i\in[m]$ such that 
$p^{t_i}_{v_i}\geq \frac{W}{\sqrt{m}}$, whose existence is guaranteed by 
the condition $\left\|\bp^{t_i}\right\|_2\geq W$. 
Invoking Eq.~\eqref{eq:bounded-dis-near-dual-price-lower}, 
we obtain that for all $0\leq u< t_{i+1}-t_i$,
\begin{align}
    p_{v_i}^{t_i+u}
    &\geq p_{v_i}^{t_i}
    -\frac{u\cdot\left(\bar{b}+\bar{a}\right)}{\sqrt{CT\ln{T}}} \geq 
    \frac{W}{\sqrt{m}}
    -\frac{3W\sqrt{CT\ln{T}}}
    {6\sqrt{m}\left(\bar{b}+\bar{a}\right)}
    \cdot\frac{\bar{b}+\bar{a}}{\sqrt{CT\ln{T}}}
    \geq\frac{W}{2\sqrt{m}}.
    \label{eq:bounded-dis-each-batch-lower-bound}
\end{align}
Combining Eq.~\eqref{eq:bounded-dis-bounded-dual-price-sum-leq} 
and Eq.~\eqref{eq:bounded-dis-each-batch-lower-bound}, we have
\begin{align}
     (K_2 - K_1)\bar{r}
     +\frac{m(K_2-K_1)\left(\bar{b}^2+\bar{a}^2\right)}{\sqrt{CT\ln{T}}}
     &-\frac{3W^2}{2}\sqrt{CT\ln{T}}
     \geq 
     \sum_{s=K_1+1}^{K_2-1}\left\langle \bp^s,\bm{b}_s\right\rangle \notag\\
     &\qquad\qquad \geq
     \sum_{i=1}^N\sum_{s=t_i}^{t_{i+1}-1}p_{v_i}^sb_{v_i s} 
     \geq 
     \sum_{i=1}^N\sum_{s=t_i}^{t_{i+1}-1}
     \frac{W}{2\sqrt{m}}b_{v_i s}.
     \label{eq:bounded-dis-bounded-dual-average-b-less}
\end{align}
Since $b_{js}\geq 0$, it follows that
\begin{align}
    (K_2 - K_1)\left(\bar{r}
    +\frac{m\left(\bar{b}^2+\bar{a}^2\right)}{\sqrt{CT\ln{T}}}\right)
    \geq
    \frac{3W^2}{2}\sqrt{CT\ln{T}}.
    \label{eq:bounded-dis-bound-dual-increase-batch-lower}
\end{align}
By the condition that 
{\color{black} $W\geq\frac{\bar{r}+m\left(\bar{b}^2+\bar{a}^2\right)}
{9\sqrt{m}\left(\bar{b}+\bar{a}\right)}$} 
and Eq.~\eqref{eq:bounded-dis-bound-dual-increase-batch-lower}, 
we have 
$K_2 - K_1\geq \left\lfloor
    \frac{W\sqrt{CT\ln{T}}}{6\sqrt{m}\left(\bar{b}+\bar{a}\right)}
\right\rfloor$, 
which implies that 
Eq.~\eqref{eq:bounded-dis-bounded-dual-batch-size-estimation} 
holds for each non-empty batch $i\in[N]$. 

When {\color{black}$\left\lfloor
    \frac{W\sqrt{CT\ln{T}}}{6\sqrt{m}\left(\bar{b}+\bar{a}\right)}
\right\rfloor\geq 1$}, 
we have 
$t_{i+1}-t_i\geq 
\frac{W\sqrt{CT\ln{T}}}
{12\sqrt{m}\left(\bar{b}+\bar{a}\right)}$.
Thus, under the event $\mathcal{E}$, we can derive that 
\[
\sum_{s=t_i}^{t_{i+1}-1}b_{v_is}
\geq 
\frac{12(t_{i+1}-t_i)\sqrt{m}\bar{r}}{W}.
\]
By Eq.~\eqref{eq:bounded-dis-bounded-dual-average-b-less}, 
we obtain
\begin{align}
     (K_2 - K_1)\bar{r}
     +\frac{m(K_2 - K_1)\left(\bar{b}^2+\bar{a}^2\right)}{\sqrt{CT\ln{T}}}
     -\frac{3W^2}{2}\sqrt{CT\ln{T}}
     &\geq
     \frac{W(K_2 - K_1 - 1-1)}{2\sqrt{m}}
     \cdot\frac{4\sqrt{m}\bar{r}}{W}
     \geq 
     (K_2 - K_1)\bar{r}, \notag
\end{align}
which implies that
\begin{align}
    K_2 - K_1 \geq 
    \frac{3W^2}{2m\left(\bar{b}^2+\bar{a}^2\right)}CT\ln{T}. \notag
\end{align}
This contradicts the fact $K_2 - K_1 \leq T$ 
and the condition {\color{black} $\frac{3W^2}{2m\left(\bar{b}^2+\bar{a}^2\right)}\geq 1$}.
\hfill \Halmos

\section{Omitted Proofs in Section~\ref{sec:finite-dis-log-no-de}}
\subsection{Proofs of Lemma~\ref{lem:finite-support-non-degenerate-unique} and Lemma~\ref{lem:finite-support-non-degeneracy-linear-solution}}
\label{sec:pf-finite-support-non-degenerate-unique}
We re-write LP~\eqref{eq:finite-support-fluid-lp} in the following general LP formulation:
\begin{align}
    &\text{maximize} &&~\bm{\rho}^{\top}\bm{\chi} \label{eq:finite-dis-non-degenerate-general} \\
    &\text{subject to} &&~\Gamma \bm{\chi}=\bm{\beta}  \qquad \text{and} \qquad \bm{\chi}\geq \bm{0},&\notag
\end{align}
where $\bm{\chi}\in\mathbb{R}^{n+2m}$ consists of all decision variables in LP~\eqref{eq:finite-support-fluid-lp}, $\Gamma=(\bm{\gamma}_1,\dots$,$\bm{\gamma}_{n+2m})\in\mathbb{R}^{(n+m)\times (n+2m)}$ and  $\bm{\beta}\in\mathbb{R}^{n+m}_{\geq 0}$ encode the equality constraints in LP~\eqref{eq:finite-support-fluid-lp}. By Definition~\ref{def:finite-support-non-degenerate}, let $\bm{\chi}^*$ denote the unique optimal solution of LP~\eqref{eq:finite-dis-non-degenerate-general}. Let $\mathcal{I}$ denote the set $\left\{ i\in[n+2m]: \chi^*_i>0 \right\}$, and by Definition~\ref{def:finite-support-non-degenerate} again, we have that $|\mathcal{I}|=n+m$. We further define $\Gamma_{\mathcal{I}}=(\bm{\gamma}_i)_{i\in\mathcal{I}}\in\mathbb{R}^{(n+m)\times (n+m)}$ and $\bm{\chi}^*_{\mathcal{I}}=(\chi_i)_{i\in\mathcal{I}}\in\mathbb{R}^{(n+m)}$. 

\proof{Proof of Lemma~\ref{lem:finite-support-non-degenerate-unique}.}
We first show that $\Gamma_{\mathcal{I}}$ is invertible. Suppose that there exists $\widetilde{\bm{\chi}}_{\mathcal{I}}\in\mathbb{R}^{n+m}$ such that $A_{\mathcal{I}}\widetilde{\bm{\chi}}_{\mathcal{I}}=0$ and $\widetilde{\bm{\chi}}_{\mathcal{I}}\neq\bm{0}$. We extend $\widetilde{\bm{\chi}}_{\mathcal{I}}\in\mathbb{R}^{n+m}$ to $\widetilde{\bm{\chi}}\in\mathbb{R}^{n+2m}$ by setting the elements at the $i$-th (for every $i\notin\mathcal{I}$) dimension to be $0$. We consider $\widehat{\bm{\chi}}\in\mathbb{R}^{n+2m}_{\geq 0}$ defined by 
\begin{equation}
    \widehat{\bm{\chi}}\defeq
    \begin{cases}
        \bm{\chi}+\epsilon\widetilde{\bm{\chi}}, \qquad &\text{if }\bm{\rho}^{\top}\widetilde{\bm{\chi}}\geq 0, \\[12pt]
        \bm{\chi}-\epsilon\widetilde{\bm{\chi}}, \qquad &\text{if }\bm{\rho}^{\top}\widetilde{\bm{\chi}}< 0,
    \end{cases}
\end{equation}
where $\epsilon$ is a positive and small constant such that $\widehat{\bm{\chi}}\geq \bm{0}$. Then we have that $\widehat{\bm{\chi}}$ is a feasible solution to LP~\eqref{eq:finite-dis-non-degenerate-general}, which is different from $\bm{\chi}^*$ and satisfies $\bm{\rho}^{\top}\widehat{\bm{\chi}}\geq \bm{\rho}^{\top}\bm{\chi}^*$. This contradicts to the optimality and uniqueness of $\bm{\chi}^*$.

We next consider adding a perturbation to linear program~\eqref{eq:finite-dis-non-degenerate-general}, i.e., 
\begin{align}
    &\text{maximize} &&~\bm{\rho}^{\top}\bm{\chi} \label{eq:finite-dis-non-degenerate-general-perturbation} \\
    &\text{subject to} &&~\Gamma \bm{\chi}=\bm{\beta} - \widetilde{\bm{\beta}}  \qquad \text{and} \qquad \bm{\chi}\geq \bm{0},&\notag
\end{align}
We show that when $2\Gamma_{\mathcal{I}}^{-1}\widetilde{\bm{\beta}}\leq \bm{\chi}^*_{\mathcal{I}}$ and $\bm{\beta}\geq\widetilde{\bm{\beta}}$, the optimal solution $\widetilde{\bm{\chi}}^*$ to LP~\eqref{eq:finite-dis-non-degenerate-general-perturbation} is unique and can be written in the following form: 
\begin{align}
    &(\widetilde{\chi}^*_i)_{i\in\mathcal{I}} = \bm{\chi}^*_{\mathcal{I}}-\Gamma_{\mathcal{I}}^{-1}\widetilde{\bm{\beta}} ,\qquad \qquad 
     (\widetilde{\chi}^*_i)_{i\notin\mathcal{I}}=\bm{0}. \label{eq:finite-dis-non-degenerate-exact-solution}
\end{align}
To prove the above statement, we may directly verify that $\widetilde{\bm{\chi}}^*$ defined in Eq.~\eqref{eq:finite-dis-non-degenerate-exact-solution} is a feasible solution to LP~\eqref{eq:finite-dis-non-degenerate-general-perturbation} and shares the same optimal basis for all $\widetilde{\bm{\beta}}$ whenever $2\Gamma_{\mathcal{I}}^{-1}\widetilde{\bm{\beta}}\leq \bm{\chi}^*_{\mathcal{I}}$. For the uniqueness of the optimal solution to LP~\eqref{eq:finite-dis-non-degenerate-general-perturbation}, consider any feasible solution $\widetilde{\bm{\chi}}'\in\mathbb{R}^{n+2m}_{\geq 0}$ such that $\bm{\rho}^{\top}\widetilde{\bm{\chi}}' \geq \bm{\rho}^{\top}\widetilde{\bm{\chi}}^*$, and let  $\bm{\chi}'=\bm{\chi}^*+\epsilon(\widetilde{\bm{\chi}}'-\widetilde{\bm{\chi}}^*)$, where $\epsilon=\min_{i\in[\mathcal{I}]}\frac{\chi^*_i}{\widetilde{\chi}^*_i}>0$. We then have $\bm{\chi}'\geq\bm{0}$ and $\Gamma\bm{\chi}'=\bm{\beta}$, indicating that $\widetilde{\bm{\chi}}'$ is a feasible solution to LP~\eqref{eq:finite-dis-non-degenerate-general-perturbation}. However, we have $\bm{\rho}^{\top}\bm{\chi}'=\bm{\rho}^{\top}\bm{\chi}^*+\epsilon\bm{\rho}^{\top}(\widetilde{\bm{\chi}}'-\widetilde{\bm{\chi}}^*)\geq \bm{\rho}^{\top}\bm{\chi}^*$, which, by the uniqueness of $\bm{\bm{\chi}}^*$, implies that $\bm{\chi}'=\bm{\chi}^*$. Thus, we have $\widetilde{\bm{\chi}}'=\widetilde{\bm{\chi}}^*$, i.e., the optimal solution to LP~\eqref{eq:finite-dis-non-degenerate-general-perturbation} is unique.

Finally, we set the constant $L=\frac{\min_{i\in\mathcal{I}}x^*_i}{2\|\Gamma_{\mathcal{I}}^{-1}\|_{\infty}}>0$ and verify that $\|\widetilde{\bm{b}}\|_{\infty}\leq\frac{\min_{i\in\mathcal{I}}x^*_i}{2\|\Gamma_{\mathcal{I}}^{-1}\|_{\infty}}$ implies $2\Gamma_{\mathcal{I}}^{-1}\widetilde{\bm{b}}\leq\bm{x}^*_{\mathcal{I}}$. This proves Lemma~\ref{lem:finite-support-non-degenerate-unique}.
\hfill\Halmos

\proof{Proof of Lemma~\ref{lem:finite-support-non-degeneracy-linear-solution}.} Consider $\|\widetilde{\bm{\beta}}^k\|_{\infty}\leq L$ ($k\in[K]$), $w_k\geq 0$ with $\sum_{k\in[K]}w_k=1$, and $\widetilde{\bm{\beta}}^0=\sum_{k\in[K]}w_k\widetilde{\bm{\beta}}^k$. Let $\widetilde{\bm{\chi}}^{*,k}$ ($k\in[K]\cup\{0\}$) denote the optimal solution to LP~\eqref{eq:finite-dis-non-degenerate-general-perturbation} with $\widetilde{\bm{\beta}}=\widetilde{\bm{\beta}}^k$. By Eq.~\eqref{eq:finite-dis-non-degenerate-exact-solution}, we have that
\begin{align}
    &(\widetilde{\chi}^{*,0}_i)_{i\in\mathcal{I}}=\bm{\chi}^*_{\mathcal{I}}-\Gamma_{\mathcal{I}}^{-1}\sum_{k\in[K]} w_k \widetilde{\bm{\beta}}^k=\sum_{k\in[K]}w_k (\widetilde{\chi}^{*,k}_i)_{i\in\mathcal{I}}, 
    &(\widetilde{\chi}^{*,0}_i)_{i\notin\mathcal{I}}=\bm{0}=\sum_{k\in[K]} w_k (\widetilde{\chi}^{*,k}_i)_{i\notin\mathcal{I}}. \notag
\end{align}
This proves Lemma~\ref{lem:finite-support-non-degeneracy-linear-solution}.
\hfill\Halmos

\subsection{Proof of Lemma~\ref{lem:finite-support-rare-event}}
\label{sec:pf-finite-support-rare-event}
By Hoeffding's inequality, we have
\begin{align}
     \Pr \left(\mathcal{E}_{\text{resource}}\right)&=1-\Pr\left(\exists z\in[K]\text{ and } j\in[m],~\left|  \sum_{t=(z-1)\kappa+1}^{z\kappa} \frac{b_{jt}}{\kappa}-B_j \right|> \widehat{L}\right) \notag \\
     &\geq 1-2Km\exp\left( -\frac{2\kappa \widehat{L}^2}{\bar{b}^2} \right)\geq 1-\frac{2m}{T^3}, \notag
\end{align}
where the last inequality is due to the condition that {\color{black} $\kappa\geq \frac{2\bar{b}^2}{\widehat{L}^2}\ln{T}$} and $K\leq T$.

By Hoeffding's inequality, we also have
\begin{align}
     \Pr \left(\mathcal{E}_{\text{request}}\right)&=1-\Pr\left(\exists z\in[K]\text{ and } i\in[n],~\left|  \sum_{t=(z-1)\kappa+1}^{z\kappa} \frac{\mathbb{I}(\theta_t=i)}{\kappa}-\mu_i \right|> \widehat{L}\right) \notag \\
     &\geq 1-2Kn\exp\left( -2\kappa \widehat{L}^2 \right)\geq 1-\frac{2n}{T^3}, \notag
\end{align}
where the last inequality is due to the condition that {\color{black} $\kappa\geq \frac{2}{\widehat{L}^2}\ln{T}$} and $K\leq T$.

By a union bound, we conclude that $\Pr\left( \mathcal{E} \right)\geq 1-\frac{2m+2n}{T^3} \geq 1 - \frac{4m+4n}{T^3}$.
\hfill\Halmos

\subsection{Proof of Lemma~\ref{lem:finite-dis-estimation-right-no-out-of-stock-1}}
\label{sec:pf-finite-dis-estimation-right-no-out-of-stock}
For any $z\geq W$, we prove that $x_t =\mathbb{I}\left(\theta_t\in\mathcal{I}^* ~\mathrm{or}~\Phi_{\theta_t,t}\geq 1\right)$ and $\Phi_{i,t}\leq \kappa$ hold for any $t\in[t_z+1,t_{z+1}]$ and $i\in[n]$ by induction on $z$. (Recall that we defined $t_z \defeq z \cdot \kappa$.)

\medskip
\noindent \underline{Induction basis.} When $z=W$, we have that for each $t\in[t_z+1,t_{z+1}]$ and $i\in[n]$, then
$\Phi_{i,t}\leq \Phi_{i,t_W+1} \leq \kappa$.
Moreover, under the event $\mathcal{E}$, we have that for each $j\in[m]$, it holds that
\begin{align}
    \ell_{j,t_W+1}= \sum_{k=0}^{W-1}\sum_{t=t_k+1}^{t_{k+1}} b_{jt}\geq W\kappa\left(B_j-\frac{1}{4}\underline{b}\right)\geq\frac{3W\kappa\underline{b}}{4}\geq 3n\max_{i\in[n],j\in[m]} |A_{ji}| \kappa,\label{eq:finite-dis-non-degenerate-initial}
\end{align}
where the last inequality is due to that {\color{black} $W\geq \frac{4n\max_{i\in[n],j\in[m]} |A_{ji}|}{\underline{b}}$.} We then derive that for any $j\in[m]$ and $t\in[t_W+1,t_W]$, we have
\begin{align}
    \ell_{j,t}&= \ell_{j,{t_W+1}}+\sum_{t=t_W+1}^{t-1} b_{jt} -\sum_{t=t_W+1}^{t-1} a_{jt}x_t  \geq  3\max_{i\in[n],j\in[m]} |A_{ji}| \kappa-\max_{i\in[n]} |A_{ji}| \kappa\geq \max_{i\in[n]}|A_{ji}| \geq a_{jt}, \notag 
\end{align}
where the first inequality by to Eq.~\eqref{eq:finite-dis-non-degenerate-initial}.  This proves the base case.

\medskip
\noindent\underline{Induction step.} Let $k > W$. Suppose $x_t =\mathbb{I}\left(\theta_t\in\mathcal{I}^* ~\mathrm{or}~\Phi_{\theta_t,t}\geq 1\right)$ and $\Phi_{i,t}\leq \kappa$ hold for any $t\in[t_W+1,t_k]$. We consider the period $t$ with $t\in[t_k+1,t_{k+1}]$.

We first prove that $\Phi_{i,t_{k}+1}\leq \kappa$. This is because
\begin{align}
    \Phi_{i,t}\leq \Phi_{i,t_k+1}&=\Phi_{i,t_{k-1}+1}+\Phi_{i,t_{k}+1}-\Phi_{i,t_{k-1}+1}\notag \\
    &=\Phi_{i,t_{k-1}+1}+X_{i,k+1}^*-\min\left( \Phi_{i,t_{k-1}+1}, (\widehat{\bm{\mu}}_{k+1})_{i} \right)\label{eq:finite-dis-no-degenerate-decomposition-for-i-not-in}\\
    &\leq \Phi_{i,t_{k-1}+1}+(\widehat{\bm{\mu}}_{k+1})_{i}-\min\left( \Phi_{i,t_{k-1}+1},(\widehat{\bm{\mu}}_{k+1})_{i} \right)=\max\left( \Phi_{i,t_{k-1}+1}, (\widehat{\bm{\mu}}_{k+1})_{i} \right)\leq \kappa,\notag
\end{align}
where the last inequality is due to the induction assumption that $\Phi_{i,t_{k-1}+1}\leq \kappa$ and $(\widehat{\bm{\mu}}_{k+1})_i\leq \kappa$.

We next prove that the stockout event never occurs. We have
\begin{align}
    \ell_{j,t_k+1}&=\ell_{j,t_W+1}+\sum_{z=W}^{k-1} \sum_{t'=t_z+1}^{t_{z+1}} b_{jt'}-\sum_{z=W}^{k-1} \sum_{t'=t_z+1}^{t_{z+1}}a_{jt}x_t \notag \\
    &= \ell_{j,t_W+1}+\sum_{z=W+2}^{k+1}\left(\sum_{i=1}^n A_{ji}X_{i,z}^* +S_{j,z}^*\right) - \sum_{z=W+2}^{k+1}\sum_{i\in\mathcal{I}^*}A_{ji}(\widehat{\bm{\mu}}_z)_i-\sum_{z=W}^{k-1}\sum_{i\notin\mathcal{I}^*}A_{ji}\min\left(\Phi_{i,t_z+1},(\widehat{\bm{\mu}}_{z+2})_{i}\right) \notag \\
    &\geq \ell_{j,t_W+1} + \sum_{i\notin\mathcal{I}^*}A_{ji}\left(\sum_{z=W}^{k-1} X_{i,z+2}^* - \sum_{z=W}^{k-1}\min\left(\Phi_{i,t_{z}+1}, (\widehat{\bm{\mu}}_{z+2})_{i} \right) \right),\notag \\
    &=\ell_{j,t_W+1} +\sum_{i\notin\mathcal{I}^*} A_{ji}\left(\Phi_{i,t_k+1}-\Phi_{i,t_W+1}\right)\geq \ell_{j,t_W+1}-n\max_{i\in[n]}|A_{ji}|\kappa\geq 2n\max_{i\in[n]}|A_{ji}|\kappa.
\end{align}
Here, the second equality is due to Line~\ref{line:finite-dis-no-degenerate-i-not-in-satisfy}, Line~\ref{line:finite-dis-no-degenerate-compute-optimal-batch} of the algorithm, and the induction assumption that stockout never happens before $t_k+1$; the first inequality is due to the condition that $X_{i,z}^*=(\widehat{\bm{\mu}}_{z})_{i}$ for $i\in\mathcal{I}^*$ under event $\mathcal{E}$; the third equality is due to Eq.~\eqref{eq:finite-dis-no-degenerate-decomposition-for-i-not-in}; the second inequality is due to that $\Phi_{i,t_W+1}\leq\kappa$ and $\Phi_{i,t_k+1}\leq \kappa$; the last inequality is due to Eq.~\eqref{eq:finite-dis-non-degenerate-initial}. Observing that 
\begin{align}
    \ell_{j,t}\geq \ell_{j,t_k+1}-\sum_{t'=t_k+1}^{t-1}a_{jt}x_t\geq 2n\max_{i\in[n]}|A_{ji}|\kappa-\max_{i\in[n]}|A_{ji}|\kappa\geq 
    \max_{i\in[n]}|A_{ji}|\geq a_{jt}, \notag
\end{align}
we prove that $x_t =\mathbb{I}\left(\theta_t\in\mathcal{I}^* ~\mathrm{or}~\Phi_{\theta_t,t}\geq 1\right)$ for all $t\in[t_k+1,t_{k+1}]$.
\hfill\Halmos

\section{Omitted Proofs in Section~\ref{sec:finite-hard-instance-de}}

\subsection{Proof of Theorem~\ref{thm:finite-dis-degenerate-hard-instance}}
\label{sec:pf-of-finite-dis-degenerate-hard-instance}

Without loss of generality, we prove the theorem for even $T$ and let $K = T/2$. Let $R^*$ denote the optimal hindsight reward. For clarity, we refer to the online setting in the theorem as the \emph{learning \& sequential decision} setting. Consider any policy $\pi$ in this setting, and let $R$ denote the (stochastic) reward achieved by $\pi$.

As an intermediate step toward proving the theorem, we also analyze a \emph{two-phase decision} setting defined as follows. The decision-maker fully knows the underlying distribution $\mathcal{P}_{\mathrm{finite}\text{-}\mathrm{hard}}$ at the beginning. Moreover, at time $1$ and time $K+1$, the decision-maker observes the types of the next $K$ orders, i.e., $\{\theta_t\}_{t\in[K]}$ and $\{\theta_{t}\}_{t\in \{K+1,\dots,2K\}}$, respectively. Under this setting, for each order type $i$, after observing the first block of $K$ orders, the optimal decision-maker may choose an integer $J_i$, reject the first $J_i$ arrivals of type $i$, and accept all subsequent type $i$ orders during time periods $1$ through $K$. Similarly, after observing the next $K$ orders at time $K+1$, the optimal decision-maker may choose integer $\widehat{J}_i$, reject the first $\widehat{J}_i$ arrivals of type $i$, and accept all subsequent type $i$ orders during time periods $K+1$ through $T=2K$.

\newcommand{\tp}{{\mathrm{two}\text{-}\mathrm{phase}}}
Let $\pi^*_\tp$ denote an optimal policy that achieves the maximum expected reward in the two-phase decision setting. Let $R^*_\tp$ denote its (stochastic) reward. It is immediate that $\mathbb{E} [R^*_\tp] \geq \mathbb{E} [R]$ since the decision-maker in the two-phase setting has strictly more information than in the learning \& sequential decision setting. Therefore, to prove the theorem, it suffices to show that $\mathbb{E}[R^* - R^*_\tp] \geq \Omega(\sqrt{T})$.

To prove the regret bound for $R^*_\tp$, we first analyze the constraints on $J_i$ and $\widehat{J}_i$ due to resource availability. Let $x_t \in \{0, 1\}$ denote the decision (whether accepting or rejecting the order) made by $\pi^*_\tp$ at time $t$. For each $t\in[K]$, we have the following constraints:
\begin{align}
    &\sum_{\tau=1}^t b_{1\tau}\geq \sum_{\tau=1}^t2\cdot x_t\mathbb{I}(\theta_\tau\in\{1,2,4\})\geq \sum_{\tau=1}^t2\cdot \mathbb{I}(\theta_\tau\in\{1,2,4\})-2(J_1+J_2+J_4),\notag \\
    &\sum_{\tau=1}^t b_{2\tau}\geq \sum_{\tau=1}^t2\cdot x_t\mathbb{I}(\theta_\tau\in\{1,3,5\})\geq \sum_{\tau=1}^t2\cdot \mathbb{I}(\theta_\tau\in\{1,3,5\})-2(J_1+J_3+J_5),\notag \\
    &\sum_{\tau=1}^{t+K} b_{1\tau}\geq \sum_{\tau=1}^{t+K} 2\cdot x_t\mathbb{I}(\theta_\tau\in\{1,2,4\})\geq \sum_{\tau=1}^{t+K} 2\cdot \mathbb{I}(\theta_i\in\{1,2,4\})-2(J_1+J_2+J_4+\widehat{J}_1+\widehat{J}_4+\widehat{J}_4),\notag \\
    &\sum_{\tau=1}^{t+K} b_{2\tau}\geq \sum_{\tau=1}^{t+K}2\cdot x_t\mathbb{I}(\theta_\tau\in\{1,3,5\})\geq \sum_{\tau=1}^{t+K}2\cdot \mathbb{I}(\theta_\tau\in\{1,3,5\})-2(J_1+J_3+J_5+\widehat{J}_1+\widehat{J}_3+\widehat{J}_5), \notag 
\end{align}
which indicates that 
\begin{align}
    &J_1+J_2+J_4\geq\frac{M_1}{2},\qquad &&J_1+J_3+J_5\geq\frac{M_2}{2},  \notag\\
    & J_1+J_2+J_4+ \widehat{J}_1+\widehat{J}_2+\widehat{J}_4 \geq \frac{M_1+\widehat{M}_1}{2},&&J_1+J_3+J_5+\widehat{J}_1+\widehat{J}_3+\widehat{J}_5\geq\frac{M_2+\widehat{M}_2}{2}, \label{eq:non-degerenate-lb-constraints-J}
\end{align}
where we define 
\begin{align}
    &M_1\defeq \max_{t\in[K]}\left(\sum_{\tau=1}^t2\cdot\mathbb{I}(\theta_\tau\in\{1,2,4\}) -t\right),\qquad
    &&M_2\defeq \max_{t\in[K]}\left(\sum_{\tau=1}^t2\cdot\mathbb{I}(\theta_\tau\in\{1,3,5\}) -t\right), \notag \\
    &\widehat{M}_1\defeq \max_{t\in[K]}\left(\sum_{\tau=K+1}^{2K}2\cdot\mathbb{I}(\theta_\tau\in\{1,2,4\}) -t\right),\qquad
    &&\widehat{M}_2\defeq \max_{t\in[K]}\left(\sum_{\tau=K+1}^{2K} 2\cdot\mathbb{I}(\theta_\tau\in\{1,3,5\}) -t\right). \notag
\end{align}
By the construction of  $\mathcal{P}_{\mathrm{finite}\text{-}\mathrm{hard}}$, for each $t$, we have $\frac{1}{4}=\Pr(\theta_t=1)=\Pr(\theta_t\in\{1,2,4\})\cdot\Pr(\theta_t\in\{1,3,5\})$. Thus, the two indicators $\mathbb{I}(\theta_t\in\{1,2,4\})$ and $\mathbb{I}(\theta_t\in\{1,3,5\})$ are independent, indicating that
\begin{itemize}
    \item $M_1$ and $M_2$ are independent;   $\widehat{M}_1$ and $\widehat{M}_2$ are independent; 
    \item $(M_1,M_2)$ and $(\widehat{M}_1,\widehat{M}_2)$ are independent.
\end{itemize}
We also define the total number of type $i$ orders by $N_i\defeq\sum_{t=1}^T\mathbb{I}(\theta_t=i)$ for each $i\in\{1,2,\dots,6\}$. We consider the following linear program:
\begin{align}
    &\text{maximize} &&\sum_{i=1}^5 R_i(N_i-\widetilde{J}_i'), \label{eq:finite-degenerate-offline-lp} \\
    &\text{subject to}    ~  && \widetilde{J}'_1+\widetilde{J}'_2+\widetilde{J}'_4\geq \frac{M_1+\widehat{M}_1}{2},&\qquad\widetilde{J}'_1 + \widetilde{J}'_3+\widetilde{J}'_5 \geq \frac{M_2+\widehat{M}_2}{2},\notag \\
    &&&\widetilde{J}'_i\geq 0,&\qquad \forall i\in \{1,2\dots, 5\}.\notag
\end{align}
We claim that the optimal solution to LP~\eqref{eq:finite-degenerate-offline-lp} upper bounds $R^*_\tp$. This is because we may simply set $\widetilde{J}'_i = J_i + \widehat{J}_i$ for each $i \in \{1, 2, \dots, 5\}$ to obtain a feasible solution to LP~\eqref{eq:finite-degenerate-offline-lp}, since $J_i$ and $\widehat{J_i}$ satisfy the constraints in Eq.~\eqref{eq:non-degerenate-lb-constraints-J}.

Our proof of the theorem then proceeds by two steps. First, we show that, with high probability, the difference between the hindsight optimum $R^*$ and the optimal solution to LP~\eqref{eq:finite-degenerate-offline-lp} is at most $\mathcal{O}\!\left(T^{3/8}\right)$. Specifically, let $\{\widetilde{J}^*_i\}$ denote an optimal solution to LP~\eqref{eq:finite-degenerate-offline-lp}. We have the following lemma, which will be proved in Section~\ref{sec:pf-of-finite-dis-hard-offline-opt-lower}.
\begin{lemma}
\label{lem:finite-dis-hard-offline-opt-lower}
    There exists a universal constant $T^{(1)}_{\mathrm{finite}\text{-}\mathrm{hard}} > 0$, such that when $T > T^{(1)}_{\mathrm{finite}\text{-}\mathrm{hard}}$, with probability $1-\frac{11}{T}$, we have
    $
        R^*\geq \sum_{i=1}^5 R_i\cdot(N_i-\widetilde{J}^*_i)-\mathcal{O}(T^{\frac{3}{8}})$.
\end{lemma}

Next, we establish that, with positive probability, there is a $\Omega(\sqrt{T})$ gap between $R^*_\tp$ and the optimal solution to  LP~\eqref{eq:finite-degenerate-offline-lp}. Specifically, we prove the following lemma in Section~\ref{sec:proof-lem-finite-dis-hard-two-phase-regret}.
\begin{lemma}
\label{lem:finite-dis-hard-two-phase-regret}
There exist universal constants $T^{(2)}_{\mathrm{finite}\text{-}\mathrm{hard}} > 0$ and $c_{\mathrm{finite}\text{-}\mathrm{hard}} > 0$, such that when $T > T^{(2)}_{\mathrm{finite}\text{-}\mathrm{hard}}$, with probability at least $c_{\mathrm{finite}\text{-}\mathrm{hard}} > 0$, it holds that $R^*_\tp \leq \sum_{i=1}^5 R_i\cdot(N_i-\widetilde{J}^*_i) - \sqrt{T}$.
\end{lemma}

Combining Lemma~\ref{lem:finite-dis-hard-offline-opt-lower} and Lemma~\ref{lem:finite-dis-hard-two-phase-regret}, we have that when $T> T_{\mathrm{finite}\text{-}\mathrm{hard}} \defeq \max\{T^{(1)}_{\mathrm{finite}\text{-}\mathrm{hard}}, T^{(2)}_{\mathrm{finite}\text{-}\mathrm{hard}}\}$, with probability at least $(c_{\mathrm{finite}\text{-}\mathrm{hard}} - 11/T) = \Omega(1)$, it holds that
\[
R^*_\tp \leq R^* - \sqrt{T} + \mathcal{O}({T^\frac38}) = R^* - \Omega(\sqrt{T}).
\]
Since $R^*_\tp \leq R^*$ almost surely, we further have $\mathbb{E}  [R^*_\tp] \leq \mathbb{E} [R^*] - \Omega(\sqrt{T})$, which concludes the proof of the theorem.
\hfill\Halmos

\subsection{Proof of Lemma~\ref{lem:finite-dis-hard-offline-opt-lower}}
\label{sec:pf-of-finite-dis-hard-offline-opt-lower}
We assume that the optimal solution $\{\widetilde{J}^*_i\}$ to LP~\eqref{eq:finite-degenerate-offline-lp} satisfies $\widetilde{J}^*_i\leq\max\left(M_1,M_2,\widehat{M}_1,\widehat{M}_2\right)$, 
since replacing each $\widetilde{J}^*_i$ with
$
\widetilde{J}^*_i\gets \min\left(\widetilde{J}_i^*,\max\left(M_1,M_2,\widehat{M}_1,\widehat{M}_2\right)\right)$
produces another feasible solution with no smaller objective value. Based on this optimal solution, we construct a hindsight (offline) solution using Algorithm~\ref{alg:finite-dis-hard-offline}, which achieves a total reward that is at most $\mathcal{O}(T^{3/8})$ below the optimal objective value of LP~\eqref{eq:finite-degenerate-offline-lp}.
Because the goal is to construct a hindsight solution, Algorithm~\ref{alg:finite-dis-hard-offline} has full knowledge of all realized replenishments and orders at every time period.
\begin{algorithm}[h]
\caption{Hindsight solution construction for the hard instance $\mathcal{P}_{\mathrm{finite}\text{-}\mathrm{hard}}$}
\label{alg:finite-dis-hard-offline}
\begin{algorithmic}[1]
\State Initialize $\Phi_{i,1}\gets\lceil\widetilde{J}^*_i\rceil$ for each $i\in[5]$ and $\ell_{j,1}\gets 0$ for each $j\in[2]$, $\kappa\gets \lceil T^{3/8}\rceil$.
\For{$t \gets 1$ to $T$} 
    \State Observe order type $\theta_t$ (and receive the deterministic replenishment).
    \If{$t\leq \kappa$} Set $x_t\gets 0$.
    \Else ~ Set $x_t\gets \mathbb{I}\left(\forall j\in[2],~\ell_{j,t}+1\geq A_{j\theta_t}\right)\cdot\mathbb{I}(\Phi_{i,t}\leq 0)$;
        update $\Phi_{\theta_t,t+1}\gets \Phi_{\theta_t,t}-(1-x_t)$.
    \EndIf
    \State Update $\ell_{j,t+1}\gets \ell_{j,t}+1-x_tA_{j\theta_t}$ for each $j\in[2]$.
\EndFor
\end{algorithmic}
\end{algorithm}

We define $W_i(a,b)$ for  each $i\in[2],1\leq a\leq b\leq T$ by
\begin{align}
    W_i(a,b)\defeq\max_{a\leq t\leq b}\left(\sum_{t=a}^b 2\cdot\mathbb{I}(\theta_t\in\{1,1+i,3+i\})-(b-a+1) \right). \notag
\end{align}
We can verify that $M_i=W_i(1,K)$ and $\widehat{M}_i=W_i(K+1,2K)$. Let $C_{\mathrm{finite}\text{-}\mathrm{hard}}=\frac{\lceil 4051\sqrt{T\ln{T}} \rceil}{\sqrt{T\ln{T}}}$, and consider the following events:
\begin{align}
    &\mathcal{U}_1=\left\{\forall i\in[2],~ W_i\left( \kappa+1,C_{\mathrm{finite}\text{-}\mathrm{hard}}\sqrt{T\ln{T}}  \right) \leq \kappa \right\},\notag \\
    & \mathcal{U}_2=\left\{ \forall i\in[5],\sum_{t=\kappa+1}^{C_{\mathrm{finite}\text{-}\mathrm{hard}}\sqrt{T\ln{T}}} \mathbb{I}(\theta_t=i)\geq\frac{C_{\mathrm{finite}\text{-}\mathrm{hard}}\sqrt{T\ln{T}}}{10} \right\},  \notag \\
    &\mathcal{U}_3=\left\{\forall i\in[2],~ \max\left(W_i(1,K),W_i(K+1,2K) \right)\leq\frac{C_{\mathrm{finite}\text{-}\mathrm{hard}}\sqrt{T\ln{T}}}{10} \right\}.\notag
\end{align}
We have the following lemma concerning $\Pr(\mathcal{U}_1\cap\mathcal{U}_2\cap\mathcal{U}_3)$, whose proof is deferred to Section~\ref{sec:pf-finite-dis-hard-some-event}.
\begin{lemma}
\label{lem:finite-dis-hard-some-event}
    For {\color{black} $C_{\mathrm{finite}\text{-}\mathrm{hard}}\geq 4050$} and {\color{black}$T$ satisfying $T^{\frac{1}{4}}\geq 4C_{\mathrm{finite}\text{-}\mathrm{hard}}(\ln{T})^{\frac{3}{2}}$}, we have
    $
        \Pr(\mathcal{U}_1\cap\mathcal{U}_2\cap\mathcal{U}_3)\geq 1-\frac{11}{T}$.
\end{lemma}
Under the event $\mathcal{U}_1\cap\mathcal{U}_2\cap\mathcal{U}_3$, we prove that we have $x_t=x'_t\defeq \mathbb{I}(\Phi_{i,t}\leq 0)$ for all $t\geq \kappa+1$ and $x_t=x'_t\defeq 0$ for all $t\leq \kappa$. For $t\in[\kappa+1,C_{\mathrm{finite}\text{-}\mathrm{hard}}\sqrt{T\ln{T}}]$ and $j\in[2]$, we have 
\begin{align}
    \sum_{i=1}^t b_{ji}- 2\sum_{i=1}^t x'_i\cdot \mathbb{I}\left( \theta_i\in\{1,1+j,3+j\} \right)&  \geq t-2\sum_{i=\kappa+1}^t \mathbb{I}\left( \theta_i\in\{1,1+j,3+j\} \right),\notag \\
    &\geq \kappa-W_{j}\left(\kappa+1,C_{\mathrm{finite}\text{-}\mathrm{hard}}\sqrt{T\ln{T}}\right)\geq 0, \notag
\end{align}
which indicates that stockouts do not occur, i.e., $x_t=x'_t$, when $t\leq C_{\mathrm{finite}\text{-}\mathrm{hard}}\sqrt{T\ln{T}}$. 

For $t>C_{\mathrm{finite}\text{-}\mathrm{hard}}\sqrt{T\ln{T}}$ and $j\in[2]$, we have 
\begin{align}
    &\sum_{i=1}^t b_{ji} - 2\sum_{i=1}^t x'_i\cdot \mathbb{I}\left( \theta_i\in\{1,1+j,3+j\} \right),\notag \\
    & \geq t-2\sum_{i=1}^t\mathbb{I}(\theta_i\in\{1,1+j,3+j\}) + 2\sum_{i=\kappa+1}^{C_{\mathrm{finite}\text{-}\mathrm{hard}}\sqrt{T\ln{T}}}\mathbb{I}(\theta_i\in\{1,1+j,3+j\})\cdot\mathbb{I}(\Phi_{\theta_i,t}> 0),\notag \\
    &\geq -W_{j}(1,T) + \sum_{v\in\{1,1+j,3+j\}} 2 \min\left(\left\lceil\widetilde{J}_v^*\right\rceil,\sum_{i=\kappa+1}^{C_{\mathrm{finite}\text{-}\mathrm{hard}}\sqrt{T\ln{T}}}\mathbb{I}(\theta_i=v)\right)=2\sum_{v\in\{1,1+j,3+j\}}\left\lceil\widetilde{J}_v^*\right\rceil-W_{j}(1,T),\label{eq:finite-dis-hard-after-no-out-of-stock}
\end{align}
where the second inequality is due to  $x_i=\mathbb{I}(\Phi_{\theta_i}^i\leq 0)$ and $\Phi_{\theta_i}^{i+1}=\Phi_{\theta_i}^i-1+x_i$, and the last equality is due to  
\begin{align*}
\widetilde{J}_v^*\leq \max(W_1(1,K),W_1(K+1,2K),W_2(1,K),W_2(K+1,2K))&\leq\frac{C_{\mathrm{finite}\text{-}\mathrm{hard}}\sqrt{T\ln{T}}}{10}\\
&\leq\sum_{i=\kappa+1}^{C_{\mathrm{finite}\text{-}\mathrm{hard}}\sqrt{T\ln{T}}}\mathbb{I}(\theta_i=v).
\end{align*}
We verify that $W_j(1,T)\leq W_j(1,K)+W_j(K+1,2K)$. By Eq.~\eqref{eq:finite-dis-hard-after-no-out-of-stock} and $\widetilde{J}_v^*$ is a feasible solution to LP~\eqref{eq:finite-degenerate-offline-lp}, for $t\in[\kappa+1,C_{\mathrm{finite}\text{-}\mathrm{hard}}\sqrt{T\ln{T}}]$ and $j\in[2]$, we have that
\[
\sum_{i=1}^t b_{ji} - 2\sum_{i=1}^t x'_i\cdot \mathbb{I}\left( \theta_i\in\{1,1+j,3+j\}\right)\geq 0,
\]
indicating that stockouts do not happen, i.e., $x_i=x'_i$ when $t> C\sqrt{T\ln{T}}$.

Finally, we compute the reward of Algorithm~\ref{alg:finite-dis-hard-offline}. Under the event $\mathcal{U}_1\cap\mathcal{U}_2\cap\mathcal{U}_3$, we have
\begin{align}
    R^*\geq \sum_{v\in[5]}R_v(N_v-\widetilde{J}_v^*-1-\kappa)=\sum_{v\in[5]}R_v(N_v-\widetilde{J}_v^*)-\mathcal{O}(T^{\frac{3}{8}}). \notag \Halmos
\end{align}

\subsection{Proof of Lemma~\ref{lem:finite-dis-hard-some-event}}
\label{sec:pf-finite-dis-hard-some-event}
We prove the following lower bounds for the probability of each of the events $\mathcal{U}_1$, $\mathcal{U}_2$, and $\mathcal{U}_3$. Then, the lemma follows by a union bound.

\medskip
\noindent\underline{$\Pr(\mathcal{U}_1)\geq 1-\frac{2}{T}$}. It suffices to show that for each $t\in\left[\kappa+1, C_{\mathrm{finite}\text{-}\mathrm{hard}}\sqrt{T\ln{T}}\right],i\in[2]$, we have
\begin{align}
    \Pr\left( \sum_{v=\kappa+1}^{t} 2\cdot\mathbb{I}(\theta_v\in\{1,1+i,3+i\})-t+\kappa > \kappa \right)\leq \frac{1}{T^2}. \notag
\end{align}
By Hoeffding's inequality, we have 
\begin{align}
     \Pr\left( \sum_{v=\kappa+1}^{t} 2\cdot\mathbb{I}(\theta_v\in\{1,1+i,3+i\})-t+\kappa > \kappa \right)&\leq \exp\left( -\frac{2(\kappa)^2}{4t} \right),\notag \\
     & \leq \exp\left(-\frac{T^{\frac{3}{4}}}{2C_{\mathrm{finite}\text{-}\mathrm{hard}}\sqrt{T\ln{T}}}\right) \leq \frac{1}{T^2}, \notag
\end{align}
where the last inequality is due to {\color{black}$T^{\frac{1}{4}}\geq 4C_{\mathrm{finite}\text{-}\mathrm{hard}}\left(\ln{T}\right)^{\frac{3}{2}}$}.

\medskip
\noindent\underline{$\Pr\left( \mathcal{U}_2 \right)\geq 1-\frac{5}{T}$}. By Hoeffding's inequality, we have 
\begin{align}
    &\Pr\left( \sum_{t=\kappa+1}^{C_{\mathrm{finite}\text{-}\mathrm{hard}}\sqrt{T\ln{T}}}\mathbb{I}(\theta_t=i)-\frac{C_{\mathrm{finite}\text{-}\mathrm{hard}}\sqrt{T\ln{T}}}{9}\leq -\frac{C_{\mathrm{finite}\text{-}\mathrm{hard}}\sqrt{T\ln{T}}}{90} \right)
   \notag\\
   &\qquad\qquad\qquad\qquad\qquad\qquad\qquad\qquad\qquad\qquad\leq \exp\left(-\frac{2C_{\mathrm{finite}\text{-}\mathrm{hard}}^2T\ln{T}}{8100C_{\mathrm{finite}\text{-}\mathrm{hard}}\sqrt{T\ln{T}}}\right)\leq\frac{1}{T}, \notag
\end{align}
where the first inequality is due to {\color{black} $\kappa\leq T^{\frac{3}{8}}+1\leq \frac{C_{\mathrm{finite}\text{-}\mathrm{hard}}\sqrt{T\ln{T}}}{72}$}, and the second inequality is due to {\color{black} $C_{\mathrm{finite}\text{-}\mathrm{hard}}\geq 4050$} and $T\geq\ln{T}$.

\medskip
\noindent\underline{$\Pr(\mathcal{U}_3)\geq 1-\frac{4}{T}$}. It suffices to show that for each $i\in[2],j\in\{0,1\}$, and $t\in[K]$, we have
\begin{align}
    \Pr\left( \sum_{v=1+jK}^{t+jK} 2\cdot\mathbb{I}(\theta_v\in\{1,1+i,3+i\}) - t>\frac{C_{\mathrm{finite}\text{-}\mathrm{hard}}\sqrt{T\ln T}}{10} \right)\leq \frac{1}{T^2}. \notag
\end{align}
By Hoeffding's inequality, we have
\begin{align}
     \Pr\left( \sum_{v=1+jK}^{t+jK} 2\cdot\mathbb{I}(\theta_v\in\{1,1+i,3+i\}) - t>\frac{C_{\mathrm{finite}\text{-}\mathrm{hard}}\sqrt{T\ln T}}{10} \right) &\leq\exp\left( -\frac{2C_{\mathrm{finite}\text{-}\mathrm{hard}}\cdot T\ln{T}}{400t} \right)\leq\frac{1}{T^2}, \notag
\end{align}
where the last inequality is due to $C_{\mathrm{finite}\text{-}\mathrm{hard}}\geq 2000$.
\hfill\Halmos

\subsection{Proof of Lemma~\ref{lem:finite-dis-hard-two-phase-regret}}
\label{sec:proof-lem-finite-dis-hard-two-phase-regret}
In this proof, we will use the following technical lemma, the proof of which can be found in Section~8 of \citet{billingsley2013convergence}.
\begin{lemma}
\label{lem:finite-max-partial-sum-abs-N01}
Let $\{X_i\}_{i=1}^\infty$ be a sequence of \emph{i.i.d.}~random variables with zero mean $ \mathbb{E}[X_1]=0$ and finite variance $\mathrm{Var}(X_1)=1$. Consider the maximum of the first $T$ centered partial sums $M_T\defeq\max_{1 \leq k \leq T}  \sum_{i=1}^k X_i$, then, as $T \to \infty$,
\[
\frac{M_T}{\sqrt{T}} \xrightarrow{d} |Z|,
\qquad \text{where } Z \sim \mathcal{N}(0,1).
\]
In other words, the suitably normalized maximum of the partial sums converges in distribution to the absolute value of a standard normal random variable.
\end{lemma}

We now start to prove Lemma~\ref{lem:finite-dis-hard-two-phase-regret}. Let $J_i$ and $\widehat{J}_i$ be the number of rejections made by $\pi^*_\tp$, we have that the reward achieved by $\pi^*_\tp$ is 
$
    R^*_\tp=\sum_{i=1}^5 R_i(N_i-J_i-\widehat{J}_i)$.
We consider the following three events:
\begin{align}
    &\mathcal{E}_1\defeq \left\{ M_1\in[30\sqrt{T},40\sqrt{T}] \right\}\bigcap \left\{ M_2\in[\sqrt{T},2\sqrt{T}] \right\}, \notag \\
    &\mathcal{E}_2\defeq \left\{ \widehat{M}_1,\widehat{M}_2\in[0.1\sqrt{T},0.2\sqrt{T}] \right\},\notag \\
    &\mathcal{E}_3\defeq \left\{ \widehat{M}_1\in[\sqrt{T},2\sqrt{T}] \right\}\bigcap \left\{ \widehat{M}_2\in[120\sqrt{T},130\sqrt{T}] \right\}. \notag
\end{align}
By Lemma~\ref{lem:finite-max-partial-sum-abs-N01}, we have that there exist absolute constants $T_{\min}$ and $c_{\text{min}}>0$ such that for {\color{black} $T>T_{\min}$},
\begin{align}
    \Pr\left( \mathcal{E}_i \right)>c_{\text{min}},\qquad \forall i\in \{1,2,3\}. \notag
\end{align}
Now we set $T^{(2)}_{\mathrm{finite}\text{-}\mathrm{hard}} = T_{\min}$ and $c_{\mathrm{finite}\text{-}\mathrm{hard}} = c_{\min}^2 / 2$. When $T > T^{(2)}_{\mathrm{finite}\text{-}\mathrm{hard}}$, we consider $\Pr\left( J_1\geq 10\sqrt{T}~\middle|~ \mathcal{E}_1 \right)$ and discuss the following two cases.

\medskip
\noindent\underline{Case 1: $\Pr\left( J_1\geq 10\sqrt{T}~\middle|~ \mathcal{E}_1 \right)\geq \frac{1}{2}$}. Notice that $(J_1,\mathcal{E}_1)$ is independent of $\mathcal{E}_2$. Thus, we have 
\begin{align}
    \Pr\left(\left\{J_1\geq 10\sqrt{T}\right\}\bigcap\mathcal{E}_1\bigcap\mathcal{E}_2\right)= \Pr\left(\left\{J_1\geq 10\sqrt{T}\right\}\bigcap\mathcal{E}_1\right)\cdot\Pr\left( \mathcal{E}_2 \right)\geq \frac{c_{\text{min}}^2}{2} = c_{\mathrm{finite}\text{-}\mathrm{hard}}. \notag
\end{align}
Under the event $\left\{J_1\geq 10\sqrt{T}\right\}\bigcap\mathcal{E}_1\bigcap\mathcal{E}_2$, we consider a feasible solution $\widetilde{J}'_i$ to LP~\eqref{eq:finite-degenerate-offline-lp} defined by
\begin{align}
    \widetilde{J}'_1=J_1+\widehat{J}_1-\sqrt{T};\qquad \widetilde{J}'_2=J_2+\widehat{J}_2+\sqrt{T};\qquad \widetilde{J}'_{i}=J_{v}+\widehat{J}_i \text{ for each }i\in\{3,4,5\}. \notag
\end{align}
We verify that 
\begin{align}
\sum_{i=1}^5 R_i \left(N_i - \widetilde{J}_i^*\right) - R^*_{\tp}  \geq \sum_{i=1}^5 R_i \left(N_i - \widetilde{J}_i'\right) - R^*_{\tp} =
    \sum_{i=1}^5 R_i\left( J_i+\widehat{J}_i-\widetilde{J}'_i \right) =2\sqrt{T} . \notag
\end{align}

\medskip
\noindent\underline{Case 2: $\Pr\left( J_1\geq 10\sqrt{T}~\middle|~ \mathcal{E}_1 \right)< \frac{1}{2}$}. In this case, we have $\Pr\left( J_1< 10\sqrt{T}~\middle|~ \mathcal{E}_1 \right)> \frac{1}{2}$. Notice that $(J_1,\mathcal{E}_1)$ is independent of $\mathcal{E}_3$. Thus, we have 
\begin{align}
    \Pr\left(\left\{J_1 < 10\sqrt{T}\right\}\bigcap\mathcal{E}_1\bigcap\mathcal{E}_3\right)= \Pr\left(\left\{J_1 < 10\sqrt{T}\right\}\bigcap\mathcal{E}_1\right)\cdot\Pr\left( \mathcal{E}_3 \right)\geq \frac{c_{\text{min}}^2}{2} = c_{\mathrm{finite}\text{-}\mathrm{hard}}. \notag
\end{align}
Under the event $\left\{J_1 < 10\sqrt{T}\right\}\bigcap\mathcal{E}_1\bigcap\mathcal{E}_3$, we have 
$     J_2 + J_4 \geq 5\sqrt{T}$.
Next we consider the following two sub-cases:
\begin{itemize}
    \item \underline{Case 2a: $J_3+\widehat{J}_3+J_5+\widehat{J}_5\leq \sqrt{T}$.} We have $J_1+\widehat{J}_1\geq 50\sqrt{T}$. Consider a feasible solution $\widetilde{J}'_v$ to LP~\eqref{eq:finite-degenerate-offline-lp} defined by
\begin{align}
    \widetilde{J}'_1=J_1+\widehat{J}_1;\qquad \widetilde{J}'_2=\widetilde{J}'_4=0;\qquad \widetilde{J}'_{i}=J_{i}+\widehat{J}_v \text{ for }i\in\{3,5\}. \notag
\end{align}
We verify that
\begin{align}
    \sum_{i=1}^5 R_i \left(N_i - \widetilde{J}_i^*\right) - R^*_{\tp}  \geq \sum_{i=1}^5 R_i \left(N_i - \widetilde{J}_i'\right) - R^*_{\tp} =\sum_{i=1}^5 R_i\left( J_i+\widehat{J}_i-\widetilde{J}'_i \right)\geq 5\sqrt{T}. \notag
\end{align}
\item \underline{Case 2b: $J_3+\widehat{J}_3+J_5+\widehat{J}_5> \sqrt{T}$.} Consider a feasible solution $\widetilde{J}'_v$ to LP~\eqref{eq:finite-degenerate-offline-lp} defined by
\begin{align}
    \widetilde{J}'_1=J_1+\widehat{J}_1+\sqrt{T};\qquad &\widetilde{J}'_2=J_2+\widehat{J}_2-\min(J_2+\widehat{J}_2,\sqrt{T});\qquad \widetilde{J}'_4=J_4+\widehat{J}_4-\left[ \sqrt{T}-J_2-\widehat{J}_2 \right]^+ ;\notag \\
    &\widetilde{J}'_3=J_3+\widehat{J}_3-\min(J_3+\widehat{J}_3,\sqrt{T});\qquad \widetilde{J}'_5=J_5+\widehat{J}_5-\left[ \sqrt{T}-J_3-\widehat{J}_3 \right]^+ .\notag
\end{align}
We verify that
\begin{align}
    \sum_{i=1}^5 R_i \left(N_i - \widetilde{J}_i^*\right) - R^*_{\tp}  \geq \sum_{i=1}^5 R_i \left(N_i - \widetilde{J}_i'\right) - R^*_{\tp} =\sum_{i=1}^5 R_i\left( J_i+\widehat{J}_i-\widetilde{J}'_i \right)\geq\sqrt{T}. \notag
\end{align}
\end{itemize}
In all, we prove the lemma.
\hfill\Halmos

\section{Omitted Proofs and Discussions related to Section~\ref{sec:non-degenerate-dis}}

\subsection{The Definition of Constants and an Assumption on $T$}
\label{sec:const-def}
We first define some constants which will be useful in Algorithm~\ref{alg:replenish-non-degenerate-dis-cts-saving}, Algorithm~\ref{alg:replenish-non-degenerate-dis-cts-consumption}, Algorithm~\ref{alg:replenish-non-degenerate-dis-cts-combine}, and their performance analysis. The constant $N$ in Eq.~\eqref{eq:cts-def-c-5} is defined in Eq.~\eqref{eq:cts-def-q-N}.
\begin{align}
    & C_5 = \max\left( \frac{  2\bar{a} + 1 + \sqrt{\left(2\bar{a}+1\right)^2  + \frac{\lambda\lambda_{\text{min}}}{8}}   }{ \lambda\lambda_{\text{min}}/ 16  }\cdot \sqrt{ 2m\ln(2N) + 10 +10m\bar{a} } , ~\frac{\sqrt{5m}\cdot\bar{b}}{\lambda\lambda_{\text{min}}} \right)\label{eq:cts-def-c-5} \\
     & C_9 = \max\left( 16\sqrt{2}\,\bar{b}+8\sqrt{2}C_2\bar{a}^2\mu , 9\bar{b}+18\bar{a}\right) \\
    & C_7 =   2(\bar{b}+2\bar{a}) \\
    & C_2 =  \max \left( 32(\bar{b}+2\bar{a}), 4C_5\mu\bar{a}^2+12\bar{b}  \right)\\
    & C_3 =  \max\left(  2\sqrt{3}C_5\mu\bar{a}^2, 3(\bar{b}+2\bar{a})  \right) \\
    & C_4 = \max \left( 3 C_9,  16 C_5\mu\bar{a}^2 \right) \\
    & C_1 = \frac{4(C_2)^2}{\delta_B^2} \\
    & C_6 = \max \left(  16C_5, \frac{8}{\lambda\lambda_{\text{min}}}\cdot\sqrt{2m\cdot ((C_2)^2 +5\bar{b}^2 )} \right) \\
    & C_0 = \max\left(  \frac{10\bar{b}^2}{\delta_B^2} ,~ \frac{20\bar{a}}{\lambda_{\text{min}}}+10 ,~ 2C_1 , \frac{1}{\underline{b}}\sqrt{C_1}\cdot C_7,  ~\frac{10\bar{b}^2}{\underline{b}^2}  \right) \\
    & C_8 = 1 + \max\left( 18C_2, C_4, 5\sqrt{6}\bar{b} C_0  \right) \\
    & C_{10} = \sqrt{  6(C_5)^2 + \frac{m}{\lambda^2\lambda^2_{\text{min}}}\cdot (20\bar{b}^2 + 18(C_8+C_9)^2) } 
\end{align}
We define the constant $C_{11}$ to be the minimal value such that for any $H\geq C_{11}$, we have 
\begin{align}
    &H> 3 + 2C_0 \ln^2 H, \label{eq:no-de-minH-1} \\
    &H > 22 +\frac{\ln H}{\delta_B^2} \cdot \max\left( 20(\bar{b})^2 , \frac{40\bar{a}\delta_B^2}{\lambda_{\text{min}}},63(\bar{b}+2\bar{a})^2,400(C_3)^2 ,  24 (C_5)^2\bar{a}^4\mu^2+2\delta_B^2,10\bar{b}^2+2\delta_B^2  \right). \label{eq:no-de-minH-2} 
\end{align}
In fact, for any positive constant $a,b>3$, we can directly check that
\begin{align}
    H>2a + 8b ^3 + 125   \Rightarrow H>a +b\ln^2 H>a+b\ln H. 
\end{align}
Therefore, we can check that $C_{11}$ is $\mathcal{O}(1)$. 

In the proof of Theorem~\ref{thm:regret-non-degenerate-cts-dis-1}, we assume that $T\geq 2C_{11}$. In fact, if Eq.~\eqref{eq:no-de-minH-1} or Eq.~\eqref{eq:no-de-minH-2} doesn't hold, we have 
\begin{align}
    &\Reg_{\mathcal{P}}(\bm{\pi}_{\mathrm{non}\text{-}\mathrm{deg}})\leq 2\bar{r}T,\notag \\
    &\leq 44+4C_0 \ln^2 T+\frac{2\ln T}{\delta_B^2} \cdot \max\left( 20(\bar{b})^2 , \frac{40\bar{a}\delta_B^2}{\lambda_{\text{min}}},63(\bar{b}+2\bar{a})^2,400(C_3)^2 ,  24 (C_5)^2\bar{a}^4\mu^2+2\delta_B^2,10\bar{b}^2+2\delta_B^2  \right),\notag \\
    &=\mathcal{O}(\ln^2 T).\notag
\end{align}

\subsection{Technical Lemmas from \citet{li2022online}}
We have the following lemma concerning the constant $\delta_B$ defined in Lemma~\ref{lem:cts-no-de-d-robust}. This statement is also briefly mentioned in \citet{li2022online} (in the proof of their Theorem 5). For completeness, we provide its proof in Section~\ref{sec:pf-cts-no-de-many-inv-for-no-binding}. 
\begin{lemma}
\label{lem:cts-no-de-many-inv-for-no-binding}
    For any non-binding resource $j\not\in I_B$, we have
    \[
        B_j\geq \delta_B + \mathbb{E}_{(r,\bm{a},\bm{b})\sim\mathcal{P}}\left[ a_{j}\cdot\mathbb{I}\left(r> \left\langle \bm{a} , \bp^* \right\rangle\right) \right].
    \]
\end{lemma}
We next introduce the property of dual‐variable convergence, stated in Lemma~\ref{lem:cts-non-degeneracy-dual-convergence}. This lemma strengthens the statement at the beginning of the proof of their Theorem 1 in \cite{li2022online}. The statement in \cite{li2022online} considers a specific vector $\bm{B}$ rather than establishing a uniform bound for all $\bm{B}' \in \Omega_b$. However, extending their statement to a uniform bound as in Lemma~\ref{lem:cts-non-degeneracy-dual-convergence} is straightforward because the proof in \cite{li2022online} considers a random variable that upper bounds $\left\| \widehat{\bp}^* - \bp^*(\bm{B}') \right\|_2$, and this random variable depends only on $\{(r_v, \bm{a}_v, \bm{b}_v)\}_{v \in [t]}$, rather than the specific choice of $\bm{B}'$. 
\begin{lemma}
\label{lem:cts-non-degeneracy-dual-convergence}
   Consider $t\in\mathbb{N}_+$, and $\epsilon>0$. We denote $\mathcal{S}=\{(r_v,\bm{a}_v,\bm{b}_v)\}_{v\in[t]}$ as \emph{i.i.d.}~samples from the non-degenerate distribution $\mathcal{P}$. Let 
   \begin{align}
        N=\left\lfloor \log_{q}\left( \frac{\underline{b}\epsilon^2}{\bar{a}\bar{r}\sqrt{m}} \right) \right\rfloor + 1,\qquad \text{where}\qquad q=\min\left\{ {1+\frac{1}{\sqrt{m}}},{1+\frac{1}{\sqrt{m}}\left( \frac{\lambda\lambda_{\min}}{8\mu\bar{a}^2} \right)^{\frac{1}{3}}} \right\}^{-1} . \label{eq:cts-def-q-N}
    \end{align}
   With probability at least 
   $
        1-2m\exp\left(-\frac{t\epsilon}{2\bar{a}^2 m}\right) - m\exp\left( -\frac{t\lambda_{\min}}{4\bar{a}^2} \right)-2(2N)^m\cdot\exp\left( -\frac{t\epsilon}{2} \right)
    $,  for any $\bm{B}'\in\times_{j=1}^m\left( B_j-\delta_B , B_j + \delta_B \right)\subset\Omega_b$, and any $\widehat{\bp}^* \in\arg\min_{\bp\geq\bm{0}}f\left(\bp;\bm{B}',\mathcal{S}\right)$, it holds that 
    \begin{align}
        \left\|\widehat{\bp}^* -\bp^*(\bm{B}')\right\|_2^2 \leq \left( \frac{2\bar{a}+1+\sqrt{\left(2\bar{a}+1\right)^2+\frac{\lambda\lambda_{{\min}}}{8}}}{\lambda\lambda_{{\min}}/16} \right)^2\cdot\epsilon . \notag
    \end{align}
\end{lemma}
Lemma~\ref{lem:cts-non-degeneracy-dual-convergence} directly implies the following corollary.
\begin{corollary}
\label{cor:cts-no-de-dual-conver-simple}
    Assume the same conditions up to Eq.~\eqref{eq:cts-def-q-N} in Lemma~\ref{lem:cts-non-degeneracy-dual-convergence}. Consider any $T$ such that $T\geq t\geq \frac{20\bar{a}}{\lambda_{{\min}}}\cdot\ln T+10$. Let the parameter $C_5$ satisfy 
    \begin{align}
        {\color{black} C_5\geq\frac{  2\bar{a} + 1 + \sqrt{\left(2\bar{a}+1\right)^2  + \frac{\lambda\lambda_{{\min}}}{8}}   }{ \lambda\lambda_{{\min}}/ 16  }\cdot \sqrt{ 2m\ln(2N) + 10 +10m\bar{a} }}. \notag
    \end{align}
    With probability at least $1-\frac{3m+2}{T^5}$, for each $\bm{B}'\in\times_{j=1}^m (B_j-\delta_B,B_j+\delta_B)\subset\Omega_b$, it holds that
    \begin{align}
        \| \widehat{\bp}_t^*(\bm{B}') - \bp^*(\bm{B}') \|_2^2 \leq (C_5)^2\cdot\frac{\ln{T}}{t}. \notag
    \end{align}
    
\end{corollary}

\subsection{Proof of Remark~\ref{rmk:cts-no-de-auto-hold-equal}}
\label{sec:proof-rmk-cts-no-de-auto-hold-equal}
For any $\bm{B}'\in\Omega_b$, we define $\bm{\Phi}\defeq\nabla f(\bm{B}',\bp^*(\bm{B}'))$. Consider any $j\in[m]$ such that $\bp_j^*(\bm{B}')>0$ and $\bm{\Phi}_j\leq 0$. By Eq.~\eqref{eq:cts-no-de-no-bingding-def}, we have $\bm{\Phi}_j < 0$. By Lemma~\ref{lem:cts-non-de-bounded-opt}, we have a small positive constant $\epsilon>0$ such that
\begin{align}
    0<\epsilon<\min\left(\frac{\bar{r}}{\underline{b}}-\left\|\bp^*(\bm{B}')\right\|_1, -\frac{2\Phi_j}{\mu\bar{a}^2}  \right). \notag 
\end{align}
Consider $\bp''=\bp^*(\bm{B}')+\epsilon\cdot\bm{e}_j\in\Omega_p$. By Lemma~\ref{lem:cts-no-de-unique}, we have
\begin{align}
    f(\bm{B}',\bp'')-f(\bm{B}',\bp^*(\bm{B}'))&\leq \left\langle \bm{\Phi} , \bp''-\bp^*(\bm{B}')\right\rangle + \frac{\mu\bar{a}^2}{2}\left\| \bp''-\bp^*(\bm{B}') \right\|_2^2 =\Phi_j\cdot \epsilon + \frac{\mu\bar{a}^2}{2}\cdot\epsilon^2 < 0. \notag
\end{align}
Thus, we establish Eq.~\eqref{eq:cts-no-de-bingding-def} assuming the first and second parts, and Eq.~\eqref{eq:cts-no-de-no-bingding-def} in Definition~\ref{def:cts-non-degeneracy}.
\hfill\Halmos

\subsection{Proof of Lemma~\ref{lem:cts-non-de-decomposition-basic-item}}
\label{sec:pf-cts-non-de-decomposition-basic-item}
We first prove the first part of Lemma~\ref{lem:cts-non-de-decomposition-basic-item}. We have
\begin{align}
    |g(\bp')|&\leq \sum_{j=1}^m B_j p^*_j + \mathbb{E}\left[\bar{r}+\bar{a}\sum_{j=1}^m p_j^*\right]  \leq \bar{b}\cdot\left\|\bp^*\right\|_1+\bar{r}+\bar{a}\left\|\bp^*\right\|_1\leq \bar{r}+\bar{r}\cdot\frac{\bar{b}+\bar{a}}{\underline{b}}, \notag
\end{align}
where the last inequality is due to Lemma~\ref{lem:cts-non-de-bounded-opt}.

We next prove the second part of Lemma~\ref{lem:cts-non-de-decomposition-basic-item}. We have
\begin{align}
    g(\bp^*)-g(\bp')&=\mathbb{E}_{(r,\bm{a},\bm{b})\sim\mathcal{P}}\left[  \left(r-\left\langle\bm{a},\bp^*\right\rangle\right) \cdot \left(\mathbb{I}(r>\left\langle\bm{a},\bp^*\right\rangle) - \mathbb{I}(r>\left\langle\bm{a},\bp'\right\rangle) \right)\right] \notag\\
    &= \mathbb{E}_{(r,\bm{a},\bm{b})\sim\mathcal{P}}\left[\left[r-\bm{a}\cdot\bp^*\right]^+ - \left(r-\bm{a}\cdot\bp^*\right)\cdot\mathbb{I}(r>\bm{a}\cdot\bp') \right]\geq 0. \notag
\end{align}
We also have
\begin{align}
    g(\bp^*)-g(\bp')&=\mathbb{E}_{(r,\bm{a},\bm{b})\sim\mathcal{P}}\left[  \left(r-\left\langle\bm{a},\bp^*\right\rangle\right) \cdot \left(\mathbb{I}(r>\left\langle\bm{a},\bp^*\right\rangle) - \mathbb{I}(r>\left\langle\bm{a},\bp'\right\rangle) \right)\right],\notag \\
    &\leq \mathbb{E}_{(r,\bm{a},\bm{b})\sim\mathcal{P}}\left[  \left|r-\left\langle\bm{a},\bp^*\right\rangle\right| \cdot\left| \mathbb{I}(r>\left\langle\bm{a},\bp^*\right\rangle) - \mathbb{I}(r>\left\langle\bm{a},\bp'\right\rangle) \right| \right],\notag \\
    & \leq \mathbb{E}_{(r,\bm{a},\bm{b})\sim\mathcal{P}}\left[  \left|\left\langle\bm{a},\bp'\right\rangle-\left\langle\bm{a},\bp^*\right\rangle\right| \cdot\left| \mathbb{I}(r>\left\langle\bm{a},\bp^*\right\rangle) - \mathbb{I}(r>\left\langle\bm{a},\bp'\right\rangle) \right| \right],\notag \\
    & \leq \mu\left|\left\langle\bm{a},\bp'\right\rangle-\left\langle\bm{a},\bp^*\right\rangle\right|^2\leq \mu \bar{a}^2\left\|\bp'-\bp^*\right\|_2^2, \notag
\end{align}
where the third inequality is due to the second part of Definition~\ref{def:cts-non-degeneracy}. \hfill \Halmos

\subsection{Proof of Lemma~\ref{lem:cts-no-de-inventory-accumulation}}
\label{sec:pf-cts-no-de-first-t-2-performance}
We first prove the three items in Lemma~\ref{lem:cts-no-de-inventory-accumulation}. Then, the lemma is proved by collecting failure probabilities.

\medskip
\noindent\underline{Item~\ref{itm:inventory-accumulation-1} of Lemma~\ref{lem:cts-no-de-inventory-accumulation}}. Consider any $t$ such that $\kappa\leq t<H$, and let $w$ be such that $V_w-1\leq t<V_{w+1}-1$. We define the random variable $\bp^{t+1,*}$ by
\begin{align}
    \bp^{t+1,*} \in\arg\min_{\bp\geq\bm{0}}  f\left(\bp; \frac{\sum_{k=1}^{V_w-1} \bm{b}_k }{V_w-1}-C_2\sqrt{\frac{\ln H}{V_{w+1}-V_w}}\cdot\bm{1}\right). \notag
\end{align}
We also define the following events: 
\begin{align}
    &\mathcal{E}_{w}^{\text{resource},1} \defeq \left\{\forall j\in[m],~\frac{\sum_{k=1}^{V_w-1} b_{jk}}{V_{w}-1}-C_2\sqrt{\frac{\ln{H}}{V_{w+1}-V_w}}\in \left(B_j-\delta_B,B_j+\delta_B\right) \right\}, \label{eq:cts-no-de-resource-1-def} \\
    &\mathcal{E}_{t}^{\text{resource},2} \defeq \left\{\forall j\in[m],~\left|\frac{\sum_{k=1}^t b_{jk}}{t} - B_j\right|\leq \sqrt{\frac{5\bar{b}^2\cdot\ln H}{t}}   \right\}.\label{eq:cts-no-de-resource-2-def}
\end{align}
By Hoeffding's inequality, we have 
\begin{align}
    1-\left(\mathcal{E}_{w}^{\text{resource},1}\right)&\leq \sum_{j\in[m]} \Pr\left( \left| \frac{\sum_{k=1}^{V_w-1} b_{jk}}{V_w-1} - B_j \right|\geq \delta_B-C_2\sqrt{\frac{\ln H}{V_{w+1}-V_w}} \right)\notag \\
    &\leq \sum_{j\in[m]} \Pr\left( \left| \frac{\sum_{k=1}^{V_w-1} b_{jk}}{V_w-1} - B_j \right|\geq \frac{\delta_B}{2} \right)\leq \sum_{j\in[m]} 2\exp\left( -\frac{\delta_B^2 \cdot (V_w-1)}{2\bar{b}} \right)\leq \frac{2m}{H^5}, \label{eq:cts-no-de-saving-period-1}
\end{align}
where the second inequality is due to $V_{w+1}-V_w\geq  C_1\ln H$ and {\color{black} $C_1\geq\frac{4C_2^2}{\delta_B^2}$}, and the last inequality is due to $V_w> \kappa=C_0\ln ^2 H$ and {\color{black} $C_0\geq \frac{10\bar{b}^2}{\delta_B^2}$}. We also have 
\begin{align}
     1-\Pr\left(\mathcal{E}_{t}^{\text{resource},2}\right)&\leq \sum_{j\in[m]} \Pr\left( \left| \frac{\sum_{k=1}^t b_{jk}}{t} - B_j \right|\geq \cdot\sqrt{\frac{5\bar{b}^2\cdot\ln H}{t}} \right)\leq m\exp\left( -\frac{10\bar{b}^2\ln H}{\bar{b}^2} \right)\leq \frac{m}{H^5}. \label{eq:cts-no-de-saving-period-2}
\end{align}
We first estimate $\left\|\bp^{t+1}-\bp^*\right\|_2^2$ by considering $\left\|\bp^{t+1}-\bp^{t+1,*}\right\|_2$ and $\left\|\bp^*-\bp^{t+1,*}\right\|_2$. Due to Corollary~\ref{cor:cts-no-de-dual-conver-simple} and {\color{black}$C_0\geq \frac{20\bar{a}}{\lambda_{\text{min}}}+10$}, we have that 
\begin{align}
    \Pr\left( \left\|\bp^{t+1,*}-\bp^{t+1}\right\|_2^2\cdot\mathbb{I}(\mathcal{E}_{w}^{\text{resource},1})\leq (C_5)^2\cdot\frac{\ln H}{V_w-1} \right) \geq 1-\frac{3m+2}{H^5}. \notag
\end{align}
Further invoking Eq.~\eqref{eq:cts-no-de-saving-period-1} and the condition that {\color{black} $C_6\geq 16C_5$ and $C_0\ln H\geq C_1$}, we have 
\begin{align}
    \Pr\left( \left\|\bp^{t+1,*}-\bp^{t+1}\right\|_2^2\leq \left(\frac{C_6}{8}\right)^2\cdot\frac{\ln H}{t+1-\kappa} \right) \geq 1-\frac{5m+2}{H^5}. \label{eq:cts-no-de-saving-period-3}
\end{align}
On the other hand, by Lemma~\ref{lem:cts-no-de-d-robust}, we have that when the events $\mathcal{E}_{w}^{\text{resource},1}\cap \mathcal{E}_{t}^{\text{resource},2}$ happen, it holds that 
\begin{align}
    \left\| \bp^{t+1,*} - \bp^* \right\|_2^2 \leq\frac{2m}{\lambda^2\lambda^2_{\text{min}}} \cdot\frac{((C_2)^2+5\bar{b}^2) \cdot \ln H}{\min\left(t,V_{w+1}-V_w \right)}\leq \left(\frac{C_6}{8}\right)^2\cdot\frac{\ln H}{\min\left(t,V_{w+1}-V_w \right)}, \notag
\end{align}
where the last inequality is due to {\color{black} $C_6\geq \frac{8}{\lambda\lambda_{\text{min}}}\cdot\sqrt{2m\cdot ((C_2)^2 +5\bar{b}^2 )}$}. Since 
$\Pr(\mathcal{E}_{t}^{\text{resource},1}\cap \mathcal{E}_{t}^{\text{resource},2})\geq 1-\frac{3m}{H^5}$ (by Eq.~\eqref{eq:cts-no-de-saving-period-1} and Eq.~\eqref{eq:cts-no-de-saving-period-2}), we have 
\begin{align}
    \Pr\left(  \left\| \bp^{t+1,*} - \bp^* \right\|_2^2 \leq \left(\frac{C_6}{8}\right)^2\cdot\frac{\ln H}{\min\left(t,V_{w+1}-V_w \right)} \right) \geq 1- \frac{3m}{H^5}. \label{eq:cts-no-de-saving-period-4}
\end{align}
Combining Eq.~\eqref{eq:cts-no-de-saving-period-3} and Eq.~\eqref{eq:cts-no-de-saving-period-4}, we have
\begin{align}
     \Pr\left( \left\| \bp^{t+1} - \bp^* \right\|_2^2 \leq \left(\frac{C_6}{4}\right)^2\cdot\frac{\ln H}{\min\left(t+1-\kappa,V_{w+1}-V_w \right)} \right) \geq 1-\frac{8m+2}{H^5}.  \notag
\end{align}
Furthermore, for any $V_w-1\leq t<V_{w+1}-1$, by the definition of $V_w$, we have $t+1-\kappa\leq 4(V_{w+1}-V_w)$, and therefore derive that 
\begin{align}
     \Pr\left(\forall \kappa<t\leq H,~ \left\| \bp^{t} - \bp^* \right\|_2^2 \leq \left(C_6\right)^2\cdot\frac{\ln H}{t-\kappa} \right) \geq 1-\frac{8m+2}{H^4}.  
\end{align}

We next prove that $\bp^t\in\Omega_p$ with high probability. Consider any $w$ and $t=V_w-1$, we prove that when the event $\mathcal{E}_w^{\text{resource},1}$ (defined in Eq.~\eqref{eq:cts-no-de-resource-1-def}) happens, it holds that $\bp^{t+1}\in\Omega_p$. Observe that
\begin{align}
    \bar{r}>\frac{1}{t}\sum_{k=1}^t [r_t]^+ \geq \sum_{j\in[m]} p_j \left( \frac{\sum_{k=1}^t b_{jk}}{t} - C_2\sqrt{\frac{\ln H}{V_{w+1}-V_w}} \right)\geq \underline{b}\left\|\bp^{t+1}\right\|_1, \notag
\end{align}
where the second inequality is due to the fact that $\bp=\bp^{t+1}$ achieves a smaller objective value than $\bp=\bm{0}$ in Eq.~\eqref{eq:cts-no-de-choose-pt-saving-period-1}. Thus, 
\begin{align}
    \Pr\left(\forall \kappa<t \leq H ,~\bp^t\in\Omega_p\right)\geq \Pr\left( \bigcap_{w} \mathcal{E}_w^{\text{resource},1}  \right)\geq 1-\frac{2m}{H^4}.
\end{align}

\medskip
\noindent\underline{Item~\ref{itm:inventory-accumulation-2} of Lemma~\ref{lem:cts-no-de-inventory-accumulation}}. We have the following lemma, whose proof is deferred to Section~\ref{sec:pf-cts-no-de-equivalence-no-out-of-stock}. 
\begin{lemma}
\label{lem:cts-no-de-equivalence-no-out-of-stock}
    Suppose for all $\kappa < t \leq H$ and $j\in[m]$, we have 
    \begin{align}
        \sum_{k=1}^{t} b_{jk} - \sum_{k= \kappa+1}^t a_{jk}\cdot \mathbb{I}\left( r_k>\left\langle \bm{a}_k , \bp^k\right\rangle \right)\geq 0. \label{eq:cts-no-de-no-out-of-stock-basic}
    \end{align}
    Then, for all $v\in(\kappa,H]$, we have $x_v=\mathbb{I}(r_v>\left\langle \bm{a}_v , \bp^v \right\rangle)$.
\end{lemma}
By Lemma~\ref{lem:cts-no-de-equivalence-no-out-of-stock}, it suffices to show that Eq.~\eqref{eq:cts-no-de-no-out-of-stock-basic} holds. The following lemma, whose proof is deferred to Section~\ref{sec:pf-cts-no-de-saving-period-each-sub-inventory}, shows the inventory lower bounds at the end of each time period $t = V_w$.
\begin{lemma}
\label{lem:cts-no-de-saving-period-each-sub-inventory}
    With probability $1-\frac{27m^2+14m}{H^4}$, the following inequalities hold simultaneously: 
    \begin{align}
        &\sum_{k=V_w}^{V_{w+1}-1} b_{jk} - a_{jk}\cdot\mathbb{I}(r_k>\left\langle \bm{a}_k ,  \bp^k\right\rangle) \geq 4C_7\sqrt{\left( V_{w+1}-V_w \right)\cdot\ln H},&&\forall w\in[N_V],j\in[m],\notag \\
        &\sum_{k=V_w}^{V_{w+1}-1} b_{jk} - a_{jk}\cdot\mathbb{I}(r_k>\left\langle \bm{a}_k ,  \bp^k\right\rangle) \leq \frac{C_8}{6}\sqrt{\left( V_{w+1}-V_w \right)\cdot\ln H},&&\forall w\in[N_V], j\in I_B,\notag \\
        &\sum_{k=V_w}^{t} b_{jk} - a_{jk}\cdot\mathbb{I}(r_k>\left\langle \bm{a}_k ,  \bp^k\right\rangle) \geq -C_7\sqrt{\left( V_{w+1}-V_w \right)\cdot\ln H},&&\forall w\in[N_V], j\in[m], t \in [V_w, t<V_{w+1}), \notag \\
        & \sum_{k=1}^{V_1-1} b_{jk} \geq C_7\cdot \sqrt{(V_2-V_1)\cdot\ln H}, &&\forall j\in[m]. \notag 
    \end{align}
\end{lemma}
We next show that if the inequalities in Lemma~\ref{lem:cts-no-de-saving-period-each-sub-inventory} all hold, then Eq.~\eqref{eq:cts-no-de-no-out-of-stock-basic} holds for all $\kappa< t \leq H$. We consider any $t$ with $V_w\leq t< V_{w+1}$ and discuss the following two cases: $w=1$ and $w>1$.
\begin{itemize}
    \item \underline{Case 1: $w=1$.} By the third and fourth inequalities in Lemma~\ref{lem:cts-no-de-saving-period-each-sub-inventory}, we have
    \begin{align}
        \sum_{k=1}^{t} b_{jk} - \sum_{k= \kappa+1}^t a_{jk}\cdot \mathbb{I}\left( r_k>\left\langle \bm{a}_k , \bp^k \right\rangle \right) &= \sum_{k=1}^{V_1-1} b_{jk} +\sum_{k=V_1}^t  b_{jk} - a_{jk}\cdot\mathbb{I}\left(r_k>\left\langle \bm{a}_k , \bp^k\right\rangle \right)\notag \\
        &\geq C_7\cdot\sqrt{(V_2-V_1)\cdot\ln H}- C_7\cdot\sqrt{(V_2-V_1)\cdot\ln H}=0. \notag
    \end{align}
    \item \underline{Case 2: $w>1$.} By the first and third inequalities in Lemma~\ref{lem:cts-no-de-saving-period-each-sub-inventory}, we have
    \begin{align}
        &\sum_{k=1}^{t} b_{jk} - \sum_{k= \kappa+1}^t a_{jk}\cdot \mathbb{I}\left( r_k>\left\langle \bm{a}_k , \bp^k\right\rangle \right)\notag \\
        &\qquad \geq \sum_{v=1}^{w-1} \sum_{k=V_v}^{V_{v+1}-1} b_{jk} - a_{jk}\cdot\mathbb{I}\left(r_k>\left\langle \bm{a}_k , \bp^k\right\rangle\right) +\sum_{k=V_{w}}^{t} b_{jk} - a_{jk}\cdot\mathbb{I}\left(r_k>\left\langle \bm{a}_k , \bp^k\right\rangle\right)\notag \\
        &\qquad\geq \sum_{k=V_{w-1}}^{V_{w}-1} b_{jk} - a_{jk}\cdot\mathbb{I}\left(r_k>\left\langle \bm{a}_k , \bp^k\right\rangle\right) +\sum_{k=V_{w}}^{t} b_{jk} - a_{jk}\cdot\mathbb{I}\left(r_k>\left\langle \bm{a}_k , \bp^k\right\rangle\right) \notag \\
        &\qquad \geq C_7\sqrt{\ln H}\cdot\left(4\sqrt{V_w-V_{w-1}}-\sqrt{V_{w+1}-V_w}\right)\geq 0, \notag
    \end{align}
    where the last inequality is due to that $V_{w+1}-V_{w}\leq 16(V_{w}-V_{w-1})$.
\end{itemize}

\medskip
\noindent\underline{Item~\ref{itm:inventory-accumulation-3} of Lemma~\ref{lem:cts-no-de-inventory-accumulation}}. We prove that if the inequalities in Lemma~\ref{lem:cts-no-de-saving-period-each-sub-inventory} all hold, then Item~\ref{itm:inventory-accumulation-3} of Lemma~\ref{lem:cts-no-de-inventory-accumulation} holds. For any $j\in[m]$, we have that 
\begin{align}
    \ell_{j,H+1}&=\sum_{k=1}^{V_1-1}b_{jk} + \sum_{w=1}^{N_V}\sum_{t=V_w}^{V_{w+1}-1}  b_{jk} - a_{jk}\cdot\mathbb{I}\left(r_k>\left\langle \bm{a}_k , \bp^k\right\rangle\right) \notag \\
    &\geq 4C_7\sqrt{(V_{N_V+1}-V_{N_V})\cdot\ln H}\geq 2C_7\sqrt{\left(H-C_0\ln^2 H\right)\cdot\ln H}\geq C_7\sqrt{H\ln H}, \notag
\end{align}
where the last inequality is due to that {\color{black}$H>2+2C_0\ln^2 H$}. For $j\in I_B$, we have
\begin{align}
    \ell_{j,H+1}&=\sum_{k=1}^{V_1-1}b_{jk} + \sum_{w=1}^{N_V}\sum_{t=V_w}^{V_{w+1}-1}  b_{jk} - a_{jk}\cdot\mathbb{I}\left(r_k>\left\langle \bm{a}_k , \bp^k\right\rangle \right) \notag \\
    &\leq \bar{b}\cdot(V_1-1)+\sum_{w=1}^{N_V} \frac{C_8}{6}\sqrt{(V_{w+1}-V_w)\cdot\ln H}\leq \bar{b}\cdot C_0\ln^2 H+\sum_{w=1}^{N_V} \frac{C_8}{6}\sqrt{3\cdot 2^{w-1}\cdot C_1 \ln^2 H}\notag \\
    &\leq \bar{b}\cdot C_0\ln^2 H + \frac{\sqrt{3C_1}C_8\ln H}{6}\cdot\frac{2^{\frac{N_V}{2}}-1}{\sqrt{2}-1} <\bar{b}\cdot C_0\ln^2 H+ 0.8C_8\cdot\sqrt{H\ln H} \leq C_8\sqrt{H\ln H}, \notag
\end{align}
where the first inequality is due to Lemma~\ref{lem:cts-no-de-saving-period-each-sub-inventory}, the second inequality is by the definition of $V_w$, and the last one is due to that {\color{black} $C_8\sqrt{H\ln H}\geq 5\bar{b}\cdot C_0\ln^2 H$}.
\hfill\Halmos

\subsection{Proof of Lemma~\ref{lem:cts-no-de-detection-estimation}}
\label{sec:pf-cts-no-de-second-t-2-performance-pre}
We first prove the three items in Lemma~\ref{lem:cts-no-de-detection-estimation}. Then, the lemma is proved by collecting failure probabilities.

\medskip
\noindent{\underline{Item~\ref{itm:inventory-accumulation-1} of Lemma~\ref{lem:cts-no-de-detection-estimation}.}} We first prove that with probability at least $1-\frac{1}{H^5}$, it holds that $\widehat{I}_B=I_B$. In this proof, let $M=\lfloor\frac{H}{2}\rfloor$. Consider the following event, denoted by $\mathcal{E}^{\text{test,1}}$:
(1) for all $j\in[m]$, we have it holds that $
        \left|\frac{\sum_{k=1}^{M} b_{jk}}{M} -B_j\right| < \min\left( \delta_B , \sqrt{\frac{5\bar{b}^2\cdot\ln H}{M}} \right)$;
(2) for any $\bm{B}'\in\Omega_b$, denote $\bp'$ as any minimizer in $
        \arg\min_{\bp\geq\bm{0}} f\left(\bp;\bm{B}',\{r_t,\bm{a}_t,\bm{b}_t\}_{t\in[M]}\right)$, and we have
    $
        \left\| \bp'-\bp^*(\bm{B}') \right\|_2^2 \leq (C_5)^2\cdot\frac{\ln H}{M}$.
We have the following lemma, whose proof is deferred to Section~\ref{sec:pf-cts-no-de-test-index-1}.
\begin{lemma}
\label{lem:cts-no-de-test-index-1}
    $\Pr(\mathcal{E}^{\mathrm{test,1}})\geq 1-\frac{3m+4}{H^5}$.
\end{lemma}
When $\mathcal{E}^{\text{test,1}}$ happens, we have 
\begin{align}
    \left\|\widehat{\bp}^*-\bp^*\right\|_2^2 &\leq 2 \left\|\widehat{\bp}^*-\bp^*\left(\frac{\sum_{k=1}^{M} \bm{b}_{k}}{M}\right) \right\|_2^2 + 2  \left\|\bp^*\left(\frac{\sum_{k=1}^{M} \bm{b}_{k}}{M}\right)  - \bp^*\right\|_2^2 \notag \\
    &\leq (C_5)^2\cdot\frac{\ln H}{M} + \frac{2}{\lambda^2\lambda^2_{\text{min}}}\cdot\left\| \frac{\sum_{k=1}^{M} \bm{b}_{k}}{M} - \bm{B} \right\|_2^2\leq (C_5)^2\cdot\frac{\ln H}{M} + \frac{10m\bar{b}^2\cdot\ln H}{M\lambda^2\lambda^2_{\text{min}}}\leq (C_5)^2\cdot\frac{3\ln H}{M}, \notag
\end{align}
where the second inequality is due to Lemma~\ref{lem:cts-no-de-d-robust} and the definition of $\mathcal{E}^{\text{test,1}}$, and the last one is due to that {\color{black}$C_5\geq \frac{\sqrt{5m}\cdot\bar{b}}{\lambda\lambda_{\text{min}}}$}. One may also check that when $\mathcal{E}^{\text{test,1}}$ happens, it holds that $\widehat{\bp}^*\in\Omega_p$. Let $\mathcal{F}$ be $\{(r_k,\bm{a}_k,\bm{b}_k)\}_{k\leq M}$, and we consider the case where the event $\mathcal{E}^{\text{test,1}}$ occurs in $\mathcal{F}$ (denoted by $\mathcal{E}^{\text{test,1}}\in\mathcal{F}$ for convenience). Observe that, conditioned on $\mathcal{F}$, $\{(r_k,\bm{a}_k,\bm{b}_k)\}_{M<k\leq H}$ are independent and identically distributed. For any $\mathcal{F}$ such that $\mathcal{E}^{\text{test,1}}\in\mathcal{F}$, we have
\begin{align}
    &\left|\mathbb{E}\left[ b_{jk} - a_{jk}\cdot\mathbb{I}\left(r_k>\left\langle \bm{a}_k , \widehat{\bp}^*\right\rangle \right) \middle|\mathcal{F} \right]  -  B_j +  \mathbb{E}\left[ a_{jk}\cdot\mathbb{I}\left(r_k>\left\langle \bm{a}_k , \bp^*\right\rangle \right) \middle|\mathcal{F}  \right] \right| \notag \\
    &\qquad \leq \bar{a}\left| \Pr\left( r_k > \left\langle \bm{a}_k , \widehat{\bp}^*\right\rangle\middle|\mathcal{F}\right) -\Pr\left( r_k > \left\langle\bm{a}_k , \bp^*\right\rangle\middle|\mathcal{F}\right) \right|\notag \\
    &\qquad \leq \bar{a}\mu\cdot \mathbb{E}\left[ \left|  \left\langle \bm{a}_k , \widehat{\bp}^*\right\rangle- \left\langle\bm{a}_k , \bp^*\right\rangle \right| \middle|\mathcal{F}\right]\leq \bar{a}^2\mu\cdot C_5\cdot\sqrt{\frac{3\ln H}{M}}, \label{eq:cts-no-de-test-opt-est}
\end{align}
where the second inequality is by the definition of non-degenerate distributions (Definition~\ref{def:cts-non-degeneracy}). Moreover, by Lemma~\ref{lem:cts-no-de-many-inv-for-no-binding}, we have
\begin{align}
    & B_j -  \mathbb{E}\left[ a_{jk}\cdot\mathbb{I}( r_k>\left\langle \bm{a}_k , \bp^*\right\rangle) \middle|\mathcal{F}\right]= 0,\qquad&&\forall j\in I_B, \notag\\
    & B_j -  \mathbb{E}\left[ a_{jk}\cdot\mathbb{I}( r_k>\left\langle\bm{a}_k , \bp^*\right\rangle) \middle|\mathcal{F}\right]\geq \delta_B,\qquad&&\forall j\notin I_B.  \notag
\end{align}
Finally, we can estimate $\widehat{\ell}_j$ for $j\in I_B$ and $j\notin I_B$:
\begin{itemize}
    \item For each $j\in I_B$: by the third part of Definition~\ref{def:cts-non-degeneracy} and Eq.~\eqref{eq:cts-no-de-test-opt-est}, we have 
    \begin{align}
        \mathbb{E}\left[ b_{jk}-a_{jk}\cdot\mathbb{I}(r_k>\bm{a}_k\cdot\widehat{\bp}^*) \middle|\mathcal{F} \right]\leq C_5\bar{a}^2\mu\cdot\sqrt{\frac{3\ln H}{M}}. \notag
    \end{align}
    By Hoeffding's inequality and {\color{black}$C_3\geq 2\sqrt{3}C_5\bar{a}^2\mu $}, we have when $\mathcal{E}^{\text{test,1}}\in\mathcal{F}$ occurs, it holds that
    \begin{align}
         \Pr\left( \sum_{k=M+1}^{H}  b_{jk}-a_{jk}\cdot\mathbb{I}\left( r_k>\left\langle \bm{a}_k , \widehat{\bp}^*\right\rangle \right) \geq C_3\sqrt{ H\ln H} \middle|\mathcal{F}\right)  \leq \exp\left( -\frac{(C_3)^2\cdot H\ln H}{2(H-M)(\bar{b}+2\bar{a})^2} \right)\leq \frac{1}{H^5}, \notag
    \end{align}
    where the last inequality is due to {\color{black}$C_3\geq3(\bar{b}+2\bar{a})$}. Thus, we have
    \begin{align}
        \Pr\left( \sum_{k=M+1}^{H}  b_{jk}-a_{jk}\cdot\mathbb{I}\left( r_k>\left\langle \bm{a}_k , \widehat{\bp}^*\right\rangle \right) \geq C_3\sqrt{ H\ln H} \right)\leq \frac{3m+6}{H^5}. \notag
    \end{align} 
    \item For each $j\notin I_B$, by Lemma~\ref{lem:cts-no-de-many-inv-for-no-binding} and Eq.~\eqref{eq:cts-no-de-test-opt-est}, we have 
    \begin{align}
        \mathbb{E}\left[ b_{jk}-a_{jk}\cdot\mathbb{I}\left( r_k>\left\langle \bm{a}_k , \widehat{\bp}^*\right\rangle \right) \middle|\mathcal{F} \right]\geq \delta_B - C_5\bar{a}^2\mu\cdot\sqrt{\frac{3\ln H}{M}}\geq\frac{\delta_B}{2}, \notag
    \end{align}
    where the last inequality is due to {\color{black}$\lfloor\frac{H}{2}\rfloor\cdot\delta_B^2\geq 12 (C_5)^2\bar{a}^4\mu^2 \ln H$}. By  Hoeffding's inequality  and {\color{black}$ H\cdot\delta_B^2\geq 400 (C_3)^2 \ln H $}, we have that when $\mathcal{E}^{\text{test,1}}\in\mathcal{F}$, it holds that
    \begin{align}
         \Pr\left( \sum_{k=M+1}^{H}  b_{jk}-a_{jk}\cdot\mathbb{I}\left( r_k>\left\langle \bm{a}_k , \widehat{\bp}^*\right\rangle \right) \leq C_3\sqrt{ H\ln H} \middle|\mathcal{F}\right)  \leq \exp\left( -\frac{2H\delta_B^2}{25 (\bar{b}+2\bar{a})^2} \right)\leq \frac{2}{H^5}, \notag
    \end{align}
    where the last inequality is due to {\color{black}$H\cdot\delta_B^2\geq 63(\bar{b}+2\bar{a})^2\ln H$}. Thus, we have
    \begin{align}
        \Pr\left( \sum_{k=M+1}^{H}  b_{jk}-a_{jk}\cdot\mathbb{I}\left( r_k>\left\langle \bm{a}_k , \widehat{\bp}^*\right\rangle \right) \leq C_3\sqrt{ H\ln H} \right)\leq \frac{3m+6}{H^5}. \notag
    \end{align}
\end{itemize}
In all, we have that, with probability at least $1-\frac{6m+12}{H^5}$, it holds that $\widehat{I}_B=I_B$.

\medskip
\noindent{\underline{Item~\ref{itm:inventory-accumulation-2} of Lemma~\ref{lem:cts-no-de-detection-estimation}.}} By Hoeffding's inequality, we have
\begin{align}
    \Pr\left( \left| \widehat{B}_j - B_j\right|\geq \sqrt{\frac{5\bar{b}^2\ln H}{H}} \right) &= \Pr\left( \left| \frac{1}{H}\cdot\sum_{k=1}^{H} ( b_{jk} - B_j)  \right|\geq \sqrt{\frac{5\bar{b}^2\ln H}{H}}  \right) \leq 2\exp\left( -\frac{5 H \ln H}{H} \right)\leq \frac{2}{H^5}, \notag
\end{align}
where the last inequality is due to {\color{black}$H\geq 10$}.

\medskip
\noindent{\underline{Item~\ref{itm:inventory-accumulation-3} of Lemma~\ref{lem:cts-no-de-detection-estimation}.}} This part of the lemma can be directly derived from Corollary~\ref{cor:cts-no-de-dual-conver-simple} and the condition {\color{black}$H\geq \frac{20\bar{a}}{\lambda_{\text{min}}}\cdot\ln H +10$}.
\hfill\Halmos

\subsection{Proof of Lemma~\ref{lem:cts-no-de-inventory-conversion}}
\label{sec:pf-cts-no-de-second-t-2-performance}
Note that $\bp^t$ remains the same in the same batch of the algorithm. For simplicity, in this proof, we let $\bq_{s} = \bp^t$ for the all $t$ in batch $s$ (i.e., $t \in [U_s, U_{s+1})$). We also define  
$\bq^{*}_{s}\defeq \bp^*\left( {\widetilde{\bm{B}}_s}\right)$.
We will examine the inequalities in Lemma~\ref{lem:cts-no-de-inventory-conversion} sequentially for increasing values of the batch index $s$. 

Let $\mathcal{F}_s=\{(r_k,\bm{a}_k,\bm{b}_k)\}_{k<U_s}$. We also define the event $\mathcal{E}_s$ ($s\geq 1$) as follows:
\begin{align}
\mathcal{E}_s \defeq \left\{\left| \ell_{j,U_s}- C_4\sqrt{3^{N_U+2-s}\cdot\ln H} \right|\leq C_9 \sqrt{3^{N_U+1-s}\cdot\ln H},  \forall j\in I_B, \text{~and~} \ell_{j,U_s} \geq C_7\sqrt{H\ln H}, \forall j\notin I_B\right\}. \notag
\end{align}
We further let $\mathcal{E}_s\in\mathcal{F}_s$ denote the event $\mathcal{E}_s$ occurring in the filtration $\mathcal{F}_s$. We have $U_2-U_1\in \left[\frac{2H}{3}+1,\frac{8H}{9}+1\right]$ and $U_{s+1}-U_s=2\cdot 3^{N_U-s}$ for $s\geq 2$. When $\mathcal{E}_s\in\mathcal{F}_s$, we discuss the three items in Lemma~\ref{lem:cts-no-de-inventory-conversion} for $s+1$.

\medskip
\noindent \underline{Item~\ref{itm:cts-no-de-inventory-conversion-1} of Lemma~\ref{lem:cts-no-de-inventory-conversion}.} When $\mathcal{E}_s\in\mathcal{F}_s$, we have 
\begin{align}
    \left|{\widetilde{B}_{j,s}} - B_j \right| &\leq \left| \widehat{B}_j - B_j \right| + \left| \frac{\ell_{j,U_s} - C_4\sqrt{3^{N_U +1 - s}\cdot\ln H} }{U_{s+1}-U_s} \right|\cdot\mathbb{I}(j\in I_B),\notag \\
    &\qquad\qquad\qquad\qquad\qquad\leq \sqrt{\frac{5\bar{b}^2\ln H}{H}} + \frac{3(C_8 +C_9)  \sqrt{3^{N_U+1-s}\ln H} }{3^{N_U+1-s}} < \delta_B, \notag
\end{align}
where the second inequality is due to that {\color{black}$C_8\geq C_4$}, and the last inequality is due to {\color{black}$ H\cdot \delta_B^2 \geq 20(\bar{b})^2 \cdot \ln H$} and {\color{black}$3^{N_U+1-s}\cdot\delta_B^2\geq 36(C_8+C_9)^2\cdot\ln H $}. Thus ${\widetilde{\bm{B}}_s}\in\Omega_b$, and we can directly verify that $\bq^*_s,\bq_s\in\Omega_p$. When $\mathcal{E}_s\in\mathcal{F}_s$, due to the assumptions on $\widehat{f}(\cdot;\cdot)$, we have
\begin{align}
    \left\|\bq^*_s - \bq_s \right\|_2^2 \leq (C_5)^2\cdot\frac{\ln H}{H}. \notag
\end{align}
Moreover, by $\widetilde{\bm{B}}_s\in\Omega_b$ and Lemma~\ref{lem:cts-no-de-d-robust}, we have
\begin{align}
     \left\|\bq^*_s - \bp^{*}(\bm{B})\right\|_2^2 \leq \frac{1}{\lambda^2\lambda^2_{\min}} \cdot  \left\| \widetilde{\bm{B}}_s - \bm{B} \right\|_2^2 \leq \frac{m}{\lambda^2\lambda^2_{\min}}\cdot\left( \frac{10\bar{b}^2\ln H}{H} + \frac{18(C_8+C_9)^2\ln H}{3^{N_U+1-s}} \right). \notag
\end{align}
Combining the two inequalities above, we have that when $\mathcal{E}_s\in\mathcal{F}_s$, it holds that
\begin{align}
    \left\|\bq_s - \bp^{*}(\bm{B}) \right\|_2^2 &\leq  2\cdot\left\|\bq^*_s - \bq_s \right\|_2^2    + 2\cdot\left\|\bq^*_s - \bp^{*} \right\|_2^2 \notag \\
    &\leq  (C_5)^2\cdot\frac{6\ln H}{H} +\frac{m}{\lambda^2\lambda^2_{\min}}\cdot\left( \frac{20\bar{b}^2\ln H}{H} + \frac{18(C_8+C_9)^2\ln H}{3^{N_U+1-s}} \right) \leq (C_{10})^2\cdot \frac{\ln H}{3^{N_U+1-s}}, \notag
\end{align}
where the last inequality is due to that {\color{black}$(C_{10})^2\geq 6(C_5)^2 + \frac{m}{\lambda^2\lambda^2_{\min}}\cdot (20\bar{b}^2 + 18(C_8+C_9)^2) $}. 

\medskip
\noindent{\underline{Item~\ref{itm:cts-no-de-inventory-conversion-2} Lemma~\ref{lem:cts-no-de-inventory-conversion}.}} To prove this part of Lemma~\ref{lem:cts-no-de-inventory-conversion}, similar to Lemma~\ref{lem:cts-no-de-equivalence-no-out-of-stock}, it suffices to show 
\begin{align}
    \ell_{j,U_s} + \sum_{k=U_s}^{t} b_{jk} - a_{jk}\cdot\mathbb{I}\left( r_k > \left\langle\bm{a}_k , \bq_s\right\rangle\right)\geq 0,\qquad\forall t<U_{s+1}. 
\end{align}
Notice that we have that conditioned on $\mathcal{F}_s$, $\{(r_k,\bm{a}_k,\bm{b}_k)\}_{U_s\leq k<U_{s+1}}$ is independent and identically distributed, satisfying
\begin{align}
     &\left| \mathbb{E}\left[b_{jk} - a_{jk}\cdot\mathbb{I}\left( r_k > \left\langle \bm{a}_k , \bq_s\right\rangle \right)\middle|\mathcal{F}_s \right] -  \mathbb{E}\left[b_{jk} - a_{jk}\cdot\mathbb{I}\left( r_k > \left\langle \bm{a}_k , \bq^*_s\right\rangle\right)\middle|\mathcal{F}_s \right] \right|\notag \\
     &\qquad\qquad\qquad\qquad\qquad\qquad\qquad\qquad\qquad\qquad\leq \bar{a}^2\mu\left\| \bq_s - \bq^*_s \right\|_2 \leq C_5 \bar{a}^2\mu\cdot\sqrt{\frac{\ln H}{H}}.\label{eq:cts-no-de-second-performance-opt-est-on-expectation} 
\end{align}
For each $j$, we discuss the following two cases depending on whether $j\in I_B$:
\begin{itemize}
    \item For $j\in I_B$, when $\mathcal{E}_s\in\mathcal{F}_s$, we have $
     \mathbb{E}\left[b_{jk} - a_{jk}\cdot\mathbb{I}\left( r_k > \left\langle \bm{a}_k , \bq^*_s\right\rangle\right)\middle|\mathcal{F}_s \right]=B_j-{\widetilde{B}_{j,s}}$. 
Thus, we have
\begin{align}
    & \ell_{j,U_s}+\sum_{k=U_s}^{t}\mathbb{E}\left[b_{jk} - a_{jk}\cdot\mathbb{I}\left( r_k > \left\langle \bm{a}_k , \bq^*_s\right\rangle\right)\middle|\mathcal{F}_s \right] \geq \min\left(  \ell_{j,U_s}, \ell_{j,U_s}+B_j(U_{s+1}-U_s) - \widetilde{\ell}_{j,s} \right)\notag \\
    &\geq \min\left( \ell_{j,U_s}, C_4\sqrt{3^{N_U+1-s}\ln H} - |B_j-\widehat{B}_j|\cdot(U_{s+1}-U_s) \right) \notag \\
    &\geq \min\left( C_4\sqrt{3^{N_U-s+2}\ln H} - C_9\sqrt{3^{N_U-s+1}\ln H} , C_4\sqrt{3^{N_U+1-s}\ln H} - 8\bar{b}\sqrt{5\cdot 3^{N_U-s}\ln H}  \right) \notag \\
    &\geq\frac{C_4}{2}\sqrt{3^{N_U-s}\ln H}, \notag
\end{align}
where the last inequality is due to {\color{black}$C_4\geq 2C_9+16\sqrt{2}\bar{b}$}. By Eq.~\eqref{eq:cts-no-de-second-performance-opt-est-on-expectation}, we have
\begin{align}
    & \ell_{j,U_s}+\sum_{k=U_s}^{t}\mathbb{E}\left[b_{jk} - a_{jk}\cdot\mathbb{I}\left( r_k > \left\langle \bm{a}_k , \bq_s\right\rangle\right)\middle|\mathcal{F}_s \right] \notag \\
    &\qquad \geq \frac{C_4}{2}\sqrt{3^{N_U+1-s}\ln H} - C_5 \bar{a}^2\mu\cdot\sqrt{\frac{\ln H}{H}}\cdot (U_{s+1}-U_s)\geq \frac{C_4}{3}\sqrt{3^{N_U+1-s}\ln H} , \notag
\end{align}
where the last inequality is due to that {\color{black}$C_4\geq 16C_5\bar{a}^2\mu$}. By Hoeffding's inequality, we have 
\begin{align}
\Pr\left( \ell_{j,U_s} + \sum_{k=U_s}^t b_{jk} - a_{jk}\cdot\mathbb{I}\left( r_k > \left\langle \bm{a}_k , \bq_s\right\rangle\right)\leq 0 \middle| \mathcal{F}_s \right) \leq 2\exp\left( -\frac{2(C_4)^2 \cdot3^{N_U+1-s}\cdot \ln H}{9(\bar{b}+2\bar{a})^2 \cdot8\cdot3^{N_U-s}}  \right)\leq\frac{2}{H^6},   \notag
\end{align}
where the last inequality is due to that {\color{black}$C_4\geq 9(\bar{b}+2\bar{a})$}. Thus we can derive that when $\mathcal{E}_s\in\mathcal{F}_s$, for $j\in I_B$, we have
\begin{align}
    \Pr\left( \forall U_s\leq t<U_{s+1},~ \ell_{j,U_s} + \sum_{k=U_s}^t b_{jk} - a_{jk}\cdot\mathbb{I}\left( r_k > \left\langle \bm{a}_k , \bq_s\right\rangle\right)> 0 \middle| \mathcal{F}_s \right)\geq 1-\frac{2}{H^5}. \notag
\end{align}

\item For $j\notin I_B$, when $\mathcal{E}_s\in\mathcal{F}_s$, we have
\begin{align}
     \mathbb{E}\left[b_{jk} - a_{jk}\cdot\mathbb{I}\left( r_k > \left\langle \bm{a}_k , \bq_s\right\rangle\right)\middle|\mathcal{F}_s \right]\geq \delta_B - \bar{a}^2\mu\left\| \bq_s - \bp^{*}\right\|_2\geq \frac{\delta_B}{2}, \label{eq:cts-no-de-no-binding-second-expectation-big-0}
\end{align}
where the last inequality is due to that {\color{black}$3^{N_U+1-s}\cdot\delta_B^2\geq 4(C_{10})^2\bar{a}^4\mu^2\ln H $}. By Hoeffding's inequality, we have
\begin{align}
    \Pr\left( \ell_{j,U_s} + \sum_{k=U_s}^t b_{jk} - a_{jk}\cdot\mathbb{I}\left( r_k > \left\langle \bm{a}_k , \bq_s\right\rangle\right)\leq 0 \middle| \mathcal{F}_s \right) \leq 2\exp\left(-\frac{2(\ell_{j,U_s})^2}{(U_{s+1}-U_s)\cdot(\bar{b}+2\bar{a})^2} \right) \leq  \frac{2}{H^6}, \notag 
\end{align}
where the last inequality is due to {\color{black}$C_7\geq 2(\bar{b}+2\bar{a})$} and $\ell_{j,U_s}\geq C_7\sqrt{H\ln H}$. We thus derive that when $\mathcal{E}_s\in\mathcal{F}_s$, for $j\notin I_B$, it holds that
\begin{align}
    \Pr\left( \forall U_s\leq t<U_{s+1},~ \ell_{j,U_s} + \sum_{k=U_s}^t b_{jk} - a_{jk}\cdot\mathbb{I}\left( r_k > \left\langle \bm{a}_k , \bq_s\right\rangle\right)> 0 \middle| \mathcal{F}_s \right)\geq 1-\frac{2}{H^5}.
\end{align} 
\end{itemize}
Combining the above discussions, we conclude that when $\mathcal{E}_s\in\mathcal{F}_s$, it holds that
\begin{align}
    \Pr\left( \forall U_s\leq t<U_{s+1},~ x_t=\mathbb{I}\left( r_t>\left\langle \bm{a}_t , \bp^t \right\rangle\right) \middle|\mathcal{F}_s\right) \geq 1-\frac{2m}{H^5}.\label{eq:cts-no-de-second-no-out-of-stock}
\end{align}

\medskip
\noindent \underline{Item~\ref{itm:cts-no-de-inventory-conversion-3} Lemma~\ref{lem:cts-no-de-inventory-conversion}.} For each $j$, we discuss the following two cases depending on whether $j\in I_B$:
\begin{itemize}
    \item For $j\in I_B$, when $\mathcal{E}_s\in\mathcal{F}_s$, we have $
     \mathbb{E}\left[b_{jk} - a_{jk}\cdot\mathbb{I}\left( r_k > \left\langle \bm{a}_k , \bq^*_s\right\rangle\right)\middle|\mathcal{F}_s \right]=B_j-{\widetilde{B}_{j,s}}$. 
By Eq.~\eqref{eq:cts-no-de-second-performance-opt-est-on-expectation}, we have
\begin{align}
    &~~~\left|\ell_{j,U_s}- C_4\sqrt{3^{N_U+1-s}\cdot \ln H} + \sum_{k=U_s}^{U_{s+1}-1}\mathbb{E}\left[b_{jk} - a_{jk}\cdot\mathbb{I}\left( r_k > \left\langle \bm{a}_k , \bq_s\right\rangle\right)\middle|\mathcal{F}_s \right] \right| \notag \\
    &\leq\left|\ell_{j,U_s}- C_4\sqrt{3^{N_U+1-s}\cdot \ln H} + \sum_{k=U_s}^{U_{s+1}-1}\mathbb{E}\left[b_{jk} - a_{jk}\cdot\mathbb{I}\left( r_k > \left\langle \bm{a}_k , \bq^*_s\right\rangle\right)\middle|\mathcal{F}_s \right] \right|+C_5\bar{a}^2\mu\cdot(U_{s+1}-U_s)\cdot\sqrt{\frac{\ln H}{H}} \notag \\
    &=  \left|\ell_{j,U_s} - C_4\sqrt{3^{N_U+1-s}\cdot \ln H} + (B_j - \widetilde{B}_{j,s})\cdot(U_{s+1}-U_s) \right|+C_5\bar{a}^2\mu\cdot(U_{s+1}-U_s)\cdot\sqrt{\frac{\ln H}{H}}\notag \\
    &= \left(|\widehat{B}_j- B_j| + C_5\bar{a}^2\mu\sqrt{\frac{\ln H}{H}} \right)\cdot(U_{s+1}-U_s) \leq (8\bar{b}+4C_5\bar{a}^2\mu)\sqrt{2\cdot 3^{N_U-s}\ln H}. \notag
\end{align}
By Hoeffding's inequality, we have 
\begin{align}
    &~~~\Pr\left( \left|\ell_{j,U_s}- C_4\sqrt{3^{N_U+1-s}\cdot \ln H} + \sum_{k=U_s}^{U_{s+1}-1}\left[b_{jk} -  a_{jk}\cdot\mathbb{I}\left( r_k > \left\langle \bm{a}_k , \bq_s\right\rangle\right) \right] \right|\geq C_9\sqrt{3^{N_U-s}\cdot\ln H} \middle|\mathcal{F}_s\right),\notag \\
    & \leq 2\exp\left( - \frac{(C_9)^2 \cdot3^{N_U-s}\cdot\ln H}{2(U_{s+1}-U_s)\cdot(\bar{b}+2\bar{a})^2}\right)\leq\frac{2}{H^5}, \label{eq:cts-no-de-second-cntrol-stock-binding-pre}
\end{align}
where the first inequality is due to {\color{black}$ C_9 \geq16\sqrt{2}\cdot\bar{b}+8\sqrt{2}C_5\bar{a}^2\mu$}, and the second one is due to {\color{black} $C_9\geq 9\cdot(\bar{b}+2\bar{a})$}. By Eq.~\eqref{eq:cts-no-de-second-no-out-of-stock} and Eq.~\eqref{eq:cts-no-de-second-cntrol-stock-binding-pre}, we have that when $\mathcal{E}_s\in\mathcal{F}_s$, it holds that
\begin{align}
    \Pr\left(\forall j\in I_B, \left|\ell_{j,U_{s+1}} - C_4\sqrt{3^{N_U+1-s}\cdot \ln H} \right| \leq C_9\sqrt{3^{N_U-s}\cdot\ln H}  \middle|\mathcal{F}_s\right)\geq 1-\frac{4m}{H^4}. \notag
\end{align}
    \item For $j\notin I_B$, by Eq.\eqref{eq:cts-no-de-no-binding-second-expectation-big-0}, we have that when $\mathcal{E}_s\in\mathcal{F}_s$, it holds that
    \begin{align}
        \ell_{j,U_s}- C_7\sqrt{H \ln H} + \sum_{k=U_s}^{U_{s+1}-1}\mathbb{E}\left[b_{jk} - a_{jk}\cdot\mathbb{I}\left( r_k > \left\langle \bm{a}_k , \bq_s\right\rangle\right)\middle|\mathcal{F}_s \right]\geq \frac{\delta_B\cdot(U_{s+1}-U_s)}{2}. \notag
    \end{align}
    By Hoeffding's inequality, we have that when $\mathcal{E}_s\in\mathcal{F}_s$, it holds that
    \begin{align}
        &~~~\Pr\left(  \ell_{j,U_s}- C_7\sqrt{H \ln H} + \sum_{k=U_s}^{U_{s+1}-1}\left[b_{jk} - a_{jk}\cdot\mathbb{I}\left( r_k > \left\langle \bm{a}_k , \bq_s\right\rangle\right) \right] < 0 \middle|\mathcal{F}_s \right)\notag \\
        &\leq 2\exp\left( -\frac{\delta_B^2\cdot (U_{s+1}-U_s)}{2(\bar{b}+2\bar{a})^2} \right)\leq 2\exp\left(-\frac{\delta_B^2\cdot 3^{N_U-s}}{(\bar{b}+2\bar{a})^2}\right)\leq \frac{2}{H^5}, \notag
    \end{align}
    where the last inequality is due to {\color{black}$3^{N_U-s}\geq \frac{5}{\delta_B^2}\cdot(\bar{b}+2\bar{a})^2\ln H$}. By Eq.~\eqref{eq:cts-no-de-second-no-out-of-stock} and Eq.~\eqref{eq:cts-no-de-second-cntrol-stock-binding-pre}, we have that when $\mathcal{E}_s\in\mathcal{F}_s$, it holds that
\begin{align}
    \Pr\left(\forall j\notin I_B, \ell_{j,U_{s+1}} \geq C_7\sqrt{H \ln H}  \middle|\mathcal{F}_s\right)\geq 1-\frac{4m}{H^4}.
\end{align}
\end{itemize}

Finally, for each batch index $v$ we define the event $\mathcal{E}_v^*$ to be
\[
\mathcal{E}_v^* \defeq \mathcal{E}_v \cap \{\text{the three items in Lemma~\ref{lem:cts-no-de-inventory-conversion} all hold for $1\leq s\leq v-1$}\}.
\]
Note that conditioned on filtration $\mathcal{F}_s$, the event $\mathcal{E}^*_s$ is deterministic. We similarly let $\mathcal{E}^*_s\in\mathcal{F}_s$ denote the event $\mathcal{E}^*_s$ occurring in the filtration $\mathcal{F}_s$. From the discussions above, we have that under the conditions in Lemma~\ref{lem:cts-no-de-inventory-conversion}, we have $\mathcal{E}^*_1\in\mathcal{F}_1$ and
$
    \Pr\left( \mathcal{E}_{s+1}^*\in \mathcal{F}_{s+1} \middle| \mathcal{E}_{s}^*\in \mathcal{F}_{s}  \right)  \geq  1-\frac{10m}{H^4}$ holds for all $1\leq s <  N_U^\flat$.
Thus, we have
\begin{align}
    \Pr\left( \mathcal{E}_{s+1}^*\in \mathcal{F}_{s+1} \right)  &\geq  \Pr\left( \mathcal{E}_{s+1}^*\in \mathcal{F}_{s+1}  \middle|\mathcal{E}_{s}^*\in \mathcal{F}_{s} \right) \cdot\left(1-\Pr\left(\mathcal{E}_{s}^*\notin \mathcal{F}_{s} \right) \right)\notag\\
    &\qquad\qquad\qquad \geq 1-\frac{10m}{H^4}-\Pr\left(\mathcal{E}_{s}^*\notin \mathcal{F}_{s} \right)=\Pr\left(\mathcal{E}_{s}^*\in \mathcal{F}_{s} \right) -\frac{10m}{H^4}, \notag
\end{align}
which indicates that $
    \Pr\left( \mathcal{E}_{ N_U^\flat}^*\in \mathcal{F}_{ N_U^\flat} \right) \geq 1-\frac{10m}{H^3}$.
\hfill\Halmos

\subsection{Proof of Lemma~\ref{lem:cts-no-de-many-inv-for-no-binding}}
\label{sec:pf-cts-no-de-many-inv-for-no-binding}
Consider an index $j^*\notin I_B$, let $B'_j= B_j-(\delta_B-\epsilon)\cdot\mathbb{I}(j=j^*)$ for each $j$, and let $\bp'=\bp^*(\bm{B}')$. We have $I_B(\bm{B}')=I_B(\bm{B})$. By Lemma~\ref{lem:cts-no-de-d-robust}, we have $\bp'=\bp^*$, which indicates that
$
    B_j= (\delta_B-\epsilon) + B'_j \geq (\delta_B-\epsilon) +\mathbb{E}_{(r,\bm{a},\bm{b})\sim\mathcal{P}}\left[ a_{j}\cdot\mathbb{I}\left( r > \left\langle \bm{a} , \bp'\right\rangle \right) \right]=(\delta_B-\epsilon) +\mathbb{E}_{(r,\bm{a},\bm{b})\sim\mathcal{P}}\left[ a_{j}\cdot\mathbb{I}\left(r>\left\langle \bm{a} , \bp^*\right\rangle\right) \right]$. 
\hfill\Halmos

\subsection{Proof of Lemma~\ref{lem:cts-no-de-equivalence-no-out-of-stock}}
\label{sec:pf-cts-no-de-equivalence-no-out-of-stock}
We prove the lemma by induction on $v$. When $v=\kappa+1$, if $r_v>\left\langle \bm{a}_v , \bp^v\right\rangle$, we have $\sum_{k=1}^v b_{jk}\geq a_{jv}$ according to Eq.~\eqref{eq:cts-no-de-no-out-of-stock-basic}, and we can derive that $x_v=\mathbb{I}\left(r_v>\left\langle \bm{a}_v , \bp^v\right\rangle \right)$ for $v=\kappa+1$. 

Now assume that we have $x_v=\mathbb{I}(r_v>\left\langle \bm{a}_v , \bp^v\right\rangle)$ for all $\kappa<v\leq V < H$. Consider the case when $v=V+1$. By Eq.~\eqref{eq:cts-no-de-no-out-of-stock-basic} and the induction hypothesis, we have that for all $j\in[m]$,
\begin{align}
    0\leq \sum_{k=1}^{v-1} b_{jk} - \sum_{k=\kappa+1}^{v-1} a_{jk} x_k +b_{jv}- a_{jv}\cdot\mathbb{I}(r_v>\left\langle \bm{a}_v , \bp^v\right\rangle)= \ell_{j,v} + b_{jv}- a_{jv}\cdot\mathbb{I}(r_v>\left\langle \bm{a}_v , \bp^v\right\rangle). \notag
\end{align}
Thus, we have $x_v=\mathbb{I}(r_v>\left\langle \bm{a}_v , \bp^v\right\rangle)$.  \hfill\Halmos

\subsection{Proof of Lemma~\ref{lem:cts-no-de-saving-period-each-sub-inventory}}
\label{sec:pf-cts-no-de-saving-period-each-sub-inventory}
Let $\mathcal{F}_w=\{(r_k,\bm{a}_k,\bm{b}_k)\}_{k<V_w}$. For a fixed $w$, when conditioned on $\mathcal{F}_w$, the random variable $\bp^k$ is deterministic and takes the same value for all $k\in [V_w, V_{w+1})$. For convenience, we denote it by $\bu_w$ in this proof. Consider another random variable 
\begin{align}
    \bu_w^*\in\arg\min_{\bp\geq \bm{0}} f\left(\bp; \frac{\sum_{k=1}^{V_w-1} \bm{b}_k }{V_w-1}-C_2\sqrt{\frac{\ln H}{V_{w+1}-V_w}}\cdot\bm{1}\right). \notag
\end{align}
Let the event $\mathcal{E}_w^{\text{sub}}$ be
\begin{align}
    \mathcal{E}_w^{\text{sub}} \defeq \left\{\bu_w\in\Omega_p,\left\|\bu_w-\bu_w^*\right\|_2^2\leq (C_5)^2\cdot\frac{\ln H}{V_w-1}  \right\}  \bigcap  \mathcal{E}_{w}^{\text{resource},1}  \bigcap  \mathcal{E}_{V_w-1}^{\text{resource},2}, \notag
\end{align}
where $\mathcal{E}_{w}^{\text{resource},1},\mathcal{E}_{w}^{\text{resource},2}$ are defined in Eq.~\eqref{eq:cts-no-de-resource-1-def} and Eq.~\eqref{eq:cts-no-de-resource-2-def}. By Corollary~\ref{cor:cts-no-de-dual-conver-simple}, Eq.~\eqref{eq:cts-no-de-saving-period-1}, Eq.~\eqref{eq:cts-no-de-saving-period-2} and the discussion in the first part of Section~\ref{sec:pf-cts-no-de-first-t-2-performance}, we derive that 
\begin{align}
    \Pr\left( \mathcal{E}_w^{\text{sub}}\right)\geq 1-\frac{9m+2}{H^5}. \notag
\end{align}
We use the notation $\mathcal{E}_w^{\text{sub}}\in \mathcal{F}_w$ to denote the event $\mathcal{E}_w^{\text{sub}}$ occurring in the filtration $F_w$. We observe that, conditioned on $\mathcal{F}_w$, $\left\{b_{jk}-a_{jk}\cdot\mathbb{I}\left(r_k>\left\langle \bm{a}_k , \bp^k\right\rangle\right)\right\}_{k=V_w}^{V_{w+1}-1}$ is independent and identically distributed. Moreover, by Lemma~\ref{lem:cts-non-de-bounded-opt}, under the event $\mathcal{E}_w^{\text{sub}}$, it holds that $ \bu_w, \bu_w^*\in\Omega_p$. We have
\begin{align}
    &\left|\mathbb{E}\left[ b_{jk}-a_{jk}\cdot\mathbb{I}\left(r_k>\left\langle \bm{a}_k , \bp^k\right\rangle\right)\middle| \mathcal{F}_w \right]-B_j+\mathbb{E}\left[ a_{jk}\cdot\mathbb{I}\left(r_k> \left\langle \bm{a}_k , \bu_w^*\right\rangle\right) \middle| \mathcal{F}_w \right] \right| \notag \\
    &\qquad = \left|\mathbb{E}\left[ a_{jk}\cdot\mathbb{I}\left(r_k>\left\langle \bm{a}_k , \bu_w^*\right\rangle \right) \middle| \mathcal{F}_w  \right]-\mathbb{E}\left[ a_{jk}\cdot\mathbb{I}\left(r_k>\left\langle \bm{a}_k , \bu_w\right\rangle\right) \middle| \mathcal{F}_w  \right] \right|,\notag \\
    &\qquad \leq \bar{a} \left| \Pr\left( r_k > \left\langle \bm{a}_k , \bu_w\right\rangle  \middle| \mathcal{F}_w  \right) - \Pr\left( r_k > \left\langle \bm{a}_k , \bu_w^*\right\rangle  \middle| \mathcal{F}_w \right)  \right|\notag\\
    &\qquad \leq \mu\bar{a} \cdot \mathbb{E}\left[\left| \left\langle \bm{a}_k ,  \bu_w^* \right\rangle - \left\langle \bm{a}_k  \bu_w\right\rangle \right| \middle| \mathcal{F}_w  \right] \leq \mu\bar{a}^2 \left\| \bu_w - \bu_w^* \right\|_2 \leq C_5\mu\bar{a}^2\cdot \sqrt{\frac{\ln H}{V_w-1}}, \label{eq:cts-no-de-cumulate-inventory-sub-1}
\end{align}
where the second inequality is due to the non-degeneracy of the distribution $\mathcal{P}$. In the following, we prove the four inequalities in the lemma separately. The lemma is then proved by collecting failure probabilities.

\medskip
\noindent{\underline{The first inequality}.} By {\color{black} $C_0\ln H\geq 4C_1$}, we have that 
\begin{align}
    2(V_w-1)\geq V_{w+1}-V_w.\label{eq:cts-no-de-saving-sub-interval-est}
\end{align}
By the third part of Definition~\ref{def:cts-non-degeneracy}, under the event $\mathcal{E}_w^{\text{sub}}$, we have 
\begin{align}
    &\mathbb{E}\left[ b_{jk} - a_{jk}\cdot\mathbb{I}\left( r_k > \left\langle \bm{a}_k , \bu_w \right\rangle \right) \middle| \mathcal{F}_w  \right]\notag \\
    &\qquad \geq \left(B_j- \frac{\sum_{k=1}^{V_w-1}b_{jk}}{V_w-1} + C_2\sqrt{\frac{\ln H}{V_{w+1}-V_w}}\right) \notag \\
    &\qquad\qquad\qquad+ \frac{\sum_{k=1}^{V_w-1}b_{jk}}{V_{w}-1} -C_2\sqrt{\frac{\ln H}{V_{w+1}-V_w}} -\mathbb{E}\left[ a_{jk}\cdot\mathbb{I}\left( r_k > \left\langle \bm{a}_k , \bu_w^* \right\rangle\right) \right] - C_5\mu\bar{a}^2\cdot\sqrt{\frac{\ln H}{V_w-1}} \notag \\
    &\qquad \geq \left(C_2-2C_5\mu\bar{a}^2\right)\cdot\sqrt{\frac{\ln H}{V_{w+1}-V_w}} + B_j-\frac{\sum_{k=1}^{V_w-1}b_{jk}}{V_w-1} \notag\\
    &\qquad \geq \left(C_2-2C_5\mu\bar{a}^2 - 6\bar{b} \right)\cdot\sqrt{\frac{\ln H}{V_{w+1}-V_w}}\geq \frac{C_2}{2}\sqrt{\frac{\ln H}{V_{w+1}-V_w}}, \label{eq:cts-no-de-saving-each-expectation-lower}
\end{align}
where the first inequality is due to Eq.~\eqref{eq:cts-no-de-cumulate-inventory-sub-1}, the second one is due to the third part of Definition~\ref{def:cts-non-degeneracy} and Eq.~\eqref{eq:cts-no-de-saving-sub-interval-est}, the third one is due to the event $ \mathcal{E}_{V_w-1}^{\text{resource},2}$ and Eq.~\eqref{eq:cts-no-de-saving-sub-interval-est}, and the last one is due to that {\color{black} $C_2\geq 4C_5\mu\bar{a}^2 + 12\bar{b}$}. Using Hoeffding's inequality, we have
\begin{align}
    \Pr\left(  \sum_{k=V_w}^{V_{w+1}-1} b_{jk}-a_{jk}\cdot\mathbb{I}\left( r_k> \left\langle \bm{a}_k , \bu_w \right\rangle \right) < \frac{C_2}{4} \sqrt{{\ln H}\cdot({V_{w+1}-V_w }) } \middle|~\mathcal{F}_w  \right)\leq 2\exp\left( -\frac{(C_2)^2\ln H}{8\left(\bar{b}+2\bar{a}\right)^2} \right)\leq \frac{2}{H^5}, \notag
\end{align}
where the last inequality is due to that {\color{black} $C_2\geq 7(\bar{b}+2\bar{a})$}.
Thus, we have 
\begin{align}
    &\Pr\left(  \sum_{k=V_w}^{V_{w+1}-1} b_{jk}-a_{jk}\cdot\mathbb{I}\left( r_k> \left\langle \bm{a}_k , \bu_w \right\rangle \right) < \frac{C_2}{4} \sqrt{{\ln H}\cdot({V_{w+1}-V_w }) }   \right) \notag \\
    &\qquad\qquad\qquad\qquad\qquad\qquad\qquad\qquad\leq \sum_{\mathcal{E}_w^{\text{sub}}\in\mathcal{F}_w} \Pr(\mathcal{F}_w)\cdot\frac{2}{H^5} + 1- \Pr(\mathcal{E}_w^{\text{sub}})\leq \frac{9m+4}{H^5}.
    \notag
\end{align}
We conclude that 
\begin{align}
    &\Pr\left( \forall w\in[N_V],~ \sum_{k=V_w}^{V_{w+1}-1} b_{jk}-a_{jk}\cdot\mathbb{I}\left( r_k>\left\langle \bm{a}_k , \bu_w \right\rangle \right) \geq 4C_7 \sqrt{{\ln H}\cdot({V_{w+1}-V_w }) }   \right) \notag \\
    &\geq \Pr\left(\forall w\in[N_V],~  \sum_{k=V_w}^{V_{w+1}-1} b_{jk}-a_{jk}\cdot\mathbb{I}\left( r_k> \left\langle \bm{a}_k , \bu_w \right\rangle \right) \geq  \frac{C_2}{4} \sqrt{{\ln H}\cdot({V_{w+1}-V_w }) }   \right)\geq 1-\frac{9m+4}{H^4},\notag
\end{align}
where the first inequality is due to {\color{black} $C_2\geq 16 C_7$}.

\medskip
\noindent{\underline{The second inequality}.} 
By the third part of Definition~\ref{def:cts-non-degeneracy}, under the event $\mathcal{E}_w^{\text{sub}}$, for $j\in I_B$, we have 
\begin{align}
    &\mathbb{E}\left[ b_{jk} - a_{jk}\cdot\mathbb{I}\left( r_k > \left\langle \bm{a}_k , \bu_w\right\rangle \right) \middle| \mathcal{F}_w  \right],\notag \\
    &\qquad \leq \left(B_j- \frac{\sum_{k=1}^{V_w-1}b_{jk}}{V_w-1} + C_2\sqrt{\frac{\ln H}{V_{w+1}-V_w}}\right) \notag \\
    &\qquad\qquad\qquad+ \frac{\sum_{k=1}^{V_w-1}b_{jk}}{V_{w}-1} -C_2\sqrt{\frac{\ln H}{V_{w+1}-V_w}} -\mathbb{E}\left[ a_{jk}\cdot\mathbb{I}\left( r_k > \left\langle \bm{a}_k , \bu_w^* \right\rangle \right) \right] + C_5\mu\bar{a}^2\cdot\sqrt{\frac{\ln H}{V_w-1}} \notag \\
    &\qquad \leq \left(C_2+ 2C_5\mu\bar{a}^2\right)\cdot\sqrt{\frac{\ln H}{V_{w+1}-V_w}} + B_j-\frac{\sum_{k=1}^{V_w-1}b_{jk}}{V_w-1} \notag\\
    &\qquad \leq \left(C_2+2C_5\mu\bar{a}^2 + 6\bar{b} \right)\cdot\sqrt{\frac{\ln H}{V_{w+1}-V_w}}\leq 2C_2\sqrt{\frac{\ln T}{V_{w+1}-V_w}}, \notag
\end{align}
where the first inequality is due to Eq.~\eqref{eq:cts-no-de-cumulate-inventory-sub-1}, the second one is due to the third part of Definition~\ref{def:cts-non-degeneracy} (with $j\in I_B$) and Eq.~\eqref{eq:cts-no-de-saving-sub-interval-est}, the third one is due to the event $ \mathcal{E}_{V_w-1}^{\text{resource},2}$ and Eq.~\eqref{eq:cts-no-de-saving-sub-interval-est}, and the last one is due to that {\color{black} $C_2\geq 4C_5\mu\bar{a}^2 + 12\bar{b}$}. Using Hoeffding's inequality, we have
\begin{align}
    \Pr\left(  \sum_{k=V_w}^{V_{w+1}-1} b_{jk}-a_{jk}\cdot\mathbb{I}\left( r_k> \left\langle \bm{a}_k , \bu_w\right\rangle \right) \geq 3C_2 \sqrt{{\ln H}\cdot({V_{w+1}-V_w }) } \middle|~\mathcal{F}_w  \right)\leq 2\exp\left( -\frac{2(C_2)^2\ln H}{\left(\bar{b}+2\bar{a}\right)^2} \right)\leq \frac{2}{H^5}, \notag
\end{align}
where the last inequality is due to that {\color{black} $C_2\geq 2(\bar{b}+2\bar{a})$}.
Thus, we have 
\begin{align}
    &\Pr\left(  \sum_{k=V_w}^{V_{w+1}-1} b_{jk}-a_{jk}\cdot\mathbb{I}\left( r_k>\left\langle \bm{a}_k , \bu_w\right\rangle \right) \geq 3C_2 \sqrt{{\ln H}\cdot({V_{w+1}-V_w }) }   \right) \notag \\
    &\qquad\qquad\qquad\qquad\qquad\qquad\qquad\qquad\leq \sum_{\mathcal{E}_w^{\text{sub}}\in\mathcal{F}_w} \Pr(\mathcal{F}_w)\cdot\frac{2}{H^5} + 1- \Pr(\mathcal{E}_w^{\text{sub}})\leq \frac{9m+4}{H^5}. \notag
\end{align}
We conclude that 
\begin{align}
    &\Pr\left( \forall w\in[N_V],~ \sum_{k=V_w}^{V_{w+1}-1} b_{jk}-a_{jk}\cdot\mathbb{I}\left( r_k> \left\langle \bm{a}_k , \bu_w \right\rangle \right) \leq \frac{C_8}{6} \sqrt{{\ln H}\cdot({V_{w+1}-V_w }) }   \right) \notag \\
    &\geq \Pr\left(\forall w\in[N_V],~  \sum_{k=V_w}^{V_{w+1}-1} b_{jk}-a_{jk}\cdot\mathbb{I}\left( r_k>\left\langle \bm{a}_k , \bu_w\right\rangle \right) \leq  3 C_2 \sqrt{{\ln H}\cdot({V_{w+1}-V_w }) }   \right)\geq 1-\frac{9m+4}{H^4}, \notag
\end{align}
where the first inequality is due to {\color{black} $C_8\geq 18 C_2$}.

\medskip
\noindent{\underline{The third inequality}.} By Eq.~\eqref{eq:cts-no-de-saving-each-expectation-lower} and Hoeffding's inequality, we have 
\begin{align}
        \Pr\left(  \sum_{k=V_w}^{t} b_{jk}-a_{jk}\cdot\mathbb{I}\left( r_k>\left\langle \bm{a}_k , \bu_w \right\rangle \right) < - C_7 \sqrt{{\ln H}\cdot({V_{w+1}-V_w }) } \middle|~\mathcal{F}_w  \right)\leq 2\exp\left( -\frac{2(C_7)^2\ln H}{\left(\bar{b}+2\bar{a}\right)^2} \right)\leq \frac{2}{H^5}, \notag
\end{align}
where the second inequality is due to {\color{black} $C_7\geq 2(\bar{b}+2\bar{a})$}. Similarly, we can derive that 
\begin{align}
     &\Pr\left(  \sum_{k=V_w}^{t} b_{jk}-a_{jk}\cdot\mathbb{I}\left( r_k>\left\langle \bm{a}_k , \bu_w \right\rangle \right) < - C_7 \sqrt{{\ln H}\cdot({V_{w+1}-V_w }) } \right)\notag \\
     &\qquad\qquad\qquad\qquad\qquad\qquad\qquad\leq \sum_{\mathcal{E}_w^{\text{sub}}\in\mathcal{F}_w} \Pr(\mathcal{F}_w)\cdot\frac{2}{H^5} + 1- \Pr(\mathcal{E}_w^{\text{sub}})\leq \frac{9m+4}{H^5}, \notag
\end{align}
which indicates that 
\begin{align}
    \Pr\left( \forall w\in[N_V]\text{ and } V_w\leq t<V_{w+1},~ \sum_{k=V_w}^{t} b_{jk}-a_{jk}\cdot\mathbb{I}\left( r_k>\left\langle \bm{a}_k , \bu_w \right\rangle \right) \geq - C_7 \sqrt{{\ln H}\cdot({V_{w+1}-V_w }) } \right)\notag\\
    \geq 1-\frac{9m+4}{H^4}. \notag
\end{align}

\medskip
\noindent{\underline{The fourth inequality}.} By Hoeffding's inequality, we have
\begin{align}
     &\Pr\left(\sum_{k=1}^{V_1-1} b_{jk} \geq  C_7\sqrt{(V_2-V_1)\cdot\ln H} \right) \geq \Pr\left(\sum_{k=1}^{V_1-1} b_{jk} \geq \frac{\underline{b}}{2}(V_1-1)\right)  \notag\\
     &\qquad\qquad\qquad\qquad\qquad\qquad\qquad\geq 1- 2\exp\left( -\frac{\underline{b}^2(V_1-1)}{2\bar{b}^2}\right)\geq 1-2\exp\left(-\frac{\underline{b}^2 C_0 \ln^2 H}{2\bar{b}^2}\right)\leq 1-\frac{2}{H^5}, \notag
\end{align}
where the first inequality is due to {\color{black} $\underline{b}C_0\ln H\geq 2\sqrt{C_1}\cdot C_7$}, and the fourth is due to {\color{black} $C_0\geq \frac{10\bar{b}^2}{\underline{b}^2}$}. 
\hfill\Halmos

\subsection{Proof of Lemma~\ref{lem:cts-no-de-test-index-1}}
\label{sec:pf-cts-no-de-test-index-1}
For the first component of the event $\mathcal{E}^{\text{test},1}$, by Hoeffding's inequality, we have
\begin{align}
     \Pr\left(\left|\frac{\sum_{k=1}^{M} b_{jk}}{M} -B_j\right|>\sqrt{\frac{5\bar{b}^2\cdot\ln H}{M}}\right)\leq 2\exp\left(-\frac{5\bar{b}^2\ln H }{\bar{b}^2}\right)\leq\frac{2}{H^5}. \notag
\end{align}
By  {\color{black} $\delta_B^2\cdot\lfloor\frac{H}{2}\rfloor\geq 5\bar{b}^2\cdot\ln H$}, we have 
\begin{align}
    \Pr\left(\left|\frac{\sum_{k=1}^{M} b_{jk}}{M} -B_j\right|>\min\left(\delta_B,\sqrt{\frac{5\bar{b}^2\cdot\ln H}{M}}\right)\right)=\Pr\left(\left|\frac{\sum_{k=1}^{M} b_{jk}}{M} -B_j\right|>\sqrt{\frac{5\bar{b}^2\cdot\ln H}{M}}\right)\leq\frac{2}{H^5}. \notag
\end{align}
The failure probability of the second component of the event $\mathcal{E}^{\text{test},1}$ be directly derived from Corollary~\ref{cor:cts-no-de-dual-conver-simple}, using {\color{black}$\lfloor\frac{H}{2}\rfloor\geq\frac{20\bar{a}}{\lambda_{\text{min}}}\cdot\ln H +10$}. Finally, the lemma is proved by collecting the failure probabilities of both components of the event $\mathcal{E}^{\text{test},1}$. \hfill\Halmos

\section{Additional Experimental Results}
\label{sec:add-experiment}
This appendix contains additional figures referenced in the main text. The experimental settings and discussions are provided in Section~\ref{sec:experiment}.
\begin{figure}[htbp]
    \centering
    \begin{subfigure}{\textwidth}
        \centering
        \includegraphics[width=0.32\textwidth]{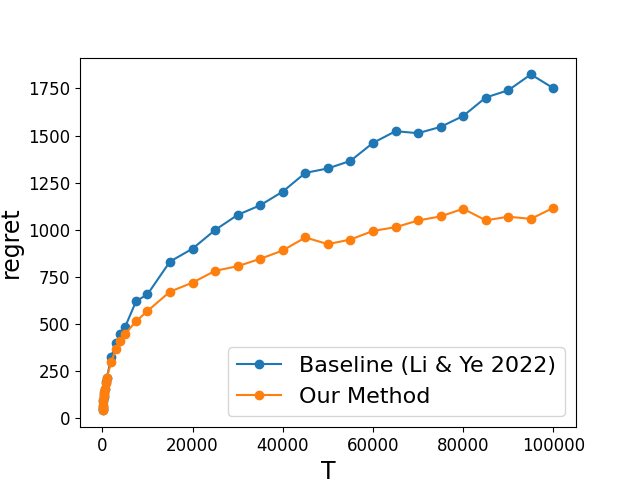}
        \includegraphics[width=0.32\textwidth]{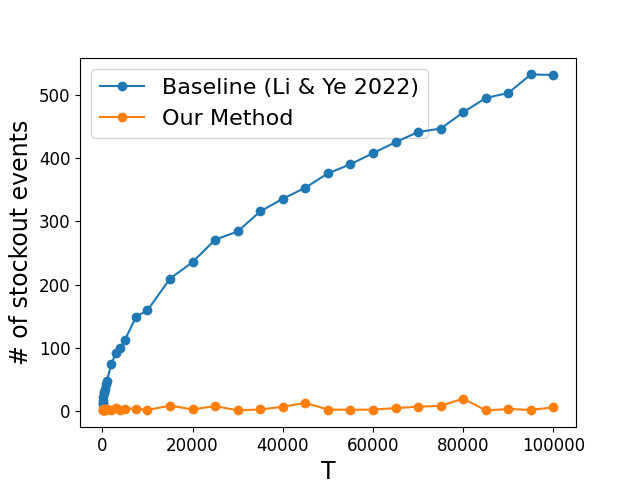}
        \includegraphics[width=0.32\textwidth]{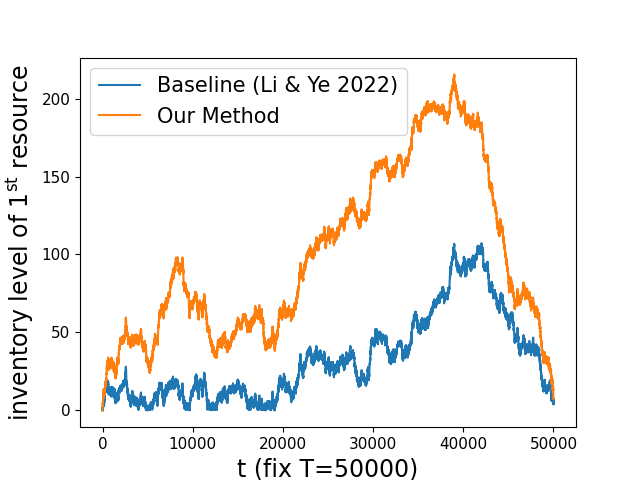}
        \caption{Random Input I with $m=10$}
    \end{subfigure}

    \vspace{0.3em}

    \begin{subfigure}{\textwidth}
        \centering
        \includegraphics[width=0.32\textwidth]{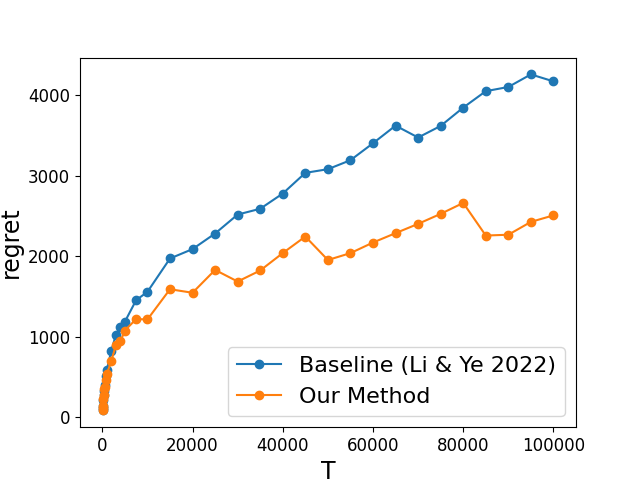}
        \includegraphics[width=0.32\textwidth]{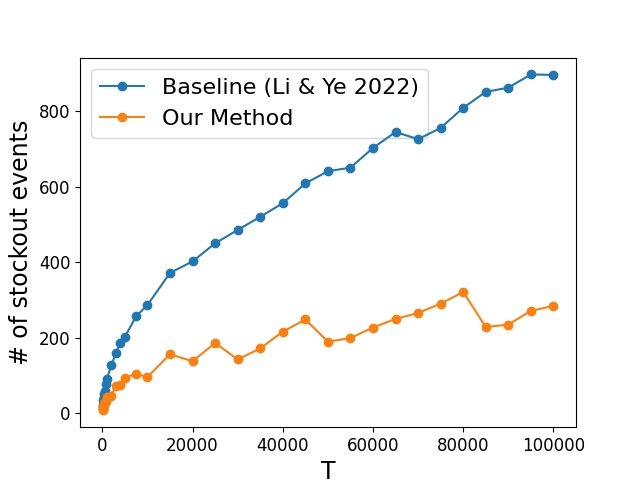}
        \includegraphics[width=0.32\textwidth]{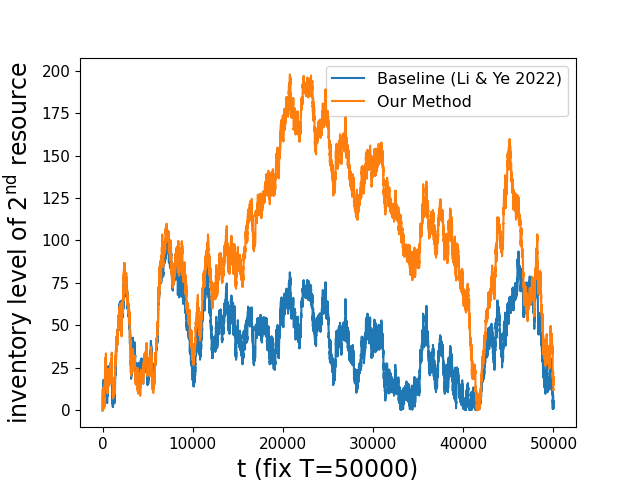}
        \caption{Random Input II with $m=10$}
    \end{subfigure}

    \caption{Additional empirical comparison between our algorithm and \cite{li2022online}}
    \label{fig:random-add}
\end{figure}

\end{document}